\documentclass[a4paper,fleqn]{cas-sc}

\usepackage[numbers]{natbib}
\usepackage{graphicx}
\usepackage{color,soul} 
\usepackage{array}
\usepackage{amsmath}
\usepackage{amsfonts}
\usepackage{amssymb}
\usepackage{psfrag}
\usepackage{epstopdf}
\usepackage{subcaption}
\usepackage{booktabs}
\usepackage{latexsym}
\usepackage{kantlipsum}
\usepackage{balance}
\usepackage{lscape}
\usepackage{rotating}
\usepackage{multicol}
\usepackage{multirow,bigdelim}
\usepackage{color, colortbl}
\usepackage{threeparttable}
\usepackage{float}

\newcommand{\bcolon}{\boldsymbol{:}}
\newcommand{\eval}[2][\right]{\relax\ifx#1\right\relax \left.\fi#2#1\rvert}

\newcommand{\bn}{\boldsymbol{n}}
\newcommand{\bt}{\boldsymbol{t}}
\newcommand{\bu}{\boldsymbol{u}}

\newcommand{\bB}{\boldsymbol{B}}
\newcommand{\bC}{\boldsymbol{C}}

\newcommand{\bG}{\boldsymbol{G}}

\newcommand{\bN}{\boldsymbol{N}}

\newcommand{\bR}{\boldsymbol{R}}
\newcommand{\bS}{\boldsymbol{S}}
\newcommand{\bT}{\boldsymbol{T}}

\newcommand{\bV}{\boldsymbol{V}}

\newcommand{\del}{\boldsymbol{\nabla}}
\newcommand{\bcross}{\boldsymbol{\times}}

\newcommand{\btau}{\boldsymbol{\tau}}

\makeatletter
\newcommand{\Rmnum}[1]{\expandafter\@slowromancap\romannumeral #1@}
\makeatother

\newcommand{\bg}{\boldsymbol{g}}

\newcommand{\bK}{\boldsymbol{K}}
\newcommand{\bP}{\boldsymbol{P}}

\newcommand{\bbeta}{\boldsymbol{\beta}}
\newcommand{\betahat}{\hat{\bbeta}}

\newcommand{\domega}{\,d\varOmega}
\newcommand{\dgamma}{\,d\varGamma}

\newcommand{\bH}{\boldsymbol{H}}

\newcommand{\tbar}{\bar{\bt}}

\newcommand{\utilde}{\tilde{\bu}}

\newcommand{\beps}{\boldsymbol{\epsilon}}

\newcommand{\bXi}{\boldsymbol{\Xi}}

\newcommand{\bepsilon}{\boldsymbol{\epsilon}}
\def\tsc#1{\csdef{#1}{\textsc{\lowercase{#1}}\xspace}}
\tsc{WGM}
\tsc{QE}
\tsc{EP}
\tsc{PMS}
\tsc{BEC}
\tsc{DE}
\begin{document}
\let\WriteBookmarks\relax
\def\floatpagepagefraction{1}
\def\textpagefraction{.001}
\shorttitle{A two-field hybrid IGA formulation to alleviate locking } 
\shortauthors{D.S.Bombarde et~al.}

\title [mode = title]{A novel hybrid isogeometric element based on two-field Hellinger-Reissner principle to alleviate different types of locking}                    

\tnotemark[1]

\tnotetext[1]{Supported by Science \& Engineering Research Board (SERB), and Department of Science \& Technology (DST), Government of India, under the project IMP/2019/000276.}

% ---------------- author 1 ----------------------
\author[]{Dhiraj S. Bombarde}
\cormark[1] % 
\ead{dhira176103023@iitg.ac.in}
% --------------- author 2 -----------------------
\author[]{Sachin S. Gautam}
\ead{ssg@iitg.ac.in}
% --------------- author 3 -----------------------
\author[]{Arup Nandy}
\ead{arupn@iitg.ac.in}
% ------------------------------------------------

\address[]{Department of Mechanical Engineering, Indian Institute of Technology Guwahati, Guwahati 781039, Assam, India}

\cortext[cor1]{Corresponding author}

\begin{abstract}
In the present work, a novel class of hybrid elements is proposed to alleviate the locking anomaly in non-uniform rational B-spline (NURBS)-based isogeometric analysis (IGA) using a two-field Hellinger-Reissner variational principle. The proposed hybrid elements are derived by adopting the independent interpolation schemes for displacement and stress field. The key highlight of the present study is the choice and evaluation of higher-order terms for the stress interpolation function to provide a locking-free solution. Furthermore, the present study demonstrates the efficacy of the proposed elements with the treatment of several two-dimensional linear-elastic benchmark problems alongside the conventional single-field IGA, Lagrangian-based finite element analysis (FEA), and hybrid FEA formulation. It is shown that the proposed class of hybrid elements performs effectively for analyzing the nearly incompressible problem domains that are severely affected by volumetric locking along with the thin plate and shell problems where the shear and membrane locking is dominant. A better coarse mesh accuracy of the proposed method in comparison with the conventional formulation is demonstrated through various numerical examples. Moreover, the formulation is not restricted to the locking-dominated problem domains but can also be implemented to solve the problems of general form without any special treatment. Thus, the proposed method is robust, most efficient, and highly effective against different types of locking.
\end{abstract}


\maketitle

% ==============================================================================
% Introduction
% ==============================================================================
\section{Introduction}

The finite element analysis (FEA) is a widely practiced numerical procedure to solve the partial differential equations governing a mathematical model for a physical problem. Over a period of time, FEA successfully distinguished itself in various problem domains and found numerous applications in a diverse set of engineering fields. Despite its widespread applications, traditional FEA has certain drawbacks. One of which is the geometry approximation. Most often, the physical domain of the problem is modeled using the computer-aided design (CAD) geometries, which are further treated to create the FE mesh. However, the resultant FE mesh is an approximate version of the actual physical CAD domain. Though refining the mesh considerably improves the approximation but with an additional computing cost and pre-processing efforts. It has been argued that these geometric irregularities can lead to significant errors in analysis. Moreover, the time spent on pre-processing, i.e.,\ to make an analysis-ready FE mesh model, is significantly high as compared to the actual analysis time \cite{Cottrell2009_book}. 

To overcome the stated limitations, Hughes et al.\ introduced the concept of IGA, which provides an efficient integration between CAD and FEA \cite{Hughes2005}. In IGA, the physical geometry remains invariant regardless of the type or number of elements. Furthermore, once the initial coarse mesh is created, IGA simplifies the mesh refinement by removing further dependency on CAD. The fundamental idea behind IGA was to reinstate the conventional Lagrangian interpolation functions with NURBS basis functions that are extensively utilized in engineering design with the existence of numerous effective and numerically stable algorithms. Furthermore, retention of the essential mathematical properties like non-negativity of basis functions, linear independence, partition of unity, and variation diminishing property \cite{Piegl1997} of the NURBS assisted the foundation for IGA background.

Since its introduction in 2005, IGA has been widely practiced in different directions and proved to be a powerful method that out-perform FEA in most of the numerical aspects. The ability of the IGA framework in retaining the exact geometry with higher inter-element continuity, smooth representation of surfaces, and tight coupling between CAD geometries and FEA model lead to effective implementation in the several application domains such as contact formulations \cite{DeLorenzis2014, Agrawal2020}, structural shape optimization \cite{Wall2008}, fluid and fluid-structure analysis \cite{Gomez2010, Bazilevs2010, Hosseini2015}, structural vibration problem \cite{Cottrell2006}, shell and plate problems  \cite{Benson2013, Echter2013, Riffnaller-Schiefer2016}, and many more. Furthermore, the effective implementation influenced researchers to make efforts to integrate NURBS based IGA in existing commercial software like LS-DYNA \cite{Hartmann2011} and Abaqus \cite{Duval2012} and develop several dedicated tools e.g.,\ GeoPDEs \cite{DeFalco2011}, PetIGA \cite{Dalcin2016}, and $\Pi$gasus \cite{Ratnani2012}.

IGA turns out to be a better and generalized version of FEA in most of the aspects. However, there is a possibility that the limitations exist in conventional FEA may pass on to its generalized form. In fact, it is true. Even though IGA performs superior to FEA in several ways, it does exhibit similar limitations sooner or later like its FE counterpart \cite{Echter2010}. One such phenomenon is locking, which appears during the analysis of thin structural geometries and incompressible material behavior. 

The term \textit{locking} is used to label particular conditions when FE solutions of the nodal variables are underestimated, or the near-infinite stiffness promotes the absurd results for the solution. The most common types of locking are volumetric, shear, and membrane locking. The volumetric locking (also known as dilatation or Poisson's locking) is associated with the Poisson's ratio ($\nu$). When the material is incompressible ($\nu=0.5$) or nearly incompressible ($\nu \approx 0.5$), IGA scheme result in unrealistic solutions or low convergence rates for a practical range of discretization \cite{Babuska1992, Prathap1993}. Whereas the shear and membrane locking is associated with the thickness of the domain ($t$). For elements having one dimension significantly smaller than the other dimensions, FE as well as IGA solutions, can lead to additional stiffening effect and wildly oscillating shear forces along the length of the element. 

To overcome the stated limitation in context of IGA, the explored contributions are limited. The classical shell theories that are extensively used in the analysis of thin structures, to handle the shear or membrane locking, are successfully incorporated into the IGA framework. Popular theories and its IGA counterpart involve the Reissner-Mindlin shell theory \cite{Benson2010, BeiraodaVeiga2012}, Kirchhoff-Love theory \cite{Kiendl2009}, and blended shell theory \cite {Benson2013}. Furthermore, the contributions are extended to the degenerated shell approach, or solid-like shell formulation \cite{Benson2010, Hosseini2013} along with the NURBS-based solid shell elements \cite{Combescure2013}. Another popular technique to alleviate different types of locking is reduced and selective reduced integration, which has been explored in the framework of IGA \cite{Hughes2010} with promising results. Furthermore, the methods like $\bar{B}$ and $\bar{F}$ projection techniques are also investigated for handling the volumetric locking in case of incompressible or nearly incompressible problem domain \cite{Elguedj2008, Zhang2018}. Along with these, several multi-field variational techniques are found to be effective in alleviating different types of locking in IGA. Popular methods include assumed natural strains (ANS) \cite{ Caseiro2014, Caseiro2015} and enhanced assumed strains (EAS) \cite{Ruip2012, Tayor2010a}, which satisfactorily handles the situations where the standard IGA is prone to locking. 

It should be pointed out that most of the mixed formulations, modified to work with IGA, are strain-based approaches and less focus has been given to the stress-based formulations. This being the motivation, the present study emphasizes on developing a relatively simple but robust two-field stress-based IGA formulation that is capable of producing locking-free solutions irrespective of the type of locking. The proposed hybrid elements are derived based on a two-field mixed variational principle where displacement and stresses are the independent field variables. The key notion is the choice and evaluation of stress interpolation functions which is inspired from the work of Jog \cite{Jog2005, Jog2010} in context of conventional FEA framework. The efficient FE implementation of the stated approach for large-deformation contact mechanics \cite{Agrawal2019_hyb}, structural acoustics \cite{Jog2015_hyb}, electromagnetic analysis \cite{Jog2014_hyb}, analysis of electromechanical systems \cite{Agrawal2017_hyb}, and coupled fluid-structure problem \cite{A.RoychowdhuryA.Nandy2014_hyb} confirms the effectiveness and robustness of the method. In the present study, the systematic evaluation of stress interpolation functions specific to the NURBS interpolations and its effective implementation in a two-dimensional linear elasticity regime is investigated. The principle concept follows the assessment of normal stress components in relation to derivatives of the respective displacement interpolations and the choice of higher-order terms in shear components are defined such that they ensure correct stiffness rank and suppress the spurious zero-energy mode which is necessary to avoid locking. Though the formulation is computationally expensive, it compensates its efficiency by providing a relatively simple formulation and high coarse mesh accuracy. Elements developed by this theory perform quite well irrespective of the type of locking. The same element can be used in situations where problems demand incompressibility, or near incompressibility of material or to solve plate/shell geometries and so on. Additionally, the same formulation can be implemented for the standard problems (absence of locking effect) without affecting the solution accuracy, making it easier to implement in coupled problem domains. Another advantage of the hybrid formulation is that no energy or work principle or variational norms are violated, and hence mathematically robust formulation is obtained. The paper further highlights the concepts involved in hybrid IGA formulation and the successful implementation of the stated method on numerous benchmark examples.

The paper is divided into five main sections, including Introduction, and organized as follows: \textit{Section 2} is focused on the mathematical preliminaries followed by \textit{Section 3} which addresses the fundamental concepts involved in the proposed two-field hybrid stress formulation in IGA framework. Furthermore, the section deals with a step-wise implementation procedure followed by \textit{section 4} where several examples with promising results are illustrated. \textit{Section 5} summarizes the results and validates the performance of the proposed formulation.
% ------------------------------------------------------------------------------
% ==============================================================================
% Mathematical preliminaries
% ==============================================================================
\section{Mathematical preliminaries}
% --------------------------- Knot vector --------------------------------------
\subsection{Knot vector}
A set of non-decreasing parametric values over a specific direction is defined as a knot vector. A typical knot vector is represented as 
$\bXi = \{\xi_1,\xi_2,\ldots,\xi_{n+p+1}\}$ where $\xi_i \in {\rm I\!R}$ (one dimensional space) with $i$ is the knot index. $i$ varies from 1 to $n+p+1$, where $p$ and $n$ are degree of polynomial function and number of basis functions respectively. The knot vector decides the division of the parametric space into finite intervals which are known as knot spans. Knot vectors are categorized into two types, uniform and non-uniform knot vectors. If the knots are equally spaced in parametric space then the knot vector is said to be uniform otherwise it is called as non-uniform knot vector. Successive knots in a knot vector can have identical values which allow having multiple knot entries at the same coordinate in the parametric space. A number of repeating knots is referred as the multiplicity of a knot. If the first and last knots have $p+1$ multiplicity (i.e. repeating $p+1$ times) then the knot vector is termed as open knot vector. Multiplicity of knot also defines the continuity at that knot which is given as $C^{p-k}$ where $k$ is the number of repeated occurrence of a knot in a knot vector. For instance, an open knot vector in which the first and last knots are repeated $p+1$ times shows $C^{-1}$ continuity at the extremities.
% ------------------------ B-spline basis functions ----------------------------
\subsection{B-spline basis functions}
A basis function can be defined as a single curve element or piecewise polynomial which is used as a basis of linear combination to describe a specific curve. For a given knot vector $\bXi$, the B-spline basis functions are defined in a recursive form as follows,
% ------------------------ De-boor recursive formula ---------------------------
\begin{align} \label{DE_boor_formula}
\text{for } p = 0, N_{i,p}(\xi) &= \begin{cases}
    1 & \text{if } \xi_i \leq \xi < \xi_{i+1}\\
    0 & \text{Otherwise }
  \end{cases}\\
\text{for } p \geq 1, N_{i,p}(\xi) &= \frac{\xi - \xi_i}{\xi_{i+p} - \xi_i}N_{i,p-1} + \frac{\xi_{i+p+1} - \xi}{\xi_{i+p+1} - \xi_{i+1}}N_{i+1,p-1}
\end{align}
% ------------------------------------------------------------------------------
This is referred to as the Cox-de Boor recursion formula \cite{Piegl1997}.  During the evaluation of these functions, the ratios of the form $\frac{0}{0}$ or $\frac{1}{0}$ are considered to be equal to $0$. The B-spline basis functions for polynomial degree $0$ and $1$ will result in standard piecewise constant and Lagrangian interpolation functions respectively. However, higher degree basis functions ($p\geqslant 2$) differ from their FEA counterparts. B-splines basis functions also exhibit important mathematical properties like non-negativity of basis functions, linear independence, partition of unity, and variation diminishing property \cite{Piegl1997}. 

Multivariate B-spline basis function are used to obtain B-spline surfaces or solids. These functions are merely the tensor product of univariate basis functions. A bivariate B-spline basis function can be defined as, 
% ---------------------------- Bivariate Basis ----------------------------------
\begin{equation}
N_{i,j}^{p,q}(\xi,\eta) = N_{i,p}(\xi) M_{j,q}(\eta)
\end{equation}
% -------------------------------------------------------------------------------
whereas, a trivariate function is given as,
% ---------------------------- trivariate Basis ---------------------------------
\begin{equation}
N_{i,j,k}^{p,q,r}(\xi,\eta,\zeta) = N_{i,p}(\xi) M_{j,q}(\eta) L_{k,r}(\zeta)
\end{equation}
% -------------------------------------------------------------------------------
where $N, M$ and $L$ are $p, q $ and $r$\textsuperscript{th} degree basis function in $\xi, \eta$ and $\zeta$ direction.
% -------------------- B-spline curves, surfaces, and solids --------------------
\subsection{B-spline curves, surfaces, and solids}
A $p$\textsuperscript{th} degree B-spline curve for a given knot vector and a set of control points (CP) is defined as, 
% ----------------------------- B-spline curve ---------------------------------
\begin{equation}
\bC(\xi) = \sum_{i=1}^{n}N_{i,p}(\xi)\bP_i
\end{equation}
% ------------------------------------------------------------------------------
whereas, the B-spline surfaces are constructed by considering a bidirectional net of CP, two knot vectors and the tensor product of two univariate B-spline basis functions and given as,
% ----------------------------- B-spline surface -------------------------------
\begin{equation}
\bS(\xi,\eta) = \sum_{i=1}^{n} \sum_{j=1}^{m}N_{i,j}^{p,q}(\xi,\eta)\bP_{i,j}
\end{equation}
% ------------------------------------------------------------------------------
Similarly, the B-spline volumes are defined as,
% ------------------------------ B-spline solid --------------------------------
\begin{equation}
\bV(\xi,\eta,\zeta) = \sum_{i=1}^{n} \sum_{j=1}^{m} \sum_{k=1}^{l} N_{i,j,k}^{p,q,r}(\xi,\eta,\zeta)\bP_{i,j,k}
\end{equation}
where $n,m$ and $l$ are number of basis functions in $\xi,\eta$ and $\zeta$ direction respectively. $\bP_{i}$, $\bP_{i,j}$, and $\bP_{i,j,k}$ are control point co-ordinates $[x\quad y\quad z]^T$ for $(i)^{th}$, $(i,j)^{th}$, and $(i,j,k)^{th}$ control point. 
% ------------------------------------------------------------------------------
% ==============================================================================
%	Section- NURBS basis functions
% ==============================================================================
\subsection{NURBS basis functions}
The evolution of B-spline interpolations to Non-Uniform Rational B-Splines (NURBS) provides an advantage of representing a wide range of objects including conic sections like circle, sphere, cylinder, etc.\ NURBS are the generalization of B-splines basis functions where Non-Uniform refers to a non-uniform knot vector whereas Rational B-splines describe the rationalization of B-spline basis functions. It inherits all the mathematical properties of B-splines with the added advantage of weights which allows them to offer great flexibility and accuracy in generation of CAD geometries. A univariate NURBS basis function is defined as,
% ============================ univariate Basis NURBS ==========================
\begin{equation}
R_{i,p}(\xi) = \frac{w_i N_{i,p}(\xi)}{W(\xi)} \quad \text{where,}\hspace{0.1cm}
W(\xi) = \sum_{i=1}^{n_{cp}} w_i N_{i,p}(\xi)
\end{equation}
% ------------------------------------------------------------------------------
where $N_{i,p}(\xi)$ is a standard B-spline basis function, $w_i$ is strictly positive ($w_i > 0$) set of weights associated with CPs, and $n_{cp}$ denotes the total number of CPs. 

A bivariate NURBS basis function can be defined as, 
% ============================= bivariate NURBS basis ==========================
\begin{equation}\label{Bivariate_NURBS_basis}
R_{i,j}^{p,q}(\xi,\eta) = R_{i,p}(\xi)R_{j,q}(\eta) = \frac{N_{i,p}(\xi) M_{j,q}(\eta) w_{i,j}}{\sum_{i=1}^{n}\sum_{j=1}^{m} N_{i,p}(\xi) M_{j,q}(\eta) w_{i,j}}
\end{equation}
% ------------------------------------------------------------------------------
whereas, a trivariate function is given as,
% ============================= trivariate NURBS basis =========================
\begin{equation}
R_{i,j,k}^{p,q,r}(\xi,\eta,\zeta)  = R_{i,p}(\xi)R_{j,q}(\eta)R_{k,r}(\zeta)  = \frac{N_{i,p}(\xi) M_{j,q}(\eta) L_{k,r}(\zeta) w_{i,j,k}}{\sum_{i=1}^{n}\sum_{j=1}^{m}\sum_{k=1}^{l} N_{i,p}(\xi) M_{j,q}(\eta) L_{k,r}(\zeta) w_{i,j,k}}
\end{equation}
% ------------------------------------------------------------------------------
where $N, M$ and $L$ are $p, q$ and $r$\textsuperscript{th} degree B-spline basis function in $\xi, \eta$ and $\zeta$ direction respectively with $w_{i,j}$ and $w_{i,j,k}$ being the weights associated with CP net over the domain $\Omega(\xi,\eta)$ and $\Omega(\xi,\eta,\zeta)$ respectively.

\subsection{NURBS curves, surfaces, and solids }
A $p^{th}$ degree NURBS curve for a given knot vector and a set of CP is defined as,
% ================================= NURBS curve ================================
\begin{equation}
\bC(\xi) = \sum_{i=1}^{n}R_{i,p}(\xi)\bP_i
\end{equation}
% ------------------------------------------------------------------------------

The NURBS surfaces are constructed by considering a bidirectional net of CP, two knot vectors and the tensor product of two univariate NURBS basis functions and given as, 
% ================================ NURBS surface ===============================
\begin{equation}
\bS(\xi,\eta) = \sum_{i=1}^{n} \sum_{j=1}^{m}R_{i,j}^{p,q}(\xi,\eta)\bP_{i,j}
\end{equation}
% ------------------------------------------------------------------------------
and the NURBS volumes are defined as,
% ================================ NURBS solid =================================
\begin{equation}
\bV(\xi,\eta,\zeta) = \sum_{i=1}^{n} \sum_{j=1}^{m} \sum_{k=1}^{l} R_{i,j,k}^{p,q,r}(\xi,\eta,\zeta)\bP_{i,j,k}
\end{equation}
% ------------------------------------------------------------------------------
% ==============================================================================
% Fundamental concepts in two-field hybrid stress formulation in context of IGA
% ==============================================================================
\section{Fundamental concepts in two-field hybrid stress formulation in context of IGA}

\subsection{A classical two-dimensional linear elasticity problem}
To understand the underlying concepts in two-field stress formulation and the differentiating features from the conventional single-field IGA formulation \cite{Hughes2005}, it will be wise to appraise the classical two-dimensional linear elasticity problem. Let $\Omega$ be the open domain with boundary $\Gamma$ which is composed of two disjoint regions such that $\Gamma = \Gamma_u \cup \Gamma_t$ where $\Gamma_u$ is displacement boundary and $\Gamma_t$ is traction boundary. The governing equations for the linear elasticity problem is given as,
% =========================== governing equation ===============================
\begin{subequations}
\begin{align}
\del\cdot\btau + \boldsymbol{f} &= \bf{0} \hspace{0.25cm}\text{on}\hspace{0.25cm}\Omega, \label{differenctial_equation}\\
\btau \bn &= \bt \hspace{0.25cm}\text{on}\hspace{0.25cm}\Gamma, \label{cauchys_equation}\\
\bu &= \bu_0 \hspace{0.25cm}\text{on}\hspace{0.25cm}\Gamma_u, \label{displacement_boundary_condition}\\
\bt &= \tbar \hspace{0.25cm}\text{on}\hspace{0.25cm}\Gamma_t \label{traction_boundary_condition}
\end{align}
\end{subequations}
% -------------------------------------------------------------------------------
where $\btau$ is the Cauchy's stress tensor, $\bn$ is the unit outward normal to $\Gamma$, $\tbar$ is traction defined on the boundary $\Gamma_t$ and $\boldsymbol{f}$ is the body force vector. The small strain tensor $\beps$ is defined as,
% ========================== small strain u ====================================
\begin{equation} \label{strain_tensor_u}
\beps(\bu) = \frac{1}{2} [(\del \bu)+(\del\bu)^T] 
\end{equation}
% ------------------------------------------------------------------------------
and the stress-strain relation is given as $\btau= \mathcal{C}\bcolon \beps$ where $\mathcal{C}$ is material constitutive tensor.
% ==============================================================================
% A two-field variational statement
% ==============================================================================
\subsection{A two-field variational statement}

The proposed NURBS based hybrid elements are developed with the two-field Hellinger-Reissner variational formulation where stress and displacement are the field variables. The variational form of the stated governing equations is evaluated using the method of weighted residuals. Involvement of two field variables will necessitate the respective variation for each field. Let $\delta\bu$ and $\delta\btau$ be the variation of the displacement field $\bu$ and stress field $\btau$ respectively, in such a way that,
% ============== Variation of displacement and stress =========================
\begin{align}
V_{u} =& \{\delta\bu \in H^{1}(\Omega) : \delta\bu =\mathbf{0} \quad\text{on}\quad\Gamma_u \}\\
V_{\tau}=& \{\delta\btau \in L^{2}(\Omega) : \delta\btau=\delta\btau^{T}\}
\end{align}
% -----------------------------------------------------------------------------

If $(\delta\bu,\delta\btau)\in(V_u\bcross V_{\tau})$, then the weak form for the given governing equation can be written as,
% ============================= strong form ====================================
\begin{equation} \label{2_d_variation_form}
%\begin{split}
\int_\Omega \delta\bu\cdot(\del\cdot\btau + \boldsymbol{f})  \domega + \int_{\Gamma_t} \delta\bu\cdot(\tbar-\bt)  \dgamma + \int_\Omega \delta \btau\bcolon [\bepsilon^{\bu}- \bepsilon^{\btau}] \domega= 0
%\end{split}
\end{equation} 
% ------------------------------------------------------------------------------
where $\bepsilon^{\btau}$ denotes the strain tensor derived with the help of stress-strain relation and $\bepsilon^{\bu}$ is the strain tensor derived from the displacements using the Eqn.~\ref{strain_tensor_u}.

After satisfying the two conditions i.e.\ $(\delta\bu,\delta \btau)= (\delta\bu,\mathbf{0})$ and $(\delta\bu,\delta \btau)= (\mathbf{0},\delta \btau)$, the above equation (Eqn.~\ref{2_d_variation_form}) simplifies to the following forms, 
% =================== Final variational statements ============================
\begin{gather}
\int_\Omega [\bar{\beps}_c(\delta\bu)]^T\btau_c \domega = \int_\Omega \delta\bu^T \boldsymbol{f} \domega + \int_{\Gamma_t} \delta\bu^T \tbar \dgamma \quad\quad \forall\quad \delta\bu \in V_u \label{variational_statement_two_field_final_1} \\
\int_\Omega \delta\btau_c^T \left[ \bar{\beps}_c(\bu) - \bC^{-1}\btau_c \right] \domega = 0 \quad\quad \forall \quad \delta\btau_c \label{variational_statement_two_field_final_2}
\end{gather}
% ==============================================================================
where $\btau_c$ and $\bar{\beps}_c$ are the vector form of stress and strain tensor.

\subsection{Approximating functions}

\subsubsection{Interpolation functions for displacement field}

The current work uses the isoparametric concept where the displacement field ($\bu$) is approximated using the NURBS basis functions that are capable of maintaining the exact geometry. The displacement field ($\bu$) is interpolated as, 
% ===================== displacement approximation =============================
\begin{equation}
\bu = \sum_{I=1}^{n_{cp}^e} R_I \tilde{u}_I = \bR \utilde \label{displacement_hybrid}, \quad
\delta\bu = \sum_{I=1}^{n_{cp}^e} R_I \delta\tilde{u}_I = \bR\delta\utilde 
\end{equation}
% ------------------------------------------------------------------------------
where $R_I$ is the NURBS basis function, $I$ denotes the global numbering assigned to the CP, and $n_{cp}^e$ is the total number of CP per element.

\subsubsection{Stress interpolation functions} \label{section_stress_interpolation}

The next entity of interest in the present formulation is the stress interpolation matrix that benefits in approximating the second independent field i.e.\ $\btau_c$. Let $\btau_c$ and its variation ($\delta\btau_c$) be interpolated as,  
% ================== stress approximations =======================================
\begin{equation} \label{stress_interpolations}
\btau_{c}= \bP\betahat, \quad \delta\btau_{c}= \bP\delta\betahat
\end{equation}
% --------------------------------------------------------------------------------
where $\betahat$ is the vector consisting of the stress parameters for the respective element, $\delta\betahat$ is the vector of stress variation parameters, and $\bP$ is the stress interpolation matrix. The accuracy of the solution is highly sensitive towards the choice of $\bP$, hence, to ensure the efficient derivation of $\bP$ it is mandatory to evaluate the NURBS basis functions in its symbolic design. However, due to the recursive nature of the NURBS basis functions, special treatment needs to be followed to derive the symbolic expressions. At present, the widely recognized NURBS toolbox can only evaluate the basis function values at a given parametric point and does not provide the desired symbolic expressions. However, a mathematical computational software, Mathematica, provides a platform to develop a code to fulfill this provision. Taking advantage of this, a Mathematica code has been developed which is capable of deriving the symbolic expressions for the basis functions. All basis function and respective derivative in the section~\ref{Derivation_of_P} are evaluated using the Mathematica.

\subsection{Derivation of $\bP$ matrix for two-dimensional elements} \label{Derivation_of_P}

The $\bP$ matrices for consistently used element types are derived below. Furthermore, it should be noted that the matrix $\bP$ is independent of the number of elements even though the expressions are different in each knot span.

\subsubsection{A two-dimensional element with linear basis along $\xi$ and $\eta$ direction}
Let a single bi-linear IGA element be modeled using the geometric data provided in Table~\ref{surface_data_bilinear_element}, where $\boldsymbol{\Xi}$, $\bH$ and $p$, $q$ are the knot vectors and degree of basis function along $\xi$ and $\eta$ direction respectively.
% ===================== table: Rectangular plate ==============================
\begin{table}[pos=h]
\begin{minipage}{0.45\textwidth}
\centering
\caption{Surface data to construct a bi-linear element} \label{surface_data_bilinear_element}
\begin{tabular*}{1\textwidth}{@{} LL@{} }
\toprule 
\bf {Surface property} & \bf {Geometric data} \\
\midrule
$p$ & 1\\ 
$q$ & 1 \\ 
$\boldsymbol{\Xi}$ & $\left\{ \begin{array}{cccc} 0 & 0 & 1 &1 \end{array}\right\}$\\
$\bH$ & $\left\{ \begin{array}{cccc} 0 & 0 & 1 &1 \end{array}\right\}$ \\
\bottomrule
\end{tabular*}
\end{minipage}
\hfill
\begin{minipage}{0.45\textwidth}
\centering
\caption{Global numbering for one bi-linear element} \label{global_numbering}
\begin{tabular*}{1\textwidth}{@{} L|LLLL@{} }
\toprule 
$i$ & $1$ &$2$ & $1$ & $2$\\
$j$ & $1$ &$1$ & $2$ & $2$\\
\midrule
Global number ($I$) & $1$ &$2$ & $3$ & $4$\\
\bottomrule
\end{tabular*}
\end{minipage}
\end{table}
% -----------------------------------------------------------------------------
The CP associated with the geometric description of a model are treated as constants, which will only affect the physical domain but not the parametric space. The displacement $\bu$ is interpolated using the Eqn.~\ref{displacement_hybrid}. For the stated bi-linear element considered, there will be two basis functions in each $\xi$ and $\eta$ direction. Let $R_{i,p}(\xi)$ and $R_{j,q}(\eta)$ be the basis functions in $\xi$ and $\eta$ direction where $i=1,2$ and $j=1,2$. Global numbering is decided as mentioned in Table~\ref{global_numbering}.

$R_I$, represents the two-dimensional basis function for modeling the desired surface, is evaluated using Eqn.~\ref{Bivariate_NURBS_basis} in such a way that,
% ======================== formula for R_i ====================================
\begin{equation}\label{R_I_formula}
R_I =  R_{i,p}(\xi)R_{j,q}(\eta) = \frac{N_{i,p}(\xi) M_{j,q}(\eta) w_{i,j}}{\sum_{i=1}^{n}\sum_{j=1}^{m} N_{i,p}(\xi) M_{j,q}(\eta) w_{i,j}}
\end{equation}
% -----------------------------------------------------------------------------
For evaluation of $\bP$, weights are considered to be equal to 1 which replaces the Eqn.~\ref{R_I_formula} to,
% ===================== formula for R_I without weight ========================
\begin{equation}\label{R_I_evaluation}
R_I = R_{i,j}^{p,q}(\xi,\eta) = R_{i,p}(\xi)R_{j,q}(\eta) = N_{i,p}(\xi) M_{j,q}(\eta)
\end{equation}
%-------------------------------------------------------------------------------

The expressions for $R_I$ are evaluated using Mathematica as follows,
% ======================= N and M from Mathematica =============================
% For N ---------------------------
\begin{multicols}{2}
\noindent
\begin{align*}
N_{1,1} &= \begin{array}{cc}
 \bigg\{ 
\begin{array}{cc}
 1-\xi  & 0\leq \xi \leq 1 \\
 0 & \text{Otherwise} \\
\end{array}
 \\
\end{array}\\
N_{2,1} &= \begin{array}{cc}
 \bigg\{ 
\begin{array}{cc}
 \xi  & 0\leq \xi \leq 1 \\
 0 & \text{Otherwise} \\
\end{array}
 \\
\end{array}
\end{align*}
% For M ------------------------
\begin{align*}
M_{1,1} &= \begin{array}{cc}
 \bigg\{  
\begin{array}{cc}
 1-\eta  &  0\leq \eta \leq 1 \\
 0 & \text{Otherwise} \\
\end{array}
 \\
\end{array} \\
M_{2,1} &= \begin{array}{cc}
 \bigg\{ 
\begin{array}{cc}
 \eta  & 0\leq \eta \leq 1 \\
 0 & \text{Otherwise} \\
\end{array}
 \\
\end{array}
\end{align*}
\end{multicols}
\vspace{-0.8cm}
% ------------------------------------------------------------------------------
% ============== R = N*M ====================================================
\begin{subequations} \label{R_for_bi_linear_element}
\begin{align}
R_1 &= N_{1,1}M_{1,1} = \begin{array}{cc}
 \bigg\{  
\begin{array}{cc}
  (1-\eta ) (1-\xi ) &  0\leq \xi \leq 1; 0\leq \eta \leq 1\\
 0 & \text{Otherwise} \\
\end{array}
 \\
\end{array} \label{R_1_bilinear_single_element}\\
R_2 &= N_{2,1}M_{1,1} = \begin{array}{cc}
 \bigg\{  
\begin{array}{cc}
 (1-\eta ) \xi  &  0\leq \xi \leq 1; 0\leq \eta \leq 1 \\
 0 & \text{Otherwise} \\
\end{array}
 \\
\end{array} \label{R_2_bilinear_single_element}\\
R_3 &= N_{1,1}M_{2,1} = \begin{array}{cc}
 \bigg\{ 
\begin{array}{cc}
 \eta  (1-\xi ) &  0\leq \xi \leq 1; 0\leq \eta \leq 1 \\
 0 & \text{Otherwise} \\
\end{array}
 \\
\end{array} \label{R_3_bilinear_single_element}\\
R_4 &= N_{2,1}M_{2,1} = \begin{array}{cc}
 \bigg\{ 
\begin{array}{cc}
 \eta  \xi &  0\leq \xi \leq 1; 0\leq \eta \leq 1 \\
 0 & \text{Otherwise} \\
\end{array}
 \\
\end{array}\label{R_4_bilinear_single_element}
\end{align}
\end{subequations}

Substituting the expressions for $R_I$ from Eqn.~\ref{R_for_bi_linear_element} into Eqn.~\ref{displacement_hybrid} will lead to a bi-linear expression of $\bu$ where the coefficients associated with each constant can be collected as $\left\{ 1,\: \eta, \: \xi, \: \xi\eta \right\}$. The normal stress components are obtained by differentiating the displacement interpolation functions with respect to natural co-ordinates such that the expressions for $\tau_{\xi\xi}$ and $\tau_{\eta\eta}$ will be the linear combination of the terms $\left\{ 1, \: \eta \right\} $ and  $ \left\{  1, \: \xi \right\}$ respectively. Finally, the shear stress components are obtained so that they will suppress any spurious zero-energy mode \cite{Jog2010} which will lead to an expression of $\tau_{\xi\eta}$ as a constant.

Introduce the constant $\beta_i$ ($i = 1,2, \hdots,5$) in such a way that no $\beta_i$ term is shared, which will result in following expression of stresses,
% ================== components of stress =============================
\begin{subequations} \label{stress_equations_for_bi_linear_element_1}
\begin{multicols}{3}
\noindent
\begin{equation*}
\tau_{\xi\xi} = \beta_1 + \beta_2 \eta, 
\end{equation*}
\begin{equation*}
\tau_{\eta\eta} = \beta_3 + \beta_4 \eta, 
\end{equation*}
\begin{equation*}
\tau_{\xi\eta} = \beta_5 
\end{equation*}
\end{multicols} 
\end{subequations}
% ----------------------------------------------------------------
And can be further written as,
% ================== Tau_xi_eta =============================
\begin{equation}
\underbrace{\begin{Bmatrix} \tau_{\xi\xi} \\ \tau_{\eta\eta} \\\tau_{\xi\eta} \end{Bmatrix}}_{\btau_c(\xi,\eta)} = 
\underbrace{\begin{bmatrix} 1 & \eta & 0 & 0 & 0 \\0 & 0 & 1 & \xi & 0 \\ 0 & 0 & 0 & 0 & 1 \end{bmatrix}}_{\bP(\xi,\eta)}
\underbrace{\begin{Bmatrix} \beta_1 \\ \beta_2 \\ \beta_3 \\ \beta_4 \\ \beta_5 \end{Bmatrix}}_{\boldsymbol{\beta}} 
\end{equation} 
% ----------------------------------------------------------------

Considering the fact that, for evaluation of integrals, it is important to define the stress components in the master space. Parametric space and master space are related with a linear mapping. Hence, stress components in a master space will have the same form but with different constants. Therefore, the stress components in a master space ($\tilde{\xi}-\tilde{\eta}$) are given as,
% ============================ tau in xi and eta =====================
\begin{equation}
\underbrace{\begin{Bmatrix} \tau_{\tilde{\xi}\tilde{\xi}} \\ \tau_{\tilde{\eta}\tilde{\eta}} \\\tau_{\tilde{\xi}\tilde{\eta}} \end{Bmatrix}}_{\btau_c(\tilde{\xi},\tilde{\eta})} = 
\underbrace{\begin{bmatrix} 1 & \tilde{\eta} & 0 & 0 & 0 \\0 & 0 & 1 & \tilde{\xi} & 0 \\ 0 & 0 & 0 & 0 & 1 \end{bmatrix}}_{\bP(\tilde{\xi},\tilde{\eta})} 
\underbrace{\begin{Bmatrix} \hat{\beta}_1 \\ \hat{\beta}_2 \\ \hat{\beta}_3 \\ \hat{\beta}_4 \\ \hat{\beta}_5 \end{Bmatrix}}_{\boldsymbol{\hat\beta}} 
\quad\rightarrow\quad \btau_c(\tilde{{\xi}},\tilde{\eta)} = \bP(\tilde{\xi},\tilde{\eta})\boldsymbol{\hat{\beta}}
\end{equation}
% ---------------------------------------------------------------------

Finally, the stress components in master space are related with the physical space with the following transformation,
% ========================= tau in master space ======================
\begin{equation}
\btau_c(x,y) = \bT \btau_c(\tilde{\xi},\tilde{\eta})
\end{equation}
% -------------------------------------------------------------------- 
where $\bT$ is the transformation matrix derived from a combination of Jacobians relating the corresponding spaces. Jacobian which relates the master space and the physical space is $J = J_2 J_1$ , where $J_2$ and $J_1$ are the Jacobians for mapping master space $(\tilde{\xi},\tilde{\eta})$ to parametric space $(\xi,\eta)$ and parametric space to physical space $(x,y)$ respectively \cite{Agrawal2018}. The transformation matrix ($\bT$), for two-dimensional problems, is evaluated as follows, \\
% ===================== transformation matrix 2d =======================
\begin{equation} \label{transformation_matrix_2d}
\bT = \begin{bmatrix} J_{11}^2 & J_{21}^2 & 2J_{11}J_{21} \\
J_{12}^2 & J_{22}^2 & 2J_{12}J_{22} \\
J_{11}J_{12} & J_{21}J_{22} & J_{11}J_{22}+J_{12}J_{21} \\
\end{bmatrix}
\end{equation}
where $J_{ij}$ are the components of combined Jacobian ($J$).

It is essential to ensure that the $\bP(\tilde{\xi},\tilde{\eta})$ matrix is independent of $h$-refinement even though the expressions are different in each knot span due to the involvement of NURBS basis function. The requirement in deriving the $\bP$ matrix resides in the coefficients associated with each constant term in displacement ($\bu$). These coefficients will not change element-wise; however, the constants will change, but the form will remain the same. The changes in constants will be taken care by an appropriate transformation matrix.

\begin{table}[pos=h]
\begin{minipage}{0.45\textwidth}
\caption{Surface data to construct two bi-linear element} \label{surface_data_two_bilinear_element}
\begin{tabular*}{1\textwidth}{@{} LL@{} }
\toprule 
\bf {Surface property} & \bf {Geometric data} \\
\midrule
$p$ & 1\\ 
$q$ & 1 \\ 
$\boldsymbol{\Xi}$ & $ \left\{ \begin{array}{ccccc} 0 & 0 & 1\slash 2 & 1 & 1 \end{array}\right\}$ \\
$\bH$ & $ \left\{ \begin{array}{cccc} 0 & 0 & 1 & 1 \end{array}\right\}$ \\
\bottomrule
\end{tabular*}
\end{minipage}
\hfill
\begin{minipage}{0.45\textwidth}
\centering
\caption{Global numbering for two bi-linear elements} \label{global_numbering_two_elements}
\begin{tabular*}{1\textwidth}{@{} L|LLLLLL@{} }
\toprule 
$i$ & $1$ &$2$ &$3$ & $1$ & $2$ &$3$\\
$j$ & $1$ &$1$ & $1$ & $2$ & $2$ & $2$\\
\midrule
Global number ($I$) & $1$ &$2$ & $3$ & $4$ & $5$ & $6$\\
\bottomrule
\end{tabular*}
\end{minipage}
\end{table}

To validate the above statement, consider a mesh of two bi-linear elements obtained by performing $h$-refinement or knot insertion into the existing data given in Table~\ref{surface_data_bilinear_element}. This will result in the geometric description illustrated in Table~\ref{surface_data_two_bilinear_element}. For stated setting of bi-linear elements, there will be three basis functions in $\xi$ direction and two basis in $\eta$ direction and the global numbering is decided as per Table~\ref{global_numbering_two_elements}. The required basis functions are evaluated as follows,
% ======================= N and M from Mathematica =============================
% For N ---------------------------
\begin{equation*}
N_{1,1} = \begin{array}{cc}
 \bigg\{
\begin{array}{cc}
 1-2 \xi  & 0\leq \xi \leq \frac{1}{2} \\
 0 & \text{Otherwise} \\
\end{array}
 \\
\end{array}, \\
N_{2,1} = \begin{array}{cc}
 \Bigg\{
\begin{array}{cc}
 -2 (\xi -1) & \frac{1}{2}\leq \xi \leq 1 \\
 2 \xi  & 0\leq \xi <\frac{1}{2} \\
0 & \text{Otherwise} \\
\end{array}
 \\
\end{array}, \\
N_{3,1} = \begin{array}{cc}
 \bigg\{
\begin{array}{cc}
 2 \xi -1 & \frac{1}{2}\leq \xi \leq 1 \\
 0 & \text{Otherwise} \\
\end{array}
 \\
\end{array}
\end{equation*}
%\columnbreak
% For M ------------------------
\begin{equation*}
M_{1,1} = \begin{array}{cc}
 \bigg\{
\begin{array}{cc}
 1-\eta  &  0\leq \eta \leq 1 \\
 0 & \text{Otherwise} \\
\end{array}
 \\
\end{array}, \quad
M_{2,1} = \begin{array}{cc}
 \bigg\{ 
\begin{array}{cc}
 \eta  & 0\leq \eta \leq 1 \\
 0 & \text{Otherwise} \\
\end{array}
 \\
\end{array}
\end{equation*}
%\end{multicols}
% ------------------------------------------------------------------------------ 

\noindent \underline{\textbf{\textit{Element 1:}}} $0\leq \xi \leq \frac{1}{2}; 0\leq \eta \leq 1$, the participating basis functions are,
% ============== R = N*M ====================================================
\begin{subequations}\label{R_eqn_two_linear_elements_first_element}
\begin{multicols}{2}
\noindent
\begin{align}
R_1 &= N_{1,1}M_{1,1} = (1-\eta ) (1-2 \xi ) \label{R_eqn_1_element_1}\\
R_2 &= N_{2,1}M_{1,1} = 2 (1-\eta ) \xi \label{R_eqn_2_element_1}
\end{align}%
\begin{align}
R_4 &= N_{1,1}M_{2,1} = \eta  (1-2 \xi ) \label{R_eqn_3_element_1}\\
R_5 &= N_{3,1}M_{2,1} = 2 \eta  \xi \label{R_eqn_4_element_1}
\end{align}
\end{multicols}
\end{subequations}

Substituting the expressions for $R$'s from Eqn.~\ref{R_eqn_two_linear_elements_first_element} into Eqn.~\ref{displacement_hybrid} will lead to a bi-linear expression of $\bu$ where the coefficients associated with each constant can be collected as $\left\{ 1,\: \eta, \: \xi, \: \xi\eta \right\}$. Following the similar approach illustrated for a single element, the stress components in a master space ($\tilde{\xi}-\tilde{\eta}$) are evaluated as,
\begin{equation}
\underbrace{\begin{Bmatrix} \tau_{\tilde{\xi}\tilde{\xi}} \\ \tau_{\tilde{\eta}\tilde{\eta}} \\\tau_{\tilde{\xi}\tilde{\eta}} \end{Bmatrix}}_{\btau_c(\tilde{\xi},\tilde{\eta})} = 
\underbrace{\begin{bmatrix} 1 & \tilde{\eta} & 0 & 0 & 0 \\0 & 0 & 1 & \tilde{\xi} & 0 \\ 0 & 0 & 0 & 0 & 1 \end{bmatrix}}_{\bP(\tilde{\xi},\tilde{\eta})} 
\underbrace{\begin{Bmatrix} \hat{\beta}_1 \\ \hat{\beta}_2 \\ \hat{\beta}_3 \\ \hat{\beta}_4 \\ \hat{\beta}_5 \end{Bmatrix}}_{\boldsymbol{\hat\beta}}  
\end{equation}
% --------------------------------------------------------------------- 

\noindent \underline{\textbf{\textit{Element 2:}}} $\frac{1}{2}\leq \xi \leq 1 ; 0\leq \eta \leq 1$, the participating basis functions are,
% ============== R = N*M ====================================================
\begin{subequations}\label{R_eqn_two_linear_elements_second_element}
\begin{multicols}{2}
\noindent
\begin{align}
R_2 &= N_{2,1}M_{1,1} = -2 (1-\eta ) (\xi -1) \label{R_eqn_1_element_1}\\
R_3 &= N_{3,1}M_{1,1} = (1-\eta ) (2 \xi -1) \label{R_eqn_2_element_1}
\end{align}
\begin{align}
R_5 &= N_{2,1}M_{2,1} = -2 \eta  (\xi -1) \label{R_eqn_3_element_1}\\
R_6 &= N_{3,1}M_{2,1} =\eta  (2 \xi -1) \label{R_eqn_4_element_1}
\end{align}
\end{multicols}
\end{subequations}

Substituting the Eqns.~\ref{R_eqn_two_linear_elements_second_element} into Eqn.~\ref{displacement_hybrid} will lead to a similar bi-linear expression of $\bu$ with the coefficients as $\left\{ 1,\: \eta, \: \xi, \: \xi\eta \right\}$. And the stress components in a master space ($\tilde{\xi}-\tilde{\eta}$) are evaluated as,
% ============================ tau in xi and eta =====================
\begin{equation}
\underbrace{\begin{Bmatrix} \tau_{\tilde{\xi}\tilde{\xi}} \\ \tau_{\tilde{\eta}\tilde{\eta}} \\\tau_{\tilde{\xi}\tilde{\eta}} \end{Bmatrix}}_{\btau_c(\tilde{\xi},\tilde{\eta})} = 
\underbrace{\begin{bmatrix} 1 & \tilde{\eta} & 0 & 0 & 0 \\0 & 0 & 1 & \tilde{\xi} & 0 \\ 0 & 0 & 0 & 0 & 1 \end{bmatrix}}_{\bP(\tilde{\xi},\tilde{\eta})} 
\underbrace{\begin{Bmatrix} \hat{\beta}_1 \\ \hat{\beta}_2 \\ \hat{\beta}_3 \\ \hat{\beta}_4 \\ \hat{\beta}_5 \end{Bmatrix}}_{\boldsymbol{\hat\beta}} 
\end{equation}
% ---------------------------------------------------------------------

In both the cases, even though the participating basis functions are different but the resulting stress interpolation matrix $\bP(\tilde{\xi},\tilde{\eta})$ in master space remains the same.

% ======================= subsection quadratic  element ====================
\subsubsection{A two-dimensional element with quadratic basis along $\xi$ and $\eta$ direction}
Let a single bi-quadratic element can be modeled using the geometric data given in Table~\ref{surface_data_bi_quadratic_element}. 
\begin{table}[pos=h]
\begin{minipage}{0.45\textwidth}
\caption{Surface data to construct a bi-quadratic element} \label{surface_data_bi_quadratic_element}
\begin{tabular*}{1\textwidth}{@{} LL@{} }
\toprule 
\bf {Surface property} & \bf {Geometric data} \\
\midrule
$p$ & 2\\ 
$q$ & 2 \\ 
$\boldsymbol{\Xi}$ & $\left\{\begin{array}{cccccc} 0 & 0 & 0 & 1 & 1 & 1 \end{array}\right\}$ \\
$\bH$ &  $\left\{\begin{array}{cccccc} 0 & 0 & 0 & 1 & 1 & 1 \end{array}\right\}$ \\
\bottomrule
\end{tabular*}
\end{minipage}
\hfill
\begin{minipage}{0.45\textwidth}
\centering
\caption{Global numbering for a bi-quadratic element} \label{global_numbering_bi_quadratic_element}
\begin{tabular*}{1\textwidth}{@{} L|LLLLLLLLL@{} }
\toprule 
$i$ & $1$ &$2$ &$3$ & $1$ & $2$ &$3$ & $1$ &$2$ &$3$\\
$j$ & $1$ &$1$ & $1$ & $2$ &$2$ & $2$ & $3$ &$3$ & $3$\\
\midrule
Global number ($I$) & $1$ &$2$ & $3$ & $4$ & $5$ &$6$ & $7$ & $8$ & $9$\\
\bottomrule
\end{tabular*}
\end{minipage}
\end{table}
The NURBS basis function will be evaluated using the tensor product of univariate basis function along $\xi$ direction ($N_{i,p}(\xi)$ for $i=1,2,3$) and $\eta$ direction ($M_{j,q}(\eta)$ for $j=1,2,3$). The global numbering is decided as per Table~\ref{global_numbering_bi_quadratic_element}. $R_I$ represents the bivariate basis function for modeling the desired surface, and evaluated using Eqn.~\ref{R_I_evaluation}. The expressions for $N_{i,p}(\xi)$ and $M_{j,q}(\eta)$ are calculated as,
% ======================= N and M from Mathematica =============================
% For N ---------------------------
\begin{equation*}
N_{1,2} = \begin{array}{cc}
 \Bigg\{
\begin{array}{cc}
 \xi ^2-2 \xi +1 & 0\leq \xi \leq 1 \\
 0 & \text{Otherwise} \\
\end{array},
 \\
\end{array} \\
N_{2,2} = \begin{array}{cc}
 \Bigg\{
\begin{array}{cc}
 -2 \left(\xi ^2-\xi \right) & 0\leq \xi \leq 1 \\
 0 & \text{Otherwise} \\
\end{array},
 \\
\end{array} \\
N_{3,2} = \begin{array}{cc}
 \Bigg\{
\begin{array}{cc}
 \xi ^2 & 0\leq \xi \leq 1 \\
 0 & \text{Otherwise} \\
\end{array}
 \\
\end{array} 
\end{equation*}
% For M ------------------------
\begin{equation*}
M_{1,2}=\begin{array}{cc}
 \Bigg\{
\begin{array}{cc}
 \eta ^2-2 \eta +1 & 0\leq \eta \leq 1 \\
 0 & \text{Otherwise} \\
\end{array},
 \\
\end{array}\\
M_{2,2}= \begin{array}{cc}
 \Bigg\{ & 
\begin{array}{cc}
 -2 \left(\eta ^2-\eta \right) & 0\leq \eta \leq 1 \\
 0 & \text{Otherwise} \\
\end{array},
 \\
\end{array}\\
M_{3,2} = \begin{array}{cc}
 \Bigg\{
\begin{array}{cc}
 \eta ^2 & 0\leq \eta \leq 1 \\
 0 & \text{Otherwise} \\
\end{array}
 \\
\end{array}
\end{equation*}
% ------------------------------------------------------------------------------
% ============== R = N*M ====================================================
\begin{subequations}\label{R_I_quadratic_all_eqn}
\begin{align}
R_1 &= N_{1,2}M_{1,2} = \left(\eta ^2-2 \eta +1\right) \left(\xi ^2-2 \xi +1\right) \label{first_eqn_quad}\\
R_2 &= N_{2,2}M_{1,2} = -2 \left(\eta ^2-2 \eta +1\right) \left(\xi ^2-\xi \right) \\
R_3 &= N_{3,2}M_{1,2} = \left(\eta ^2-2 \eta +1\right) \xi ^2 \\
R_4 &= N_{1,2}M_{2,2} = -2 \left(\eta ^2-\eta \right) \left(\xi ^2-2 \xi +1\right)\\
R_5 &= N_{2,2}M_{2,2} = 4 \left(\eta ^2-\eta \right) \left(\xi ^2-\xi \right) \\
R_6 &= N_{3,2}M_{2,2} =  -2 \left(\eta ^2-\eta \right) \xi ^2 \\
R_7 &= N_{1,2}M_{3,2} =  \eta ^2 \left(\xi ^2-2 \xi +1\right) \\
R_8 &= N_{2,2}M_{3,2} =  -2 \eta ^2 \left(\xi ^2-\xi \right)  \\
R_9 &= N_{3,2}M_{3,2} = \eta ^2 \xi ^2 \label{Last_eqn_quad}
\end{align}
\end{subequations}

Substituting the expressions for $R_I$ from Eqn.~\ref{R_I_quadratic_all_eqn} into Eqn.~\ref{displacement_hybrid} will lead to a bi-quadratic expression of $\bu$ where the coefficients associated with each constant can be collected as $\left\{1, \:\xi, \: \eta , \:\xi^2, \: \xi\eta, \:\eta^2, \:\xi^2\eta, \:\eta^2\xi, \:\xi^2\eta^2 \right\}$. The normal stress components are obtained by differentiating the displacement interpolation functions with respect to natural co-ordinates such that the expressions for $\tau_{\xi\xi}$ and $\tau_{\eta\eta}$ will be the linear combination of the terms $\left\{1, \: \xi, \:\eta, \:\xi\eta, \:\eta^2, \:\eta^2\xi \right\} $ and  $\left\{1, \: \xi, \:\eta, \:\xi\eta, \:\xi^2, \: \eta\xi^2\right\} $ respectively. Finally, the shear stress components are obtained so that they will suppress any spurious zero-energy mode \cite{Jog2010} which will lead to an expression of shear stress $\tau_{\xi\eta}$ as a linear combination of the terms $ \left\{1, \: \xi, \:\eta, \:\xi\eta \right\}$.

Introducing the constant $\beta_i$ ($i = 1,2,\hdots 16$), we will have the following expression of stresses,
% ================== components of stress =============================
\begin{subequations}\label{stresses_for_quad_element_all}
\begin{gather}
\tau_{\xi\xi} = \beta_1 + \beta_2 \xi + \beta_3 \eta + \beta_4 \xi\eta + \beta_5 \eta^2 + \beta_6 \eta^2\xi \\
\tau_{\eta\eta} = \beta_7 + \beta_8 \xi + \beta_9 \eta + \beta_{10} \xi\eta + \beta_{11} \xi^2 + \beta_{12} \eta\xi^2 \\
\tau_{\xi\eta} =\beta_{13} + \beta_{14} \xi + \beta_{15} \eta + \beta_{16} \xi\eta
\end{gather} 
\end{subequations}
% ----------------------------------------------------------------
Moreover, the stress components in a master space ($\tilde{\xi}-\tilde{\eta}$) are given as,
% ================== Tau_xi_eta =============================
\begin{equation}
\begin{Bmatrix} \tau_{\tilde{\xi}\tilde{\xi}} \\ \tau_{\tilde{\eta}\tilde{\eta}} \\\tau_{\tilde{\xi}\tilde{\eta}} \end{Bmatrix} = 
\underbrace{\begin{bmatrix} 
1 & \tilde{\xi} & \tilde{\eta} & \tilde{\xi} \tilde{\eta} &  \tilde{\eta}^2 & \tilde{\eta}^2 \tilde{\xi} & 0 & 0 & 0 & 0 & 0 & 0 & 0 & 0 & 0 & 0\\
0 & 0 & 0 & 0 & 0 & 0 & 1 & \tilde{\xi} & \tilde{\eta} & \tilde{\xi} \tilde{\eta} & \tilde{\xi}^2 & \tilde{\eta} \tilde{\xi}^2 & 0 & 0  & 0 & 0 \\
0 & 0 & 0 & 0 & 0 & 0 & 0 & 0 & 0 & 0 & 0 & 0 & 1 & \tilde{\xi} & \tilde{\eta} & \tilde{\xi} \tilde{\eta}\\
\end{bmatrix}}_{\bP(\tilde{\xi},\tilde{\eta})}
\underbrace{\begin{Bmatrix} \hat{\beta}_1 \\ \hat{\beta}_2 \\ \vdots \\ \hat{\beta}_{16} \end{Bmatrix}}_{\boldsymbol{\hat{\beta}}} 
\end{equation} 

% ======================= subsection cubic  element ====================
\subsubsection{A two-dimensional element with cubic basis along $\xi$ and $\eta$ direction}

Following the similar approach, a single bi-cubic element is modeled using the geometric data given in Table~\ref{surface_data_bi_cubic_element}. 
% ===================== table: Rectangular plate ==============================
\begin{table}[width=.5\linewidth,cols=2,pos=h]
 \caption{Surface data to construct a bi-cubic element} \label{surface_data_bi_cubic_element}
\begin{tabular*}{\tblwidth}{@{} LL@{} }
\toprule 
\bf {Surface property} & \bf {Geometric data} \\
\midrule
$p$ & 3\\ 
$q$ & 3 \\ 
$\boldsymbol{\Xi}$ & $\left\{\begin{array}{cccccccc}0&0&0&0&1&1&1&1\end{array}\right\}$ \\
$\bH$ & $\left\{\begin{array}{cccccccc}0&0&0&0&1&1&1&1\end{array}\right\}$ \\
\bottomrule
\end{tabular*}
\end{table}
% -----------------------------------------------------------------------------
Univariate basis functions involved in deriving the bivariate NURBS basis functions are $N_{i,p}(\xi)$ and $M_{j,q}(\eta)$ where $i$ and $j$ varies from 1 to 4. Then, the resultant two-dimensional basis function can be evaluated as,
% ============== R = N*M ====================================================
\begin{subequations}\label{All_RI_for_bicubic_element}
\begin{align}
R_1 &= N_{1,3}M_{1,3} = \left(3 \eta ^2-\eta ^3-3 \eta +1\right) \left(3 \xi ^2-\xi ^3-3 \xi +1\right) \label{first_eqn_cubic}\\
 R_2 &= N_{2,3}M_{1,3} = 3 \left(-\eta ^3+3 \eta ^2-3 \eta +1\right) \left(\xi ^3-2 \xi ^2+\xi \right) \\
 R_3 &= N_{3,3}M_{1,3} = -3 \left(3 \eta ^2-\eta ^3-3 \eta +1\right) \left(\xi ^3-\xi ^2\right) \\
 R_4 &= N_{4,3}M_{1,3} = \left(-\eta ^3+3 \eta ^2-3 \eta +1\right) \xi ^3 \\
 R_5 &= N_{1,3}M_{2,3} = 3 \left(\eta ^3-2 \eta ^2+\eta \right) \left(3 \xi ^2-\xi ^3-3 \xi +1\right) \\
 R_6 &= N_{2,3}M_{2,3} = 9 \left(\eta ^3-2 \eta ^2+\eta \right) \left(\xi ^3-2 \xi ^2+\xi \right) \\
 R_7 &= N_{3,3}M_{2,3} = -9 \left(\eta ^3-2 \eta ^2+\eta \right) \left(\xi ^3-\xi ^2\right) \\
 R_8 &= N_{4,3}M_{2,3} = 3 \left(\eta ^3-2 \eta ^2+\eta \right) \xi ^3 \\
 R_9 &= N_{1,3}M_{3,3} = -3 \left(\eta ^3-\eta ^2\right) \left(-\xi ^3+3 \xi ^2-3 \xi +1\right) \\
 R_{10} &= N_{2,3}M_{3,3} = -9 \left(\eta ^3-\eta ^2\right) \left(\xi ^3-2 \xi ^2+\xi \right) \\
 R_{11} &= N_{3,3}M_{3,3} = 9 \left(\eta ^3-\eta ^2\right) \left(\xi ^3-\xi ^2\right) \\
 R_{12} &= N_{4,3}M_{3,3} = -3 \left(\eta ^3-\eta ^2\right) \xi ^3 \\
 R_{13} &= N_{1,3}M_{4,3} = \eta ^3 \left(-\xi ^3+3 \xi ^2-3 \xi +1\right) \\
 R_{14} &= N_{2,3}M_{4,3} = 3 \eta ^3 \left(\xi ^3-2 \xi ^2+\xi \right) \\
 R_{15} &= N_{3,3}M_{4,3} = -3 \eta ^3 \left(\xi ^3-\xi ^2\right) \\
 R_{16} &= N_{4,3}M_{4,3} = \eta ^3 \xi ^3 \label{Last_eqn_cubic}
\end{align}
\end{subequations}
% ----------------------------------------------------------------
This will lead to the stress interpolation matrix in a master space ($\tilde{\xi}-\tilde{\eta}$) as,
% ================== Tau_xi_eta =============================
\begin{equation}
\begin{aligned}
\underbrace{\begin{Bmatrix} \tau_{\tilde{\xi}\tilde{\xi}} \\ \tau_{\tilde{\eta}\tilde{\eta}} \\\tau_{\tilde{\xi}\tilde{\eta}} \end{Bmatrix}}_{\btau_c(\tilde{{\xi}},\tilde{\eta)}} &= 
\underbrace{\begin{bmatrix} 
\bP_1 & \bP_2 & \bP_3\\
\bP_4 & \bP_5 & \bP_6 \\
\bP_7 & \bP_8 & \bP_9
\end{bmatrix}}_{\bP(\tilde{\xi},\tilde{\eta})}
\underbrace{\begin{Bmatrix} \hat{\beta}_1 \\ \hat{\beta}_2 \\ \vdots \\ \hat{\beta}_{33} \end{Bmatrix}}_{\boldsymbol{\hat{\beta}}} 
\quad\rightarrow\quad \btau_c(\tilde{{\xi}},\tilde{\eta)} = \bP(\tilde{\xi},\tilde{\eta})\boldsymbol{\hat{\beta}}
\end{aligned}
\end{equation}  
% ----------------------------------------------------------------
where, $\bP_2$, $\bP_4$, $\bP_7$, $\bP_8$ are the zero vector of size $1\times12$, $\bP_3$, $\bP_6$ are the zero vectors of size $1\times 9$, and $\bP_1$, $\bP_5$, and $\bP_9$ are evaluated as, 
\begin{align*}
\bP_1 &= \left[ 1,\tilde{\xi},\tilde{\eta},\tilde{\xi}^2,\tilde{\eta}\tilde{\xi},\tilde{\eta}^2,\tilde{\eta}\tilde{\xi}^2,\tilde{\eta}^2\tilde{\xi},\tilde{\eta}^2\tilde{\xi}^2,\tilde{\eta}^3\tilde{\xi}^2,\tilde{\eta}^3\tilde{\xi},\tilde{\eta}^3  \right] \\
\bP_5 &= \left[ 1,\tilde{\xi},\tilde{\eta},\tilde{\xi}^2,\tilde{\eta}\tilde{\xi},\tilde{\eta}^2,\tilde{\eta}\tilde{\xi}^2,\tilde{\eta}^2\tilde{\xi},\tilde{\eta}^2\tilde{\xi}^2,\tilde{\eta}\tilde{\xi}^3,\tilde{\eta}^2\tilde{\xi}^3,\tilde{\xi}^3\right] \\
\bP_9 &= \left[1,\tilde{\xi},\tilde{\eta},\tilde{\xi}^2,\tilde{\eta}\tilde{\xi},\tilde{\eta}^2,\tilde{\eta}\tilde{\xi}^2,\tilde{\eta}^2\tilde{\xi},\tilde{\eta}^2\tilde{\xi}^2\right]
\end{align*}

% ==================================================================
% Elemental IGA equations for two-field variation principle
% ==================================================================
\subsection{IGA equilibrium equations for two-field variation principle}

Recollecting the expressions for $\bu$ (Eqn.~\ref{displacement_hybrid}), $\btau_c$ (Eqn.~\ref{stress_interpolations}), and if the strains are defined as $\bar\beps_c(\bu) = \bB\tilde{\bu}$ where $\bB$ is the strain displacement matrix corresponding to the NURBS interpolation functions, then the weak statement for the two-field variational statement (Eqn.~\ref{variational_statement_two_field_final_1} and~\ref{variational_statement_two_field_final_2}) will reduce to,
% =================== Weak form Matrix =============================
\begin{equation} \label{weak_form_matrix}
\begin{bmatrix} -\bH & \bG \\
\bG^T & \boldsymbol{0}
\end{bmatrix} \begin{Bmatrix} \hat\bbeta \\ \tilde{\bu} \end{Bmatrix} = \begin{Bmatrix} \hat\bg \\ \hat{\boldsymbol{f}} \end{Bmatrix}
\end{equation}
%==================================================================
where,
\begin{equation*}
\bG =\int_\Omega \bP^T \bB \domega,\quad \bH = \int_\Omega \bP^T\bS\bP \domega, \quad
\hat\bg = \boldsymbol{0}, \quad \hat{\boldsymbol{f}}  = \int_\Omega \bN^T \boldsymbol{f} \domega + \int_{\Gamma_t} \bN^T \tbar \dgamma \label{f_vector}
\end{equation*}
Evaluating the expression for $\hat\bbeta$ from Eqn.~\ref{weak_form_matrix} will lead to,
\begin{gather}
\underbrace{\bG^T \bH^{-1} \bG}_{\bK} \tilde{\bu} = \hat{\boldsymbol{f}} \rightarrow \bK\tilde{\bu} = \hat{\boldsymbol{f}}
\end{gather}
where the term $\bG^T \bH^{-1} \bG$ represents the stiffness matrix for two-field hybrid stress formulation. 

% ==============================================================================
% Section- Mapping, force vector and boundary conditions
% ==============================================================================
\subsection{Mapping associated with IGA}

\begin{figure}
\begin{subfigure}{0.35\columnwidth}
\centering
\includegraphics[width=1\columnwidth]{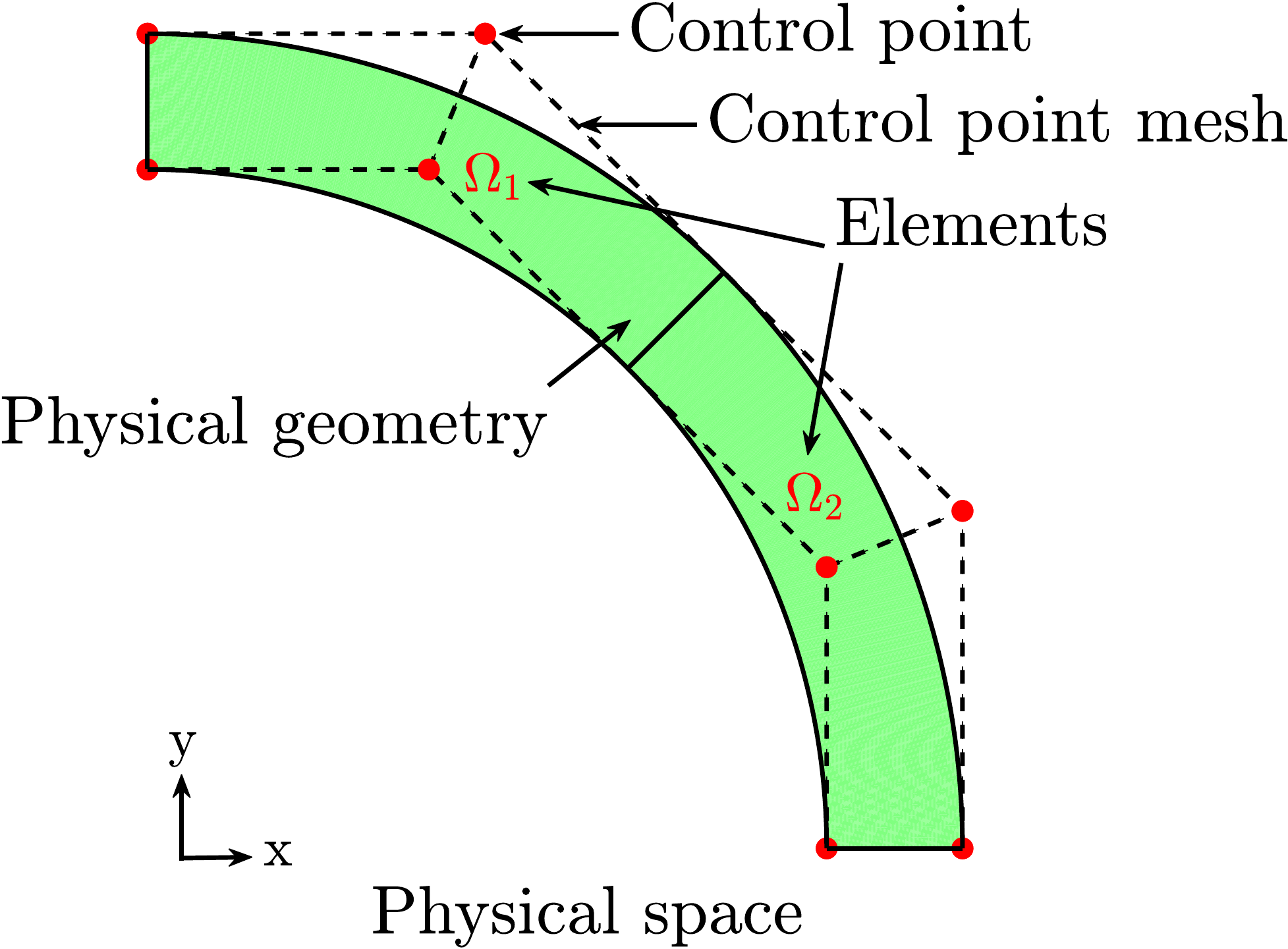}
\end{subfigure}%
\begin{subfigure}{0.3\columnwidth}
\centering
\includegraphics[width=0.9\columnwidth]{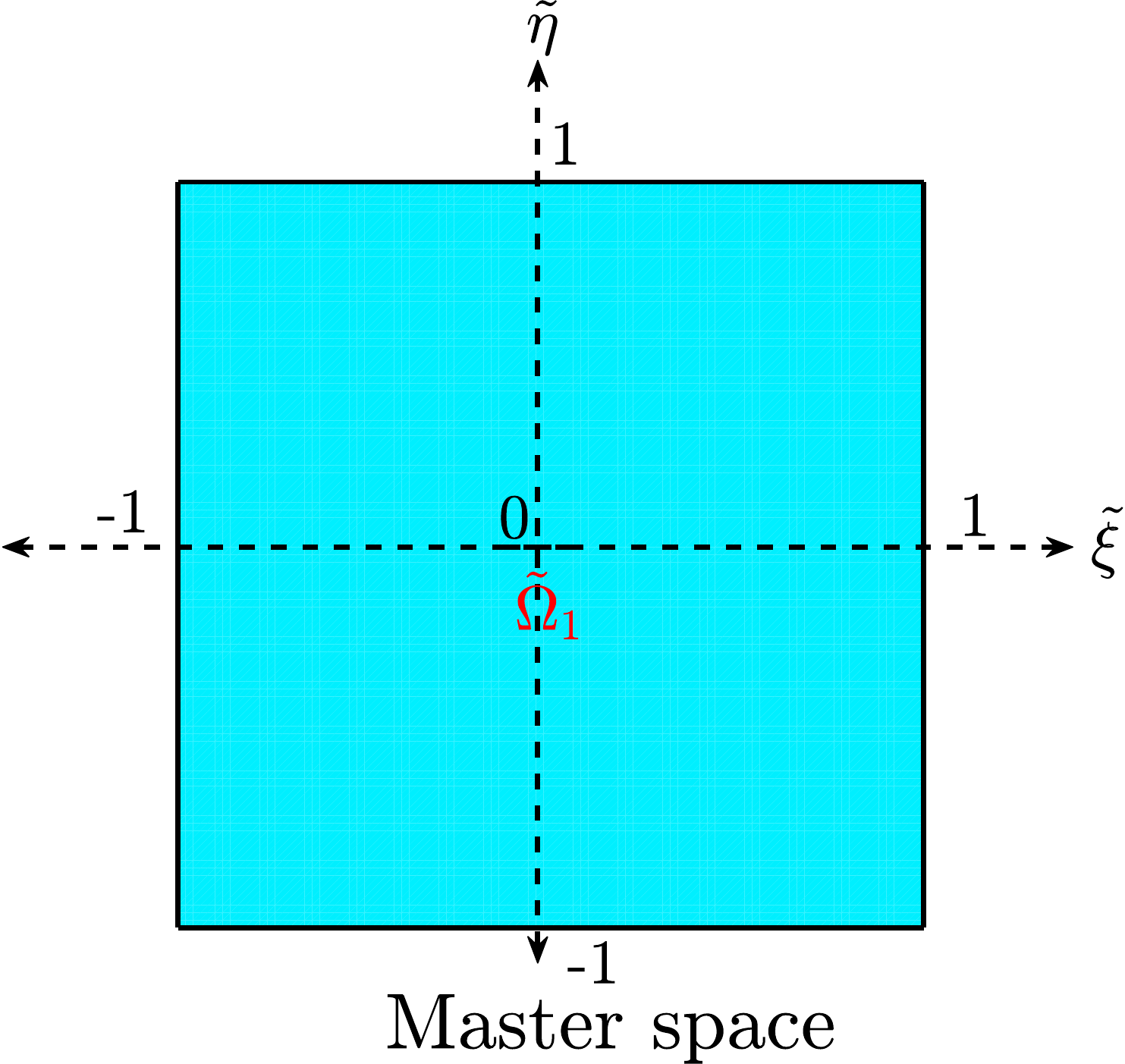}
\end{subfigure}
\newline 
\begin{subfigure}{1\columnwidth}
\centering
\includegraphics[width=0.425\columnwidth]{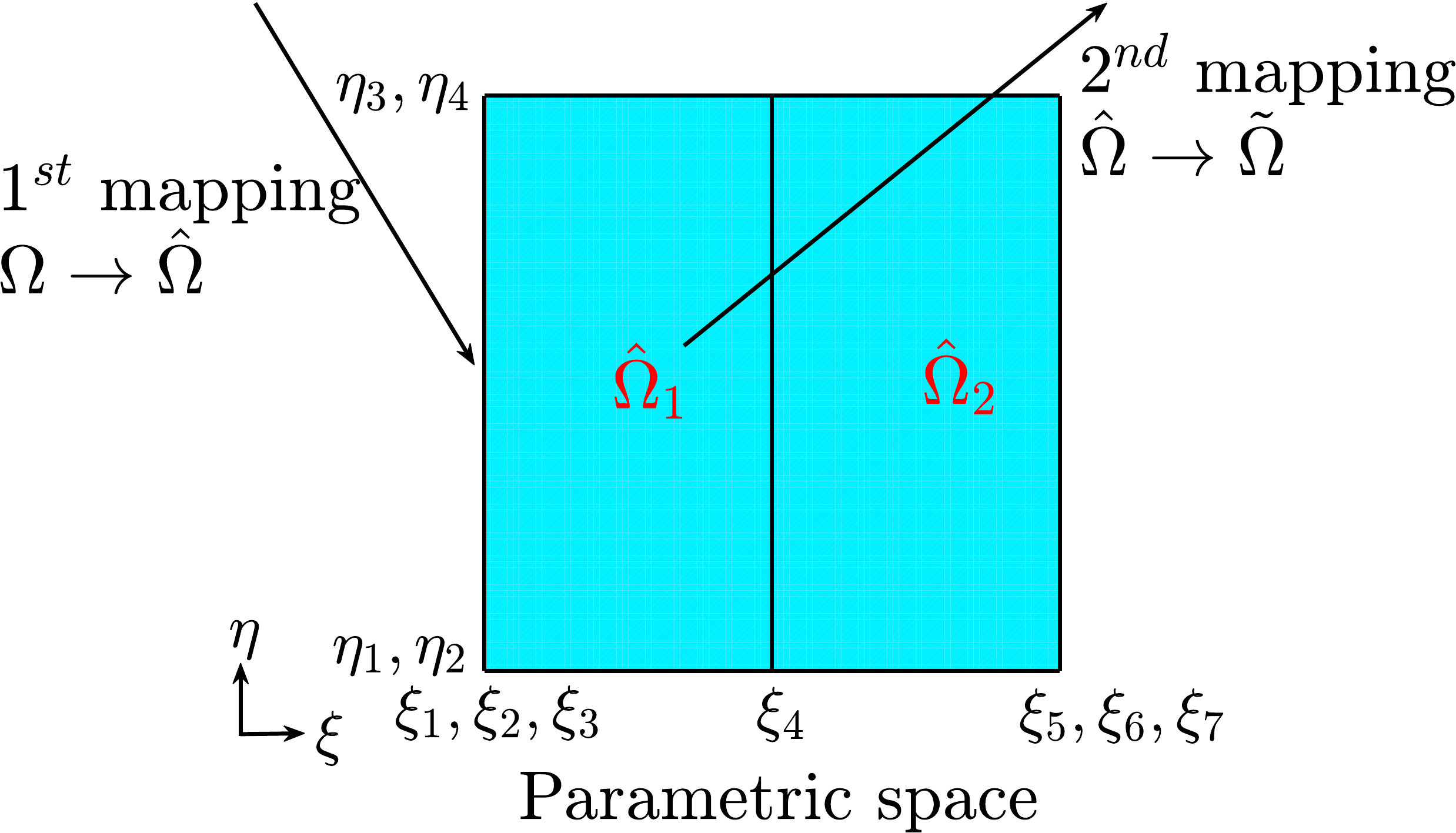}
\end{subfigure}
\caption{Mapping and different spaces involved in IGA framework}
\label{mapping}
\end{figure}

The integrals, associated with the evaluation of element stiffness matrices or the force vector, are solved using the Gauss-Legendre quadrature rules. However, due to the involvement of three mapping spaces, IGA necessitates an additional Jacobian. If $\Omega$ represents the physical space of the problem domain, $\hat\Omega$ denotes the parametric space where the NURBS basis functions are defined, and $\tilde{\Omega}$ is the master or parent space then, the mapping involved in the IGA can be illustrated using Figure~\ref{mapping}.

% ==============================================================================
%  Section- Force vector
% ==============================================================================
\subsection{Force vector and boundary conditions}

The approach for evaluation of the force vector and imposition of boundary conditions is identical to the conventional single field IGA. The traction boundary conditions ($\bt = \bar\bt$ on $\Gamma_t$) needs no special treatment, however, if the problem demands the point load at certain location then, it is mandatory to have a control point at that location which can be possible if the continuity at that point reduces to $C^{-1}$. This can be achieved by increasing the multiplicity of a knot in the knot vector. Furthermore, the procedure involved to incorporate homogeneous or non-homogeneous Dirichlet boundary conditions is indeed the same as single-field formulation. The numerical examples presented hereby deals with the homogeneous Dirichlet boundary condition i.e.\ $\bu = 0 $ on $\Gamma_u$, which can be incorporated by regulating the corresponding control variables as zero. For instance, if a particular boundary is fixed then, the displacement variable associated with the control points responsible to design that boundary are considered as zero even though the control points may or may not lie specifically on that boundary.  

% ==============================================================================
% Section- Refinement
% ==============================================================================
\subsection{Refinement}

Complimentary to conventional FEA, IGA incorporates \textit{h} and \textit{p}-refinement which can be achieved by knot insertion into the original knot vector and degree elevation of NURBS basis functions respectively. Moreover, IGA adapted an another approach by which the refinement will hold higher continuity with fewer basis functions namely \textit{k}-refinement. However, contrary to FEA, the stated refinement techniques can be practiced without further interaction with CAD model once the initial geometry is created. The same is also applicable in case of two-field hybrid IGA. The only difference is the choice of $\bP$ matrix, which will not change for a particular element during the \textit{h}-refinement, however, degree elevation techniques, either \textit{p} or \textit{k}-refinement, required to change the $\bP$ matrix corresponding to the resulting degree of NURBS basis function.

% ==============================================================================
% Section- Post-processing
% ==============================================================================
\subsection{Post-processing}

In order to visualize the deformed problem domain and CP mesh, the resultant displacement vector $\tilde{\bu}$ is added to the CP co-ordinates as follows,
\begin{align}
[\boldsymbol{CP}]_{new} = [\boldsymbol{CP}]_{old}+\tilde{\bu}
\end{align}
where $[\boldsymbol{CP}]_{new}$ is the new set of control points, $[\boldsymbol{CP}]_{old}$ are old control point values and $\tilde{\bu}$ is the displacement vector. After evaluation of updated control points, $[\boldsymbol{CP}]_{new}$ and original knot vectors ($\boldsymbol{\Xi}$ and $\bH$) are used to represent the deformed geometry and respective control point mesh \cite{Agrawal2018}.
% ------------------------------------------------------------------------------

% ==============================================================================
% Section- Numerical examples
% ==============================================================================
\section{Numerical examples}

A summary of short-hand notations used in subsequent examples is given as follows: The \textbf{IGA} and \textbf{H-IGA} denotes the conventional IGA and the proposed two-field hybrid IGA formulation. The extensions \textbf{d1}, \textbf{d2}, and \textbf{d3} states the use of linear, quadratic, and cubic basis functions along $\xi$ and $\eta$ direction, and \textbf{C0}, \textbf{C1}, and \textbf{C2} is the $C^0$, $C^1$, and $C^2$ inter-element continuity respectively. Furthermore, the notation \textbf{FEA} and \textbf{H-FEA} represents the conventional FEA and hybrid FEA formulations with extensions \textbf{Q4} and \textbf{Q9} as four and nine node quadrilateral elements respectively.

\subsection{Straight cantilever beam}

\begin{figure}[pos=!ht]
\centering
\includegraphics[width=.45\columnwidth]{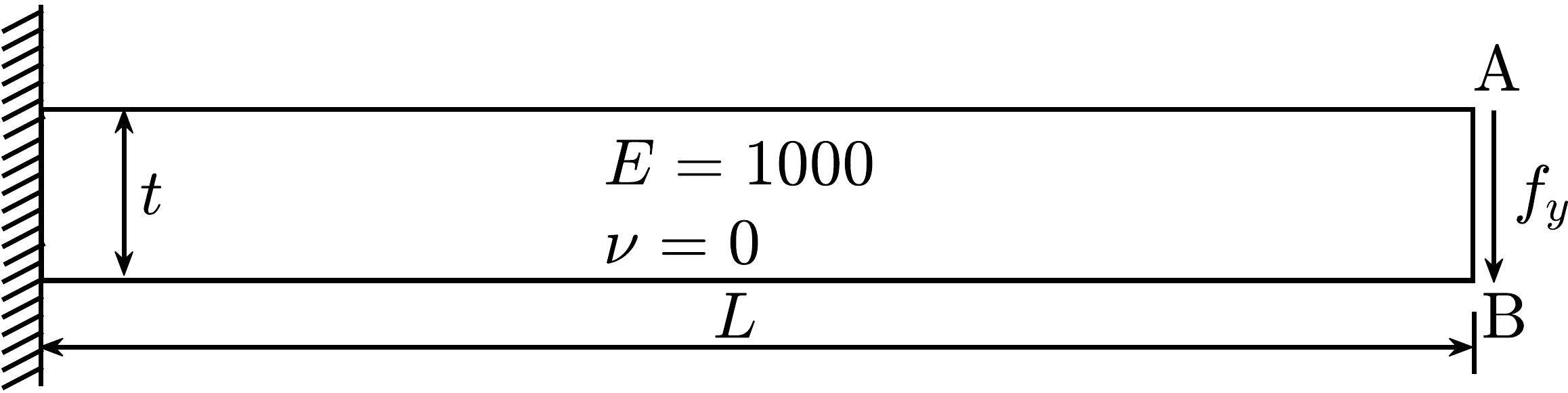}
\caption{A straight cantilever beam problem, material data, and boundary conditions}
\label{Straight_cantilever_beam_problem_definition}
\begin{subfigure}{0.45\columnwidth}
\centering
\includegraphics[width=1\columnwidth]{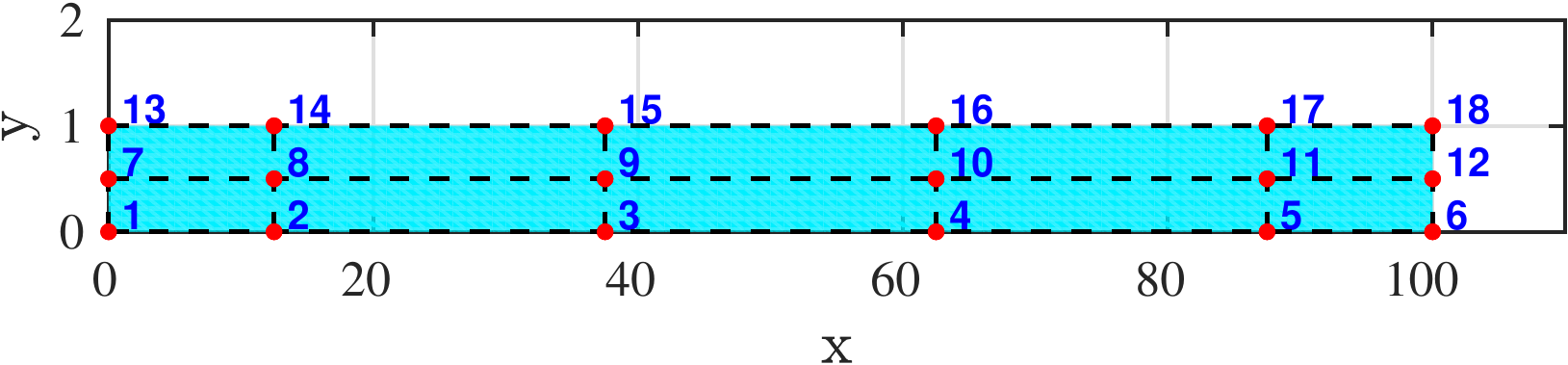}
\caption{CP mesh}
\label{Control_point_mesh_straight_beam}
\end{subfigure}
\begin{subfigure}{0.45\columnwidth}
\centering
\includegraphics[width=1\columnwidth]{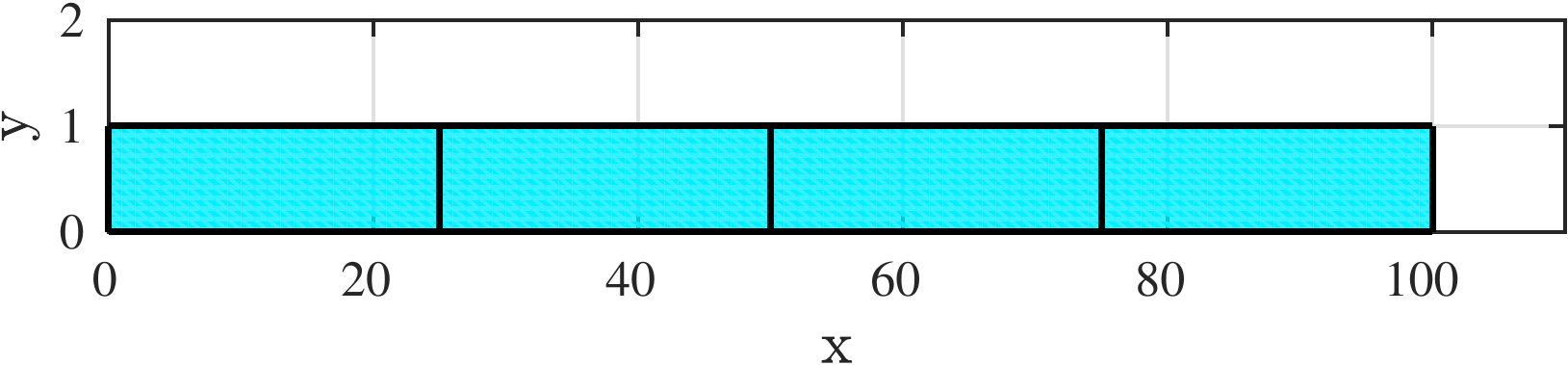}
\caption{Discretized problem domain}
\label{Discretized_problem_domain}
\end{subfigure} 
\caption{Geometric description of a straight cantilever beam problem ($L/t = 100$) for four NURBS elements with quadratic basis along $\xi$ and $\eta$ direction}\label{geometric_description_straight_cantilever_beam_problem}
\begin{subfigure}{0.45\columnwidth}
\centering
\includegraphics[width=1\columnwidth]{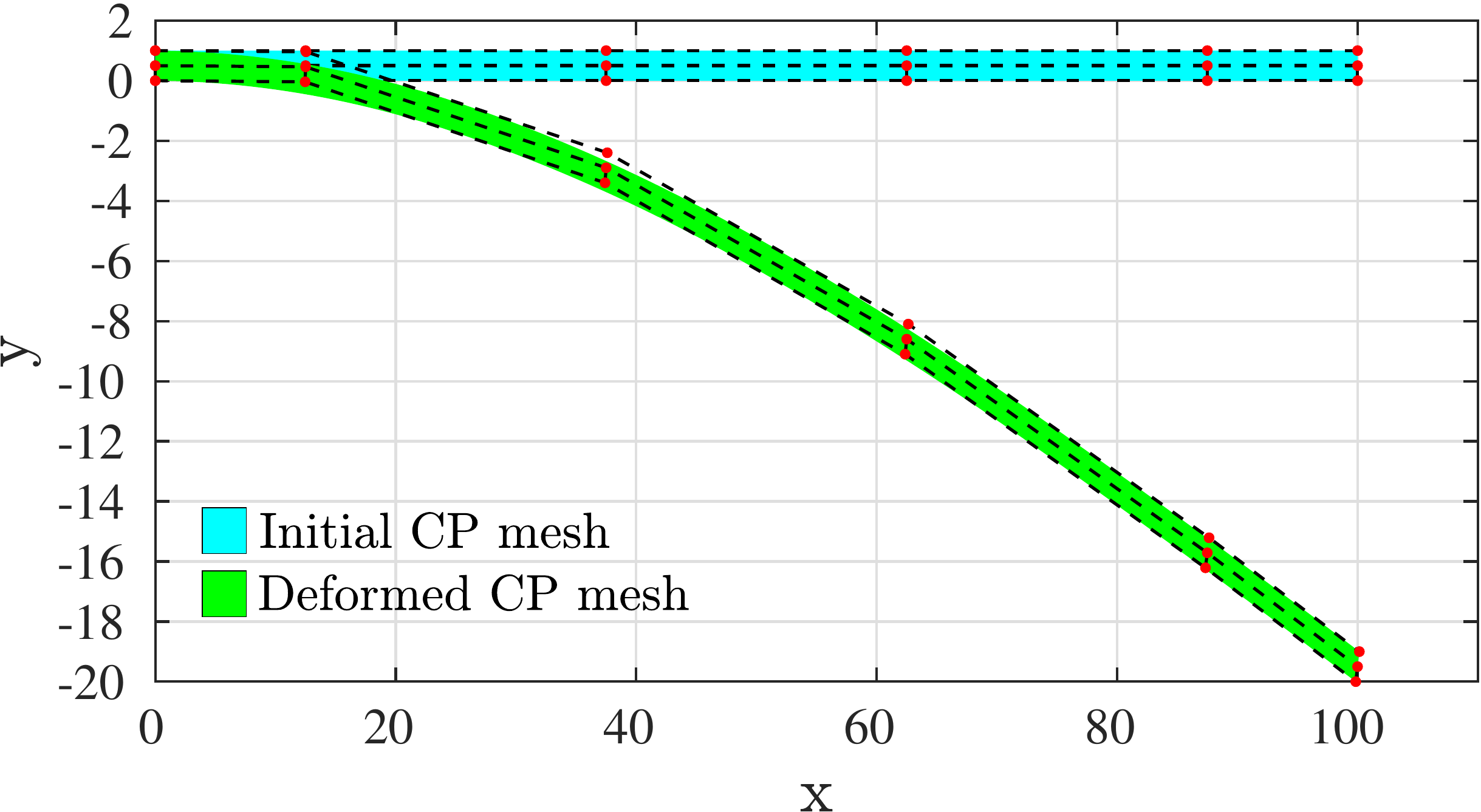}
\caption{Deformed CP mesh}
\label{Deformed_control_point_mesh_straight_beam}
\end{subfigure}
\begin{subfigure}{0.45\columnwidth}
\centering
\includegraphics[width=1\columnwidth]{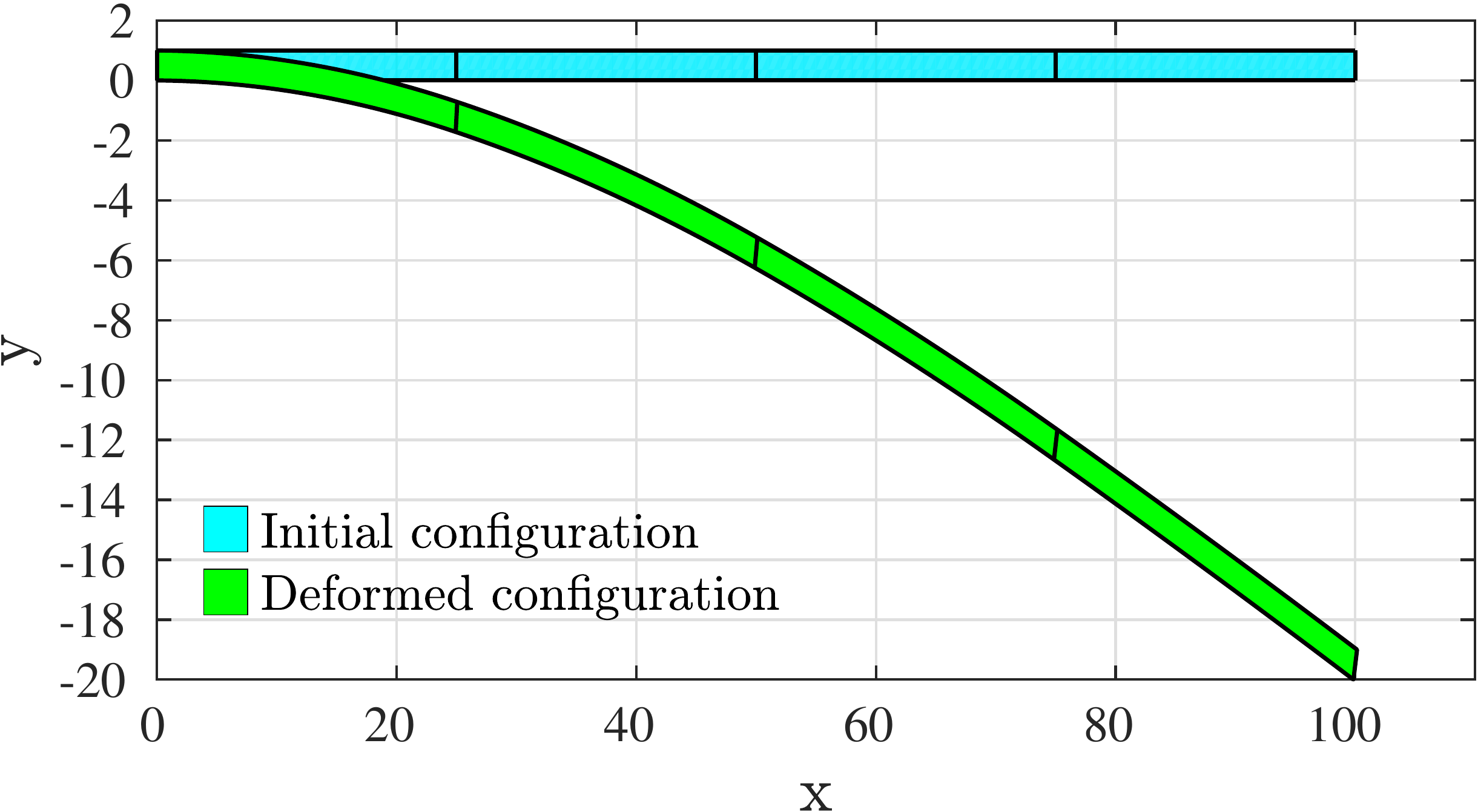}
\caption{Deformed discretized problem domain}
\label{Deformed_discretized_problem_domain}
\end{subfigure}
\newline
\begin{subfigure}{\columnwidth}
\centering
\includegraphics[trim={4cm 8cm 6cm 6cm},clip,width=0.5\columnwidth]{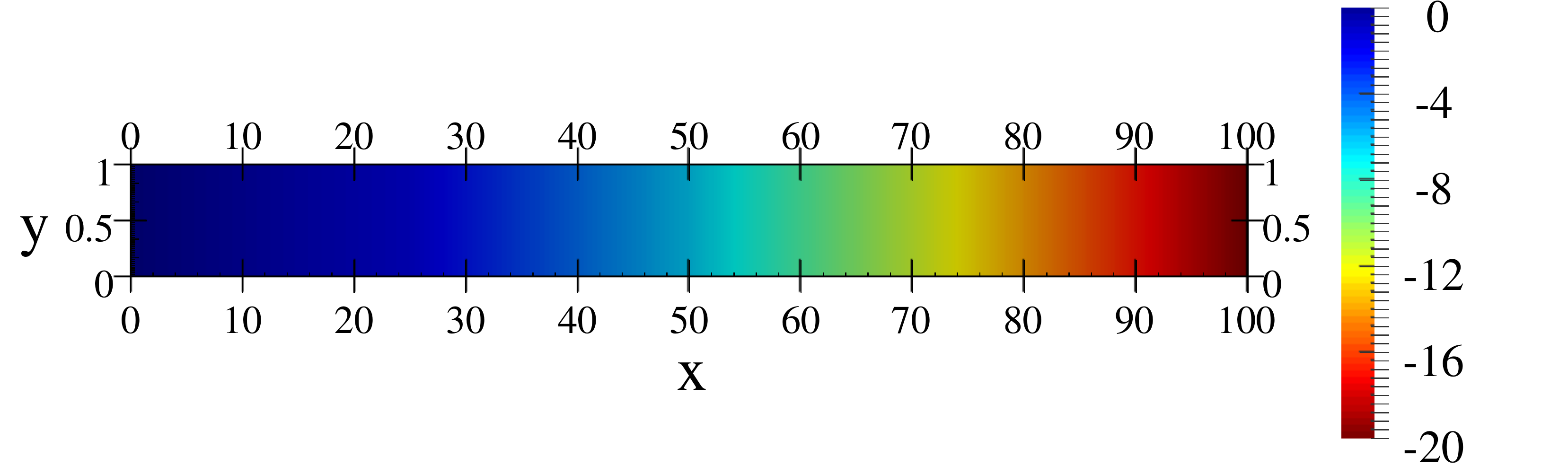}
\caption{$u_y^{H-IGA}$}
\label{contour_plot_straight_beam_u_y}
\end{subfigure}
\caption{(a-b) Deformed geometric description, and (c) contour plot for vertical displacement of a straight cantilever beam problem ($L/t = 100$) for four NURBS elements with quadratic basis along $\xi$ and $\eta$ direction}\label{comprehensive_results_straight_beam}
\end{figure}

For the first numerical example, a linear elastic behavior of a two-dimensional cantilever beam of length $L$ and thickness $t$ subjected to vertical load $f_y$ is investigated. The problem is kept simple to test the reliability of the proposed formulation under the influence of the shear locking, see Figure~\ref{Straight_cantilever_beam_problem_definition} \cite{book_FEM}.

The CP and the respective weights to model the coarsest possible mesh representing the exact geometry is provided in the appendix (Table~\ref{Table_CP_rectangular_beam}). Once the initial mesh is generated, the sequence of meshes is constructed using the $h$ and $k$-refinement. One such mesh of $4 \times 1$ quadratic NURBS elements, for the slenderness ratio ($L/t$) 100, is illustrated in Figure~\ref{geometric_description_straight_cantilever_beam_problem}, which highlights the required number of control points, the control point mesh, and the respective discretization of a domain into quadratic NURBS based elements. In the interest of embracing the proposed formulation, the elaborated results are presented in the Figure~\ref{comprehensive_results_straight_beam} which gives an idea about the deformed configuration of control point mesh and the discretized domain along with the contour plot for vertical displacement for the stated mesh.

\begin{figure}[pos=ht]
\begin{subfigure}{0.47\columnwidth}
\centering
\includegraphics[width=1\columnwidth]{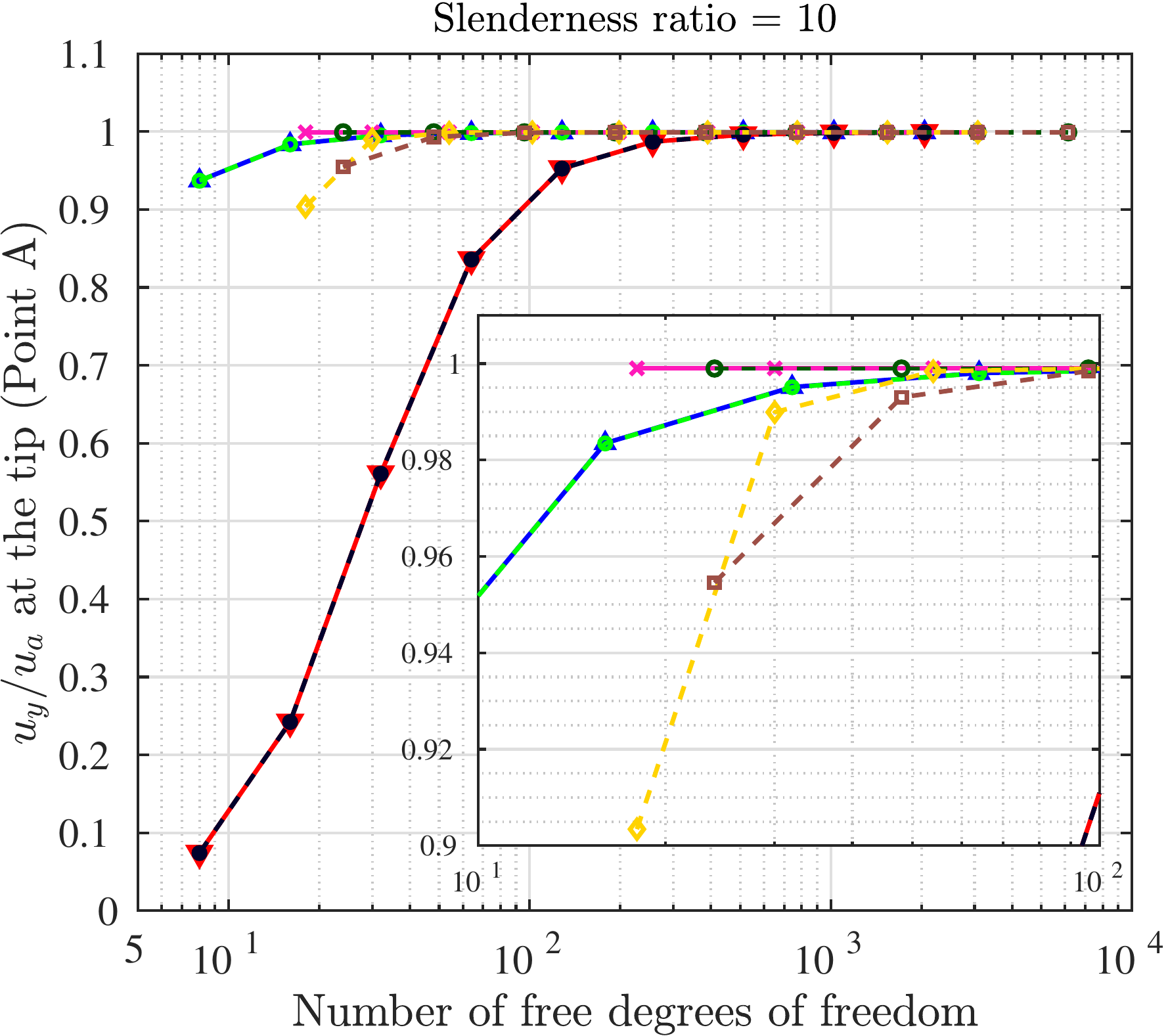}
\caption{$L/t = 10$}
\label{Normalized_disp_straight_beam_10}
\end{subfigure}
\begin{subfigure}{0.47\columnwidth}
\centering
\includegraphics[width=1\columnwidth]{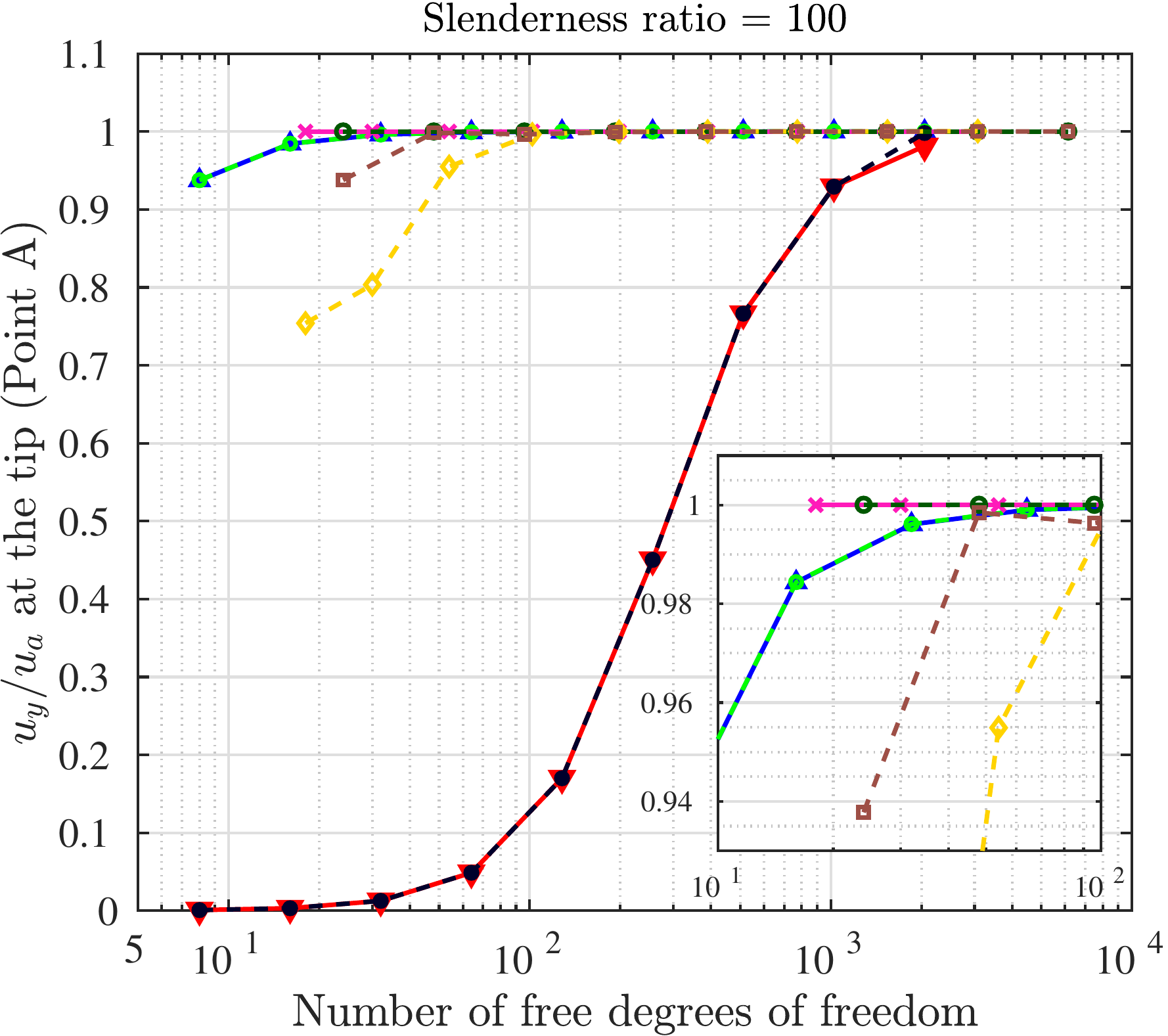}
\caption{$L/t = 100$}
\label{Normalized_disp_straight_beam_100}
\end{subfigure}
\newline
\begin{subfigure}{1\columnwidth}
\centering
\includegraphics[width=.65\columnwidth]{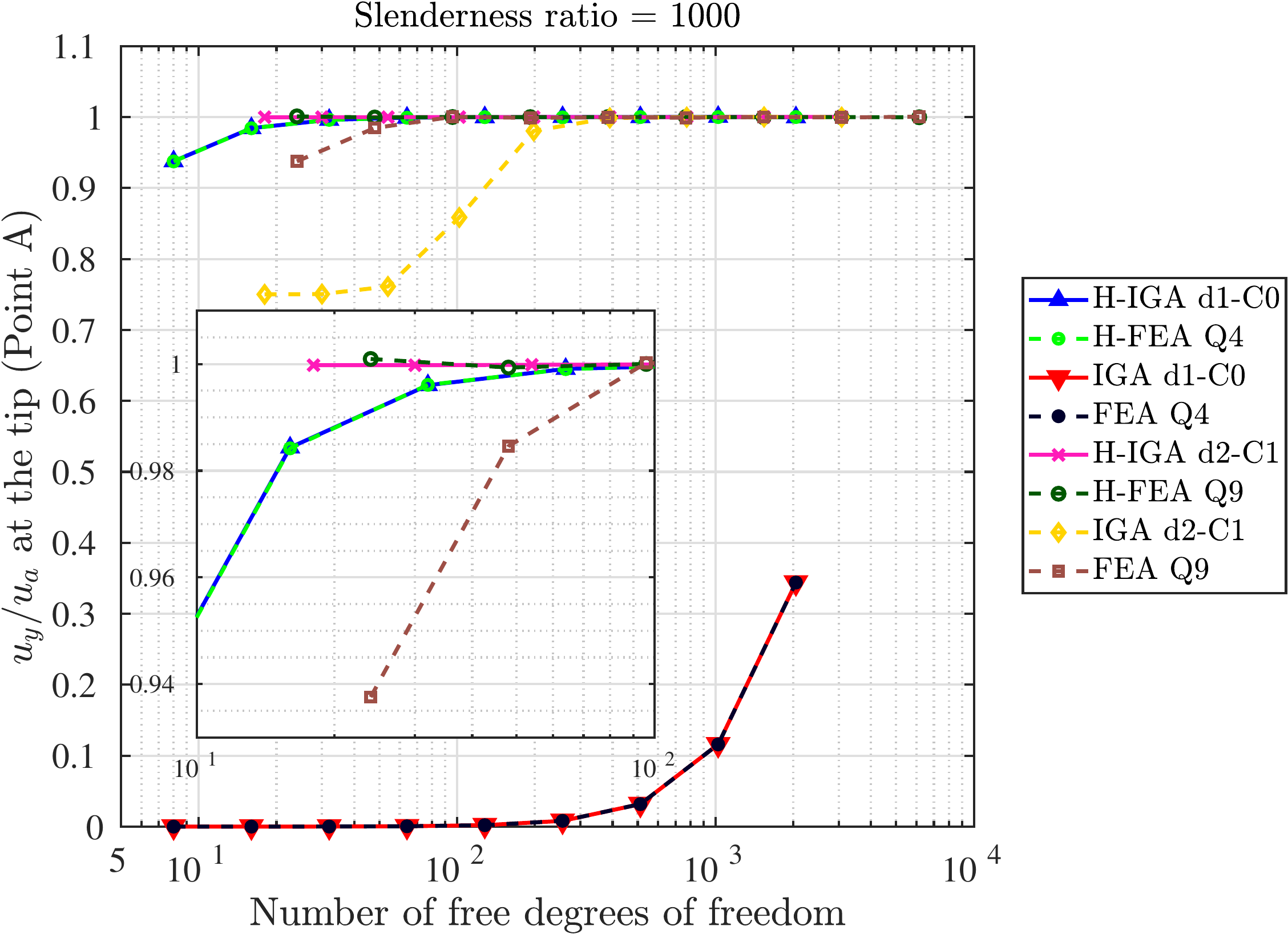}
\caption{$L/t = 1000$}
\label{Normalized_disp_straight_beam_1000}
\end{subfigure}
\caption{Normalized vertical displacement at point `A' for straight rectangular cantilever beam problem for three $L/t$ ratios}\label{Normalized_displacement_straight_beam_problem}
\end{figure}

The problem is studied by employing the FEA and IGA formulations for three different slenderness ratios (10, 100, and 1000) in order to gradually introduce the shear locking effect into the problem domain. The problem data considered for the three cases is given as follows,
\begin{enumerate}
\item $L/t = 10 $, $L = 100 $, $t = 10$, $f_y=4.97018$
\item $L/t = 100 $, $L = 100 $, $t = 1$, $f_y=4.9997\times 10^{-3}$
\item $L/t = 1000 $, $L = 100 $, $t = 0.1$, $f_y=4.9999\times 10^{-6}$
\end{enumerate}
For all the three cases, the analytical solution for the vertical displacement ($u_y$) at point `A' is 20.

The vertical displacement at point `A' is numerically evaluated for all the three cases and the corresponding results are presented in Figure~\ref{Normalized_displacement_straight_beam_problem}. For the lower value of slenderness ratio ($L/t=10$, Figure~\ref{Normalized_disp_straight_beam_10}), the locking effect is significantly low whether it is IGA or FEA formulation. However, the proposed hybrid IGA out-performs the conventional formulation with coarse mesh accuracy. For instance, with only two quadratic elements (active dof = 18), the hybrid IGA formulation is capable of providing the results which are in close approximation with the analytical solution. Whereas, with the conventional IGA, further refinement is needed to achieve a similar accuracy.

The effect of the shear locking can be distinctively observed in conventional IGA while using the lower degree basis functions with a higher slenderness ratio of the problem domain. As illustrated in Figure~\ref{Normalized_disp_straight_beam_100} and~\ref{Normalized_disp_straight_beam_1000}, conventional IGA with quadratic basis functions locks with higher value of $L/t$ and a significant refinement is needed to alleviate the anomaly. On the other hand, the hybrid IGA performs convincingly well in all the conditions by alleviating the locking and providing the superior coarse mesh accuracy. Furthermore, the use of higher degree NURBS, either with conventional or hybrid IGA, significantly reduced the locking. The results for cubic NURBS interpolations are not presented for this particular problem as the formulation, either it is conventional IGA or hybrid IGA, converge to the exact solution with minimum number of active degrees of freedom itself. For the point of interest, the results obtained by linear NURBS elements in IGA formulation and Q4 elements in FEA are identical due to the fact that the NURBS basis functions of degree 1 with weights as 1, will reduce to the conventional Lagrangian basis functions used for Q4 elements. 

\subsection{Curved beam}

\begin{figure}[pos=t]
\begin{subfigure}{0.32\columnwidth}
\centering
\includegraphics[width=0.7\columnwidth]{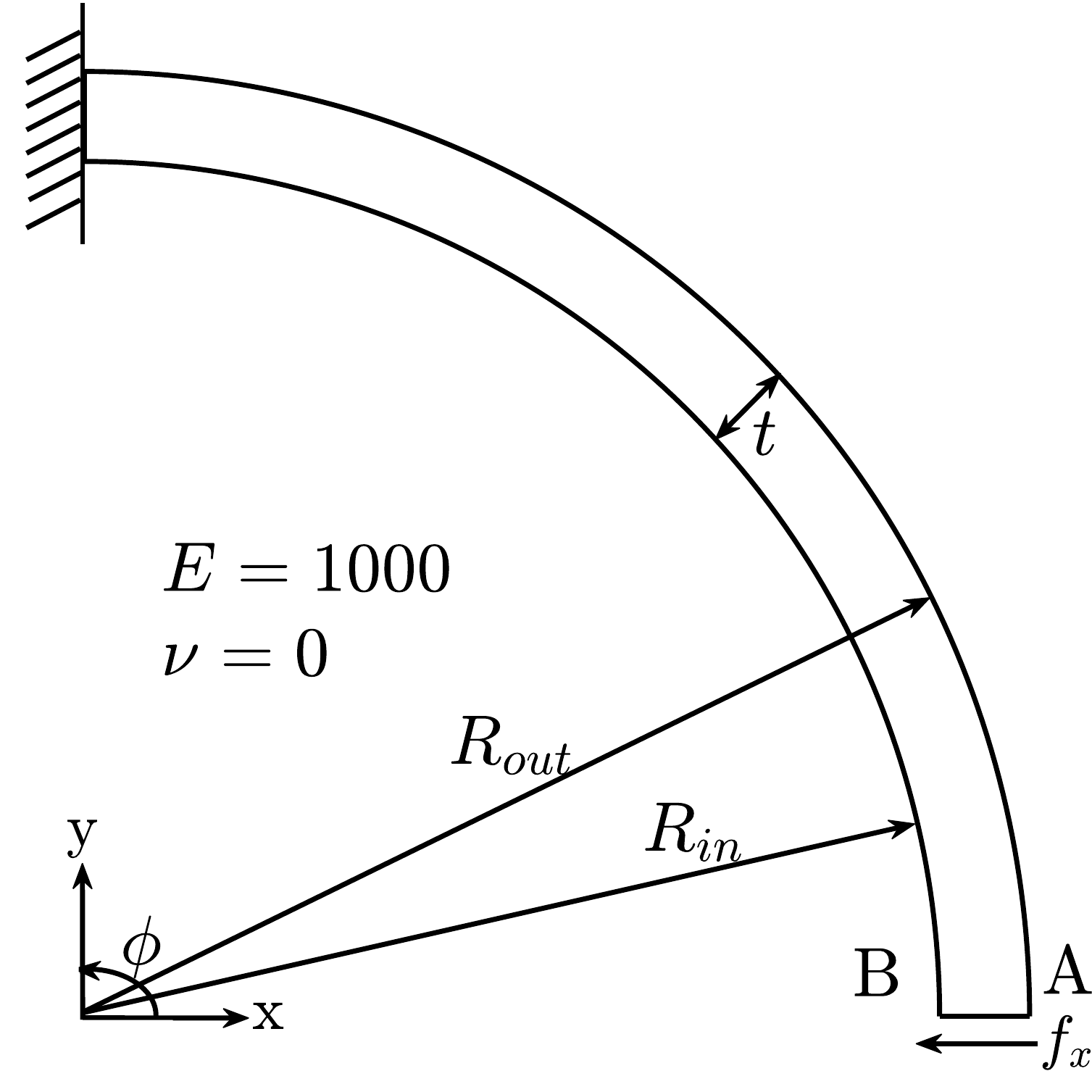}
\caption{A curved beam, material data, and boundary conditions}
\label{curved_beam_problem_definition}
\end{subfigure}
\begin{subfigure}{0.32\columnwidth}
\centering
\includegraphics[width=0.8\columnwidth]{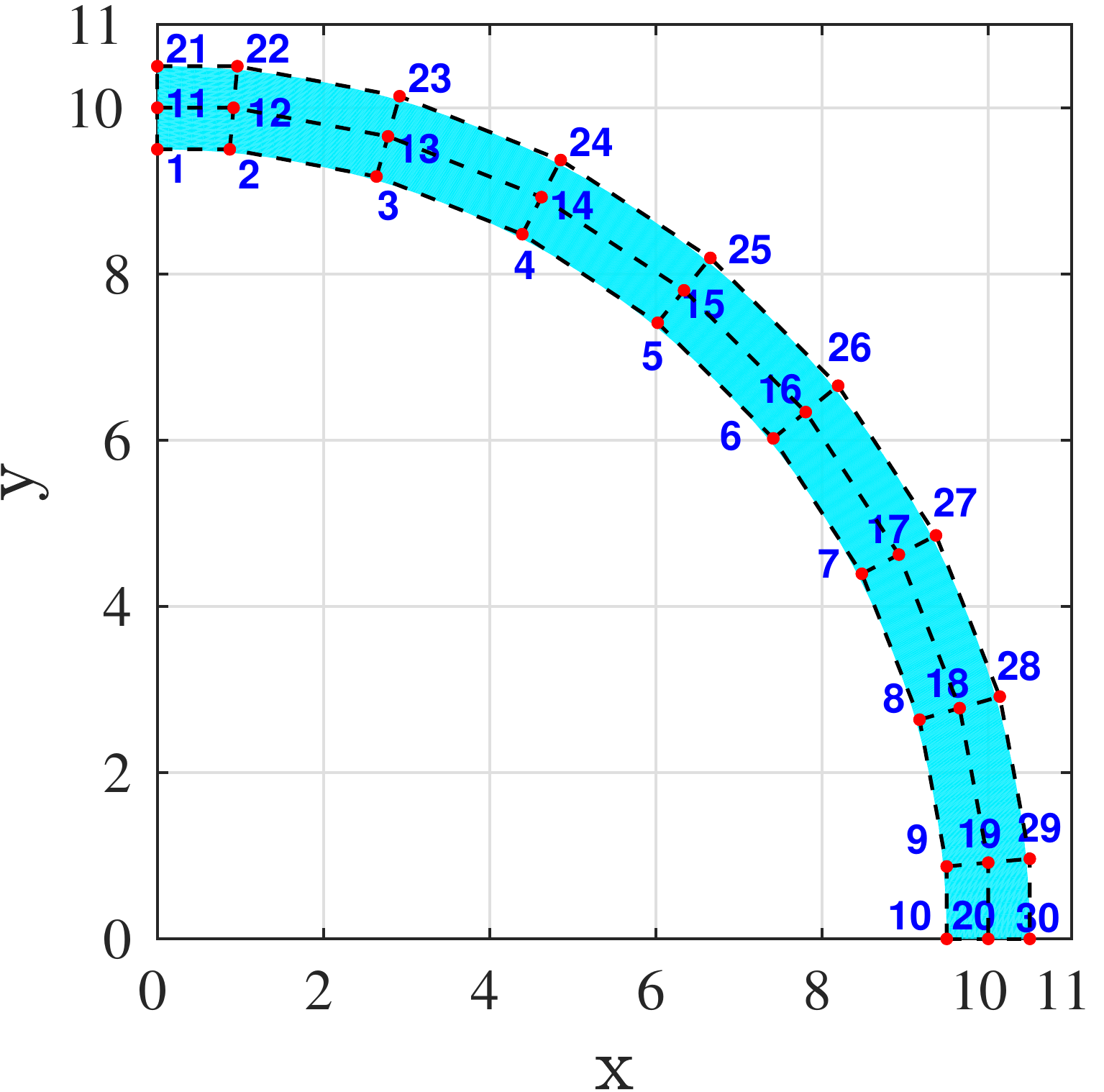}
\caption{CP mesh}
\label{control_point_mesh_curved_beam_problem}
\end{subfigure}
\begin{subfigure}{0.32\columnwidth}
\centering
\includegraphics[width=0.8\columnwidth]{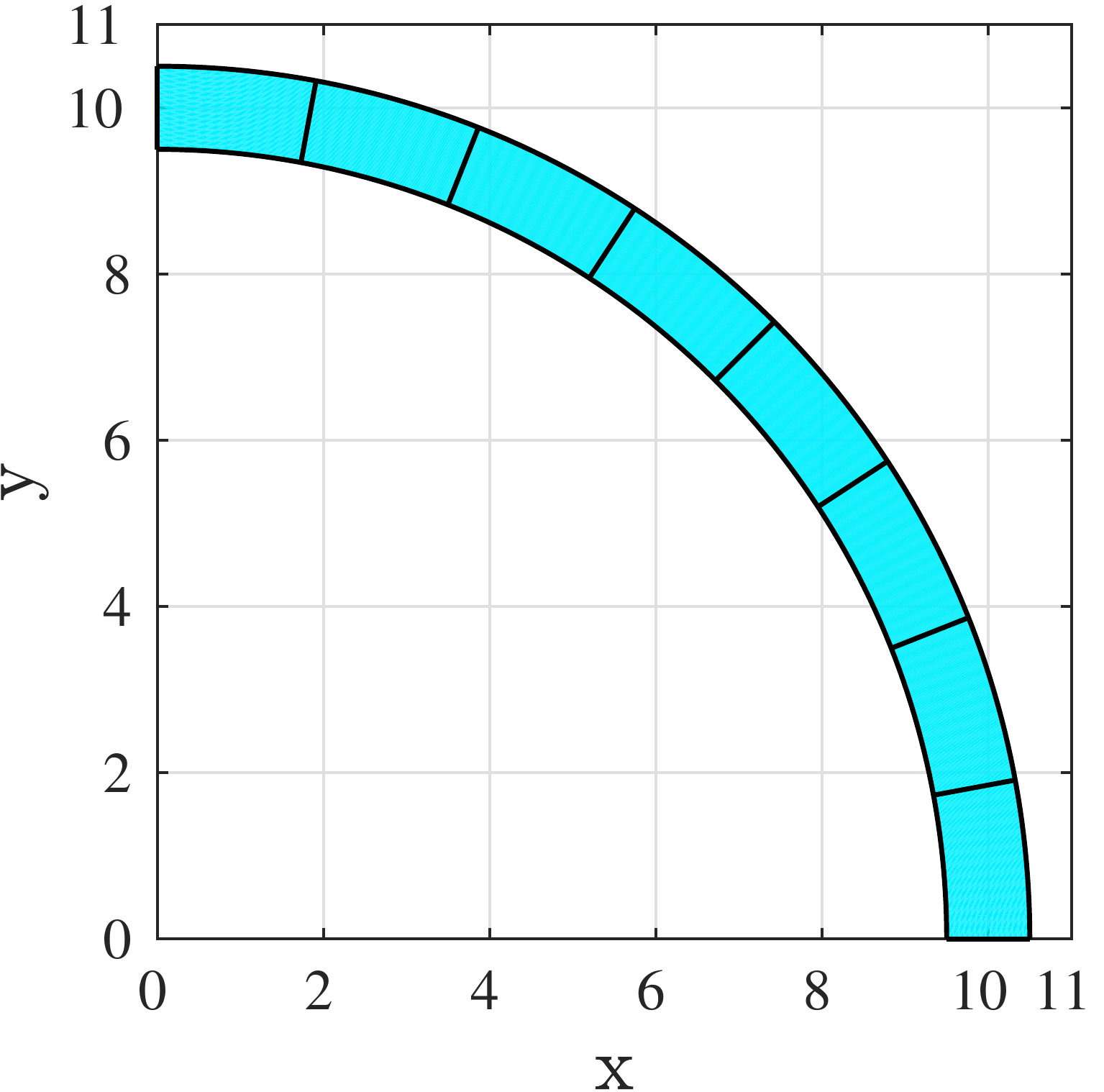}
\caption{Discretized problem domain}
\label{element_discretization_curved_beam_problem}
\end{subfigure} 
\caption{(a) The problem definition and (b-c) geometric description of an curved cantilever beam problem ($R/t = 10$) for 8$\times$1 NURBS elements with quadratic basis along $\xi$ and $\eta$ direction} \label{geometric_description_of_curved_beam_problem}
\end{figure}

\begin{figure}[pos=b]
\begin{subfigure}{0.29\columnwidth}
\centering
\includegraphics[width=0.9\columnwidth]{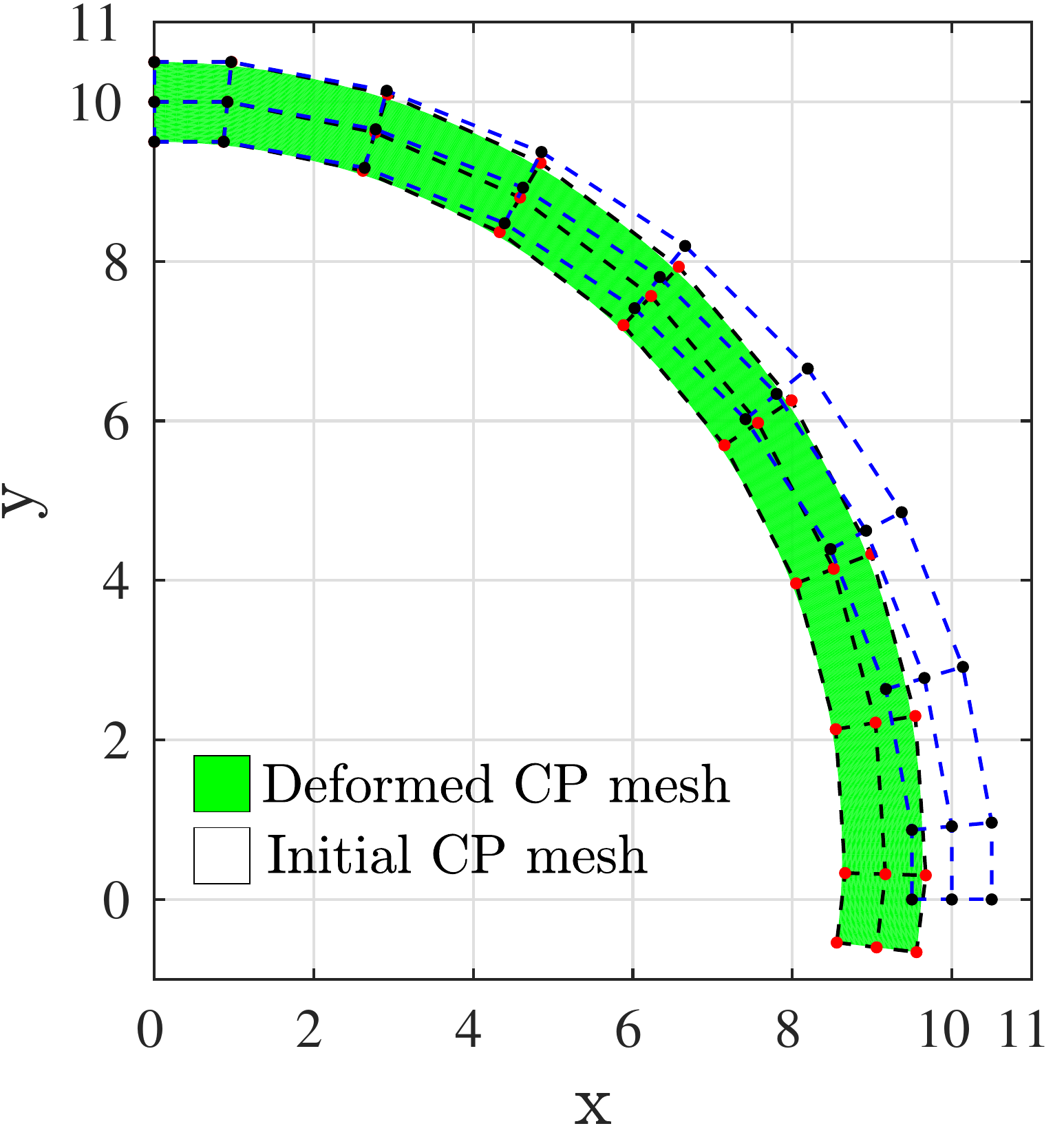}
\captionsetup{width=3cm}
\caption{Deformed CP mesh}
\label{Deformed_CP_mesh_curved_beam}
\end{subfigure}
\begin{subfigure}{0.29\columnwidth}
\centering
\includegraphics[width=0.9\columnwidth]{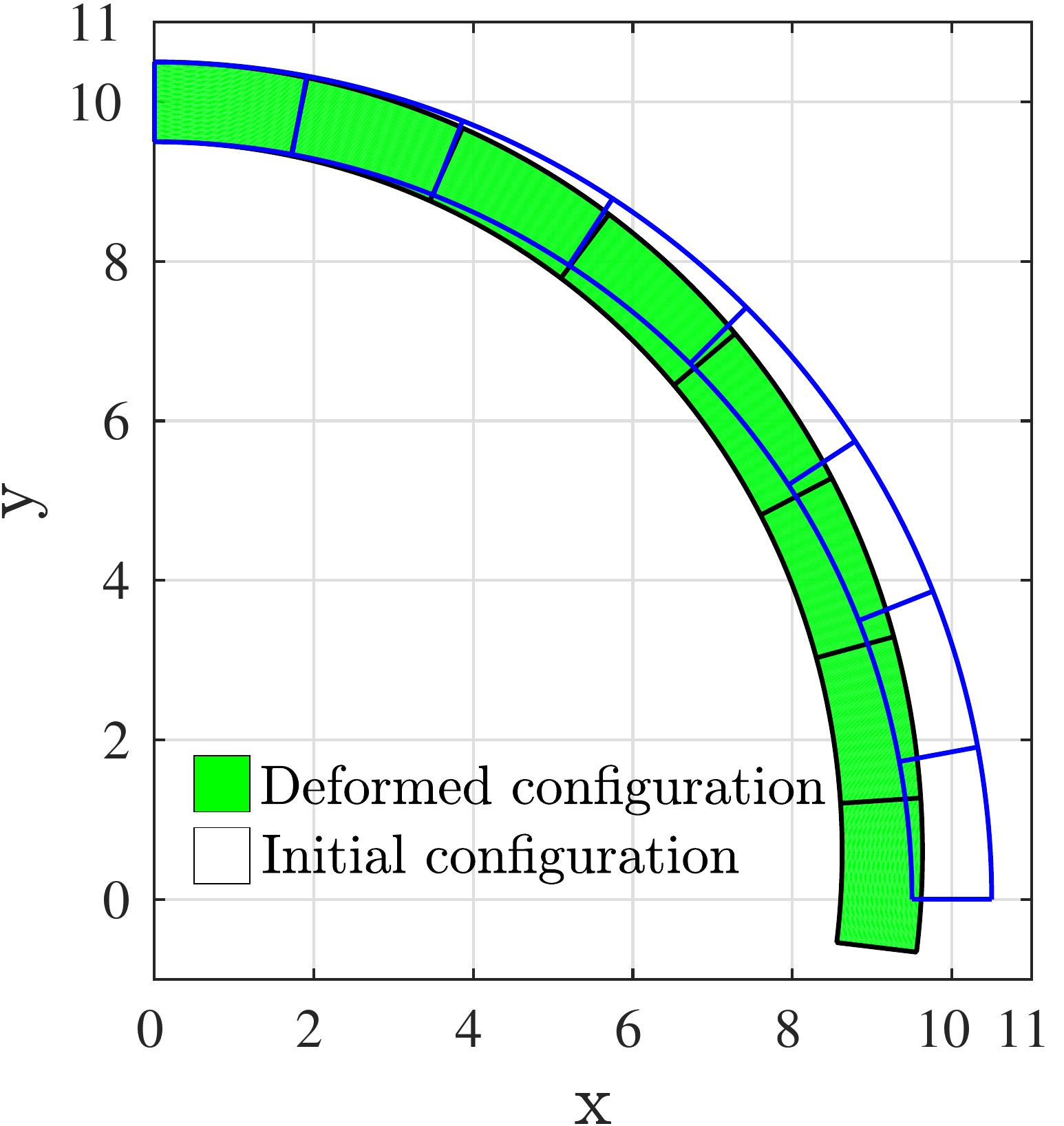}
\captionsetup{width=5.5cm}
\caption{Deformed discretized domain}
\label{deformed_element_discretization_curved_beam}
\end{subfigure} 
\begin{subfigure}{0.38\columnwidth} 
\centering
\includegraphics[trim={9cm 1cm 6.5cm 2cm},clip,width=1\columnwidth]{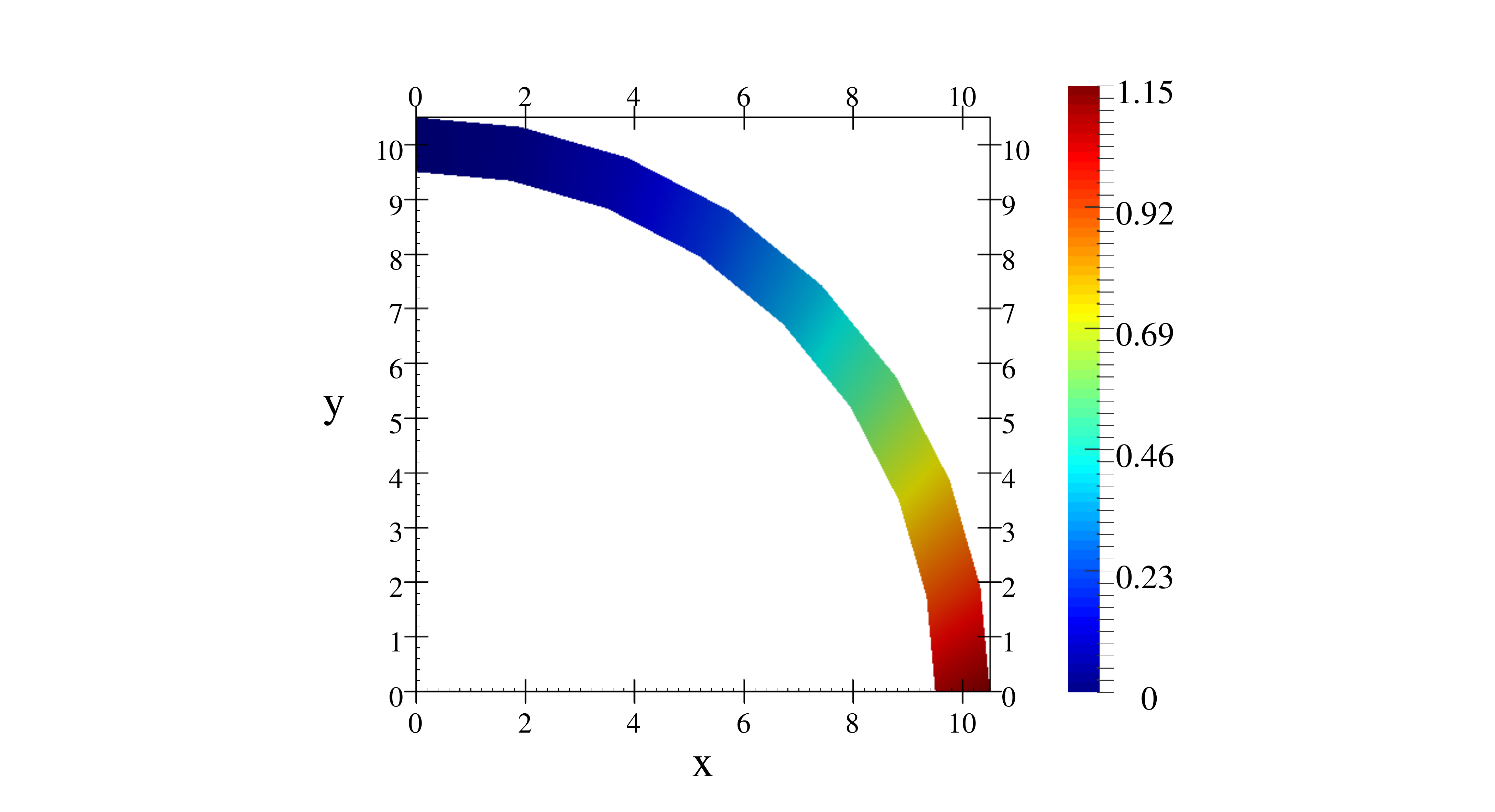}
\captionsetup{width=4cm}
\caption{$u^{H-IGA}$}
\label{Contour_plot_total_displacement_curved_beam}
\end{subfigure}
\caption{(a-b) Deformed geometric description, and (c) contour plot of magnitude of total displacement for a curved beam problem for 8$\times$1 NURBS elements with quadratic basis along $\xi$ and $\eta$ direction as illustrated in Figure~\ref{control_point_mesh_curved_beam_problem}-\ref{element_discretization_curved_beam_problem} using hybrid IGA formulation}\label{comprehensive_results_curved_beam_nu_04999}
\end{figure}

In the present example, a linear elastic behavior of a two-dimensional curved cantilever beam is investigated. The problem is composed of a curved beam subjected to horizontal load on one end and fixed on the other end. The problem setup and the boundary conditions are illustrated in Figure~\ref{curved_beam_problem_definition} where $R_{in}$ and $R_{out}$ are the inner and outer radii measured from the origin, $R$ is the mean radius, and $t$ is the thickness of the beam. $f_x$ is the magnitude of the load at the free end such that radial displacement at the tip (Point `A') is evaluated as 0.942 \cite{Caseiro2014}. It is calculated as, $f_x = 0.1t^3$. $\nu$ is Poisson's ratio, and $E$ is Young's modulus.

To exactly represent the circular edges of the problem domain; the minimum requirement is to incorporate the quadratic NURBS basis functions along the curvature. The required CP along with the respective weights to model the coarsest possible mesh representing the exact geometry is provided in the appendix (Table~\ref{control_points_curved_beam}-\ref{weights_curved_beam}). Once the initial mesh is generated, the sequence of meshes is constructed using the $h$ and $k$-refinement. One such mesh of $8 \times 1$ quadratic NURBS elements, for the slenderness ratio ($R/t$) 10, is illustrated in Figure~\ref{control_point_mesh_curved_beam_problem}-\ref{element_discretization_curved_beam_problem}, which focuses on the number of control points involved, the control point mesh, and the respective element discretization of a domain. To achieve the sense of completeness, the extensive results for the stated mesh is presented in Figure~\ref{comprehensive_results_curved_beam_nu_04999} which elaborates the deformed configuration of the problem domain (Figure~\ref{Deformed_CP_mesh_curved_beam}-~\ref{deformed_element_discretization_curved_beam}) along with the contour plots for displacement field (Figure~\ref{Contour_plot_total_displacement_curved_beam}) obtained by incorporating the hybrid IGA formulation.

\begin{figure}[pos=ht]
\begin{subfigure}{0.45\columnwidth}
\centering
\includegraphics[width=1\columnwidth]{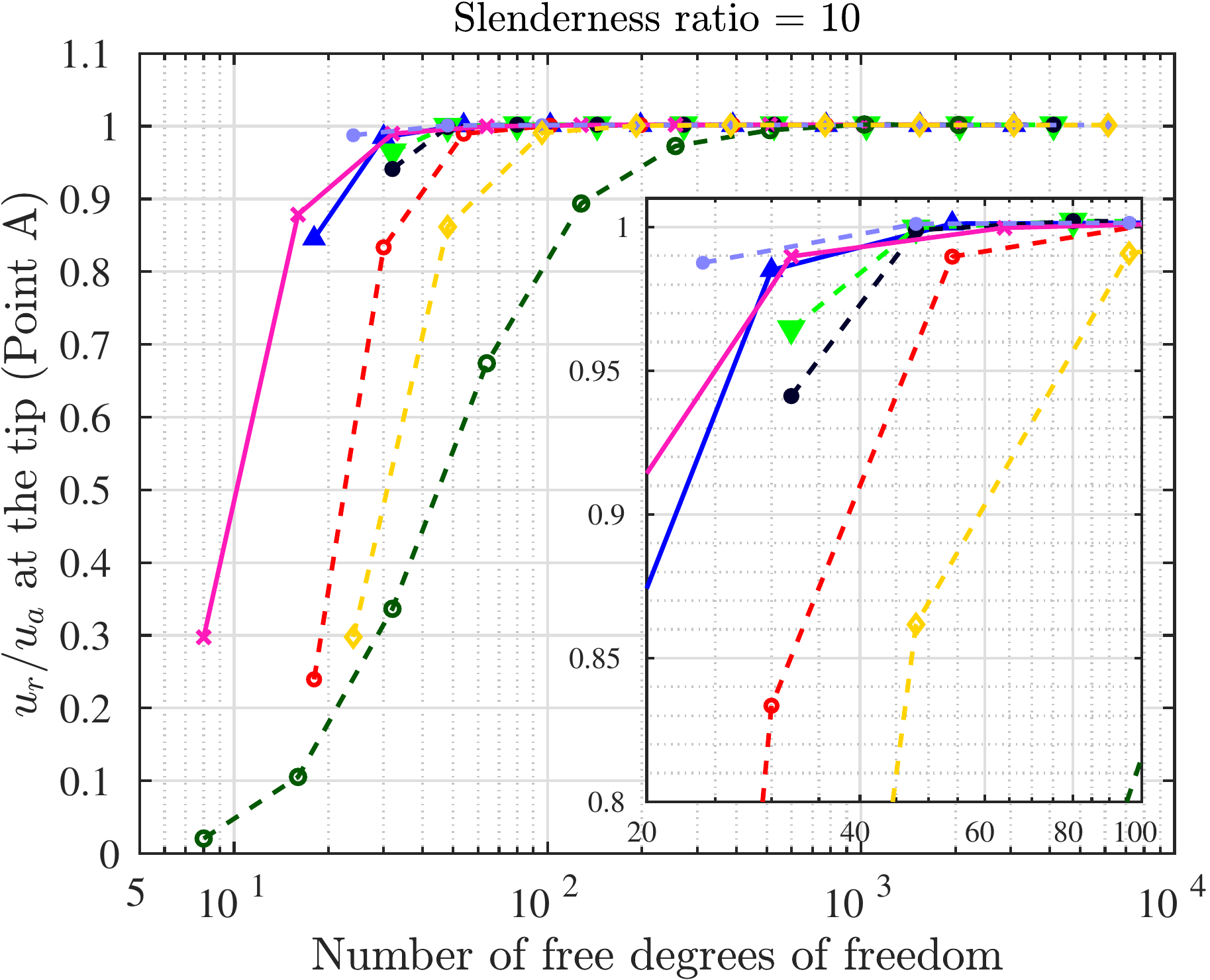}
\caption{$R/t = 10$}
\label{Normalized_displacement_curved_beam_s10}
\end{subfigure}
\begin{subfigure}{0.45\columnwidth}
\centering
\includegraphics[width=1\columnwidth]{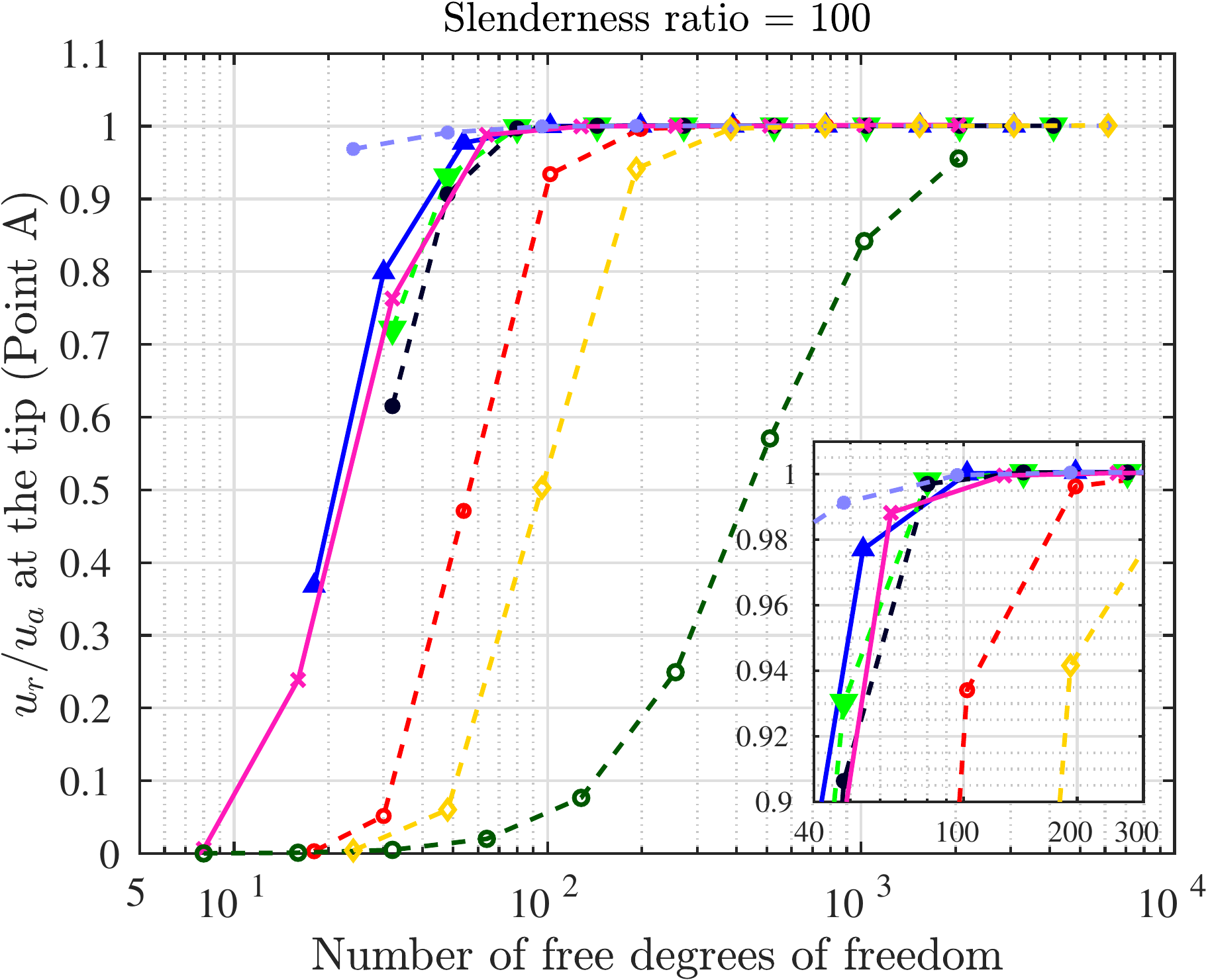}
\caption{$R/t = 100$}
\label{Normalized_displacement_curved_beam_s100}
\end{subfigure} 
\newline
\begin{subfigure}{1\columnwidth}
\centering
\includegraphics[width=.6\columnwidth]{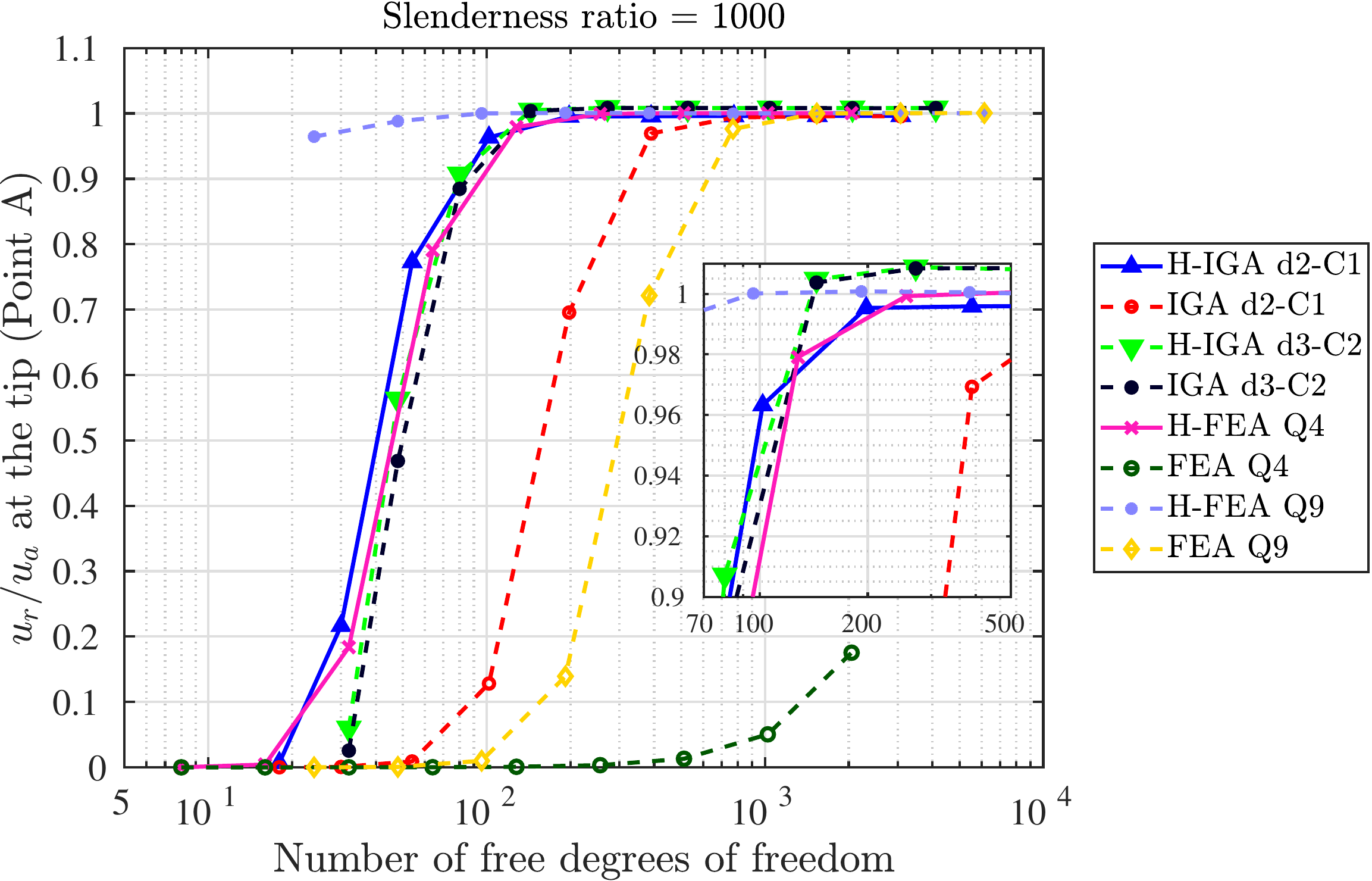}
\caption{$R/t = 1000$}
\label{Normalized_displacement_curved_beam_s1000}
\end{subfigure}
\caption{Normalized radial displacement at point `A' for a curved cantilever beam problem for three $R/t$ ratios}
\label{Normalized_radial_displacement_curved_beam_problem}
\end{figure}

The problem is solved using the FEA and IGA formulations for three different slenderness ratios (10, 100, and 1000) in order to gradually introduce the shear locking effect into the problem domain. The problem data considered for the three cases is given as follows,
\begin{enumerate}
\item $R/t = 10 $, $R_{in} = 9.5$, $R_{out} = 10.5$, $R = 10$, $t = 1$, $f_x=0.1$
\item $R/t = 100 $, $R_{in} = 9.95$, $R_{out} = 10.05$, $R = 10$, $t = 0.1$, $f_x=0.1\times 10^{-3}$
\item $R/t = 1000 $, $R_{in} = 9.995$, $R_{out} = 10.005$, $R = 10$, $t = 0.01$, $f_x=0.1\times 10^{-6}$
\end{enumerate}
For all the three cases, the analytical solution for the radial displacement ($u_a$) at point `A' is 0.942.

The radial displacement at point `A' is numerically evaluated by employing the different formulations and results are presented in Figure~\ref{Normalized_radial_displacement_curved_beam_problem}. For the lower value of slenderness ratio ($R/t=10$, Figure~\ref{Normalized_displacement_curved_beam_s10}), it can be seen that the locking effect is substantially low whether it is IGA or FEA formulation. However, the proposed hybrid IGA out-performs conventional formulation with coarse mesh accuracy. For instance, the hybrid IGA results for quadratic NURBS basis are in close approximation with the analytical solution even with merely 30 active degrees of freedom, whereas the results for conventional IGA even with the cubic basis are inferior for nearly the same degrees of freedom. 

As the slenderness ratio increases, the influence of the shear locking in the conventional IGA formulation can be observed distinctly while using the lower degree basis functions. As illustrated in Figure~\ref{Normalized_displacement_curved_beam_s100} and~\ref{Normalized_displacement_curved_beam_s1000}, conventional IGA with quadratic basis functions locks severely with higher value of $R/t$. On the other hand, the hybrid IGA performs convincingly well in all the conditions by alleviating the locking. Furthermore, the use of higher degree NURBS significantly reduced the locking, but hybrid results can be seen marginally better than the conventional formulation. From a comparative perspective, the IGA results, either conventional or hybrid, seem to be better than their FE counterparts.

\subsection{Cook's membrane problem}

\begin{figure}[pos=H]
\begin{subfigure}{0.33\columnwidth}
\centering
\includegraphics[width=0.9\columnwidth]{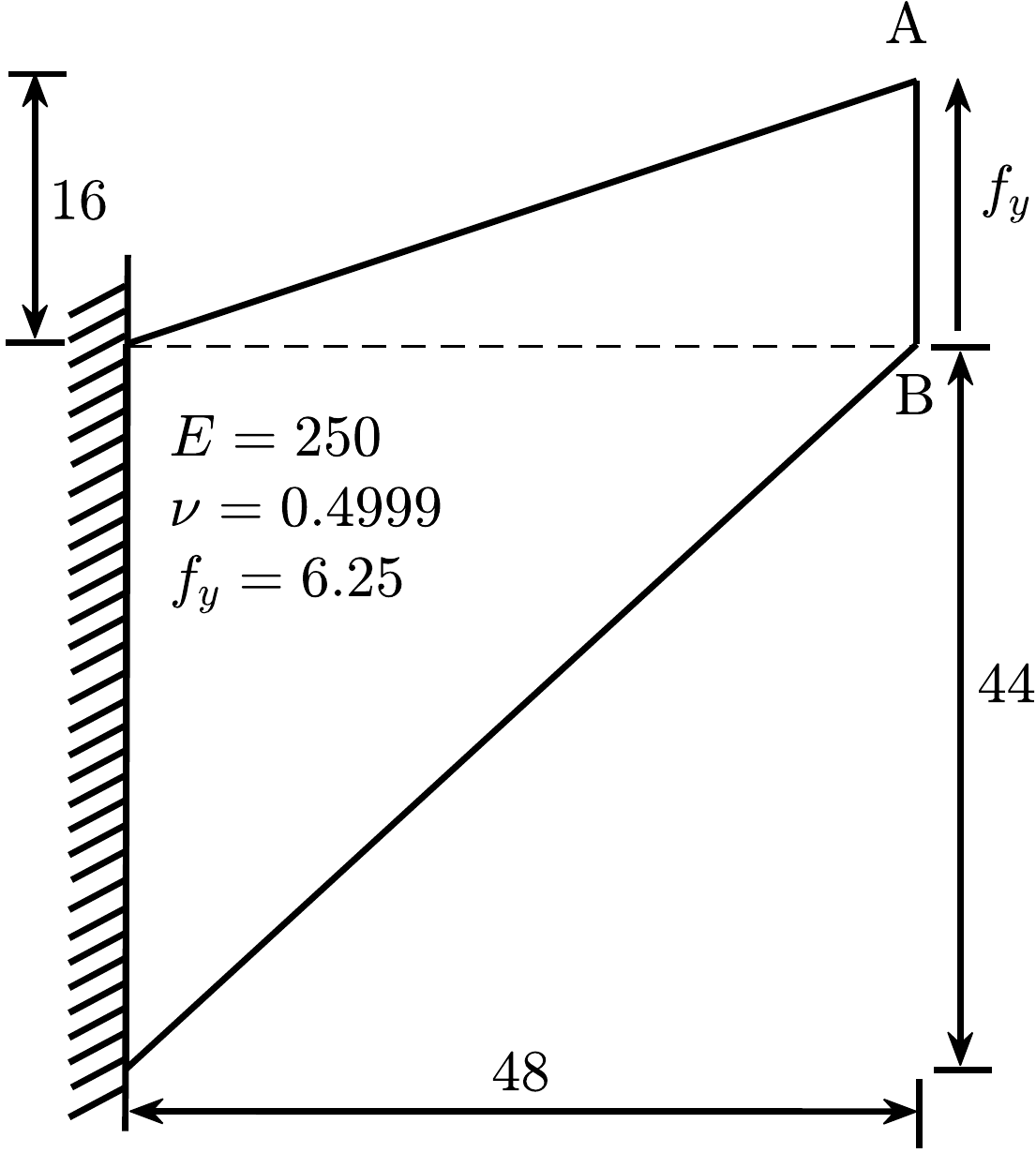}
\caption{The problem setup and boundary conditions}
\label{cooks_membrane_problem_definition}
\end{subfigure}
\begin{subfigure}{0.31\columnwidth}
\centering
\includegraphics[width=0.9\columnwidth]{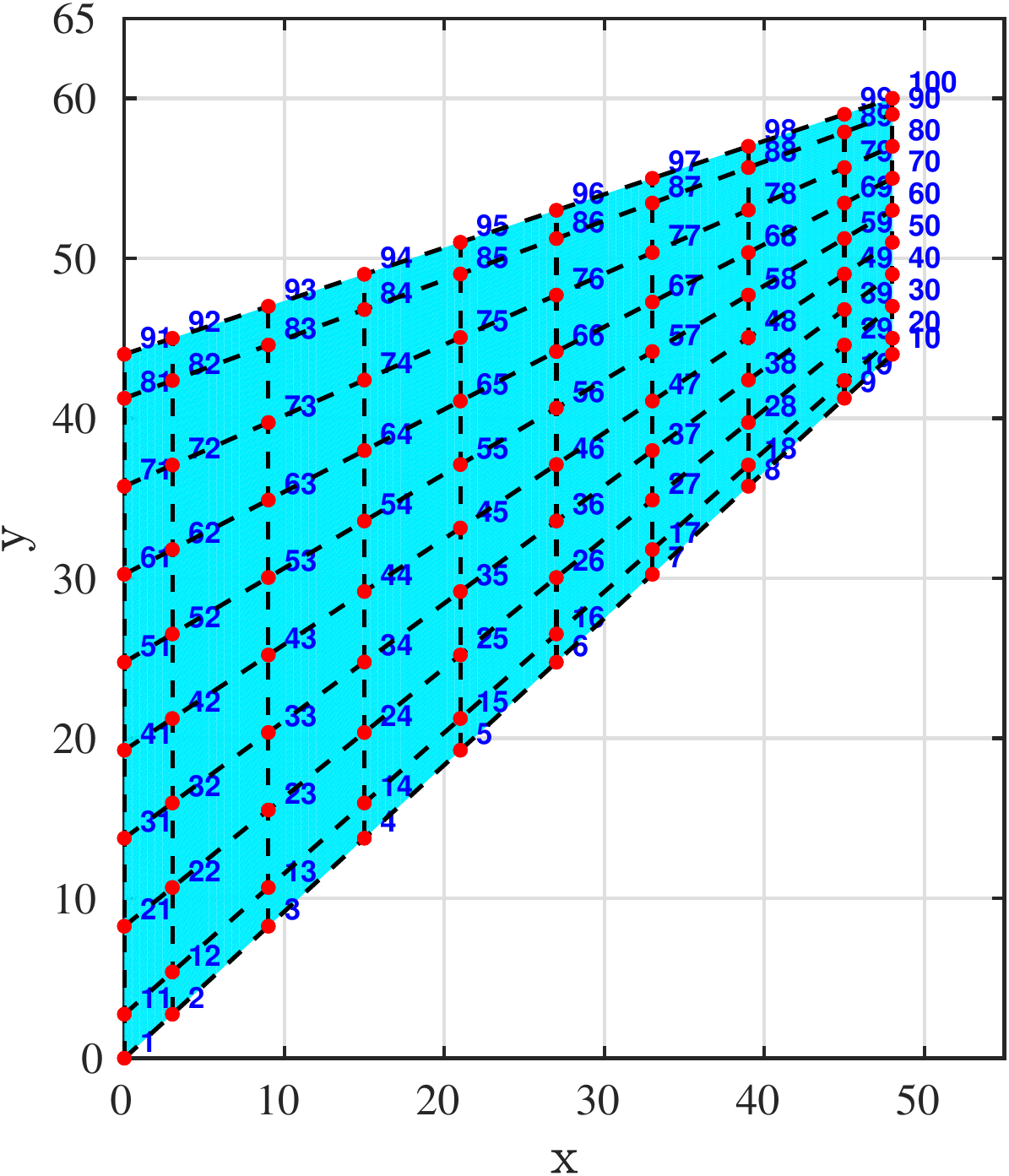}
\caption{CP mesh}
\label{control_point_mesh_cooks_membrane}
\end{subfigure}
\begin{subfigure}{0.31\columnwidth}
\centering
\includegraphics[width=0.9\columnwidth]{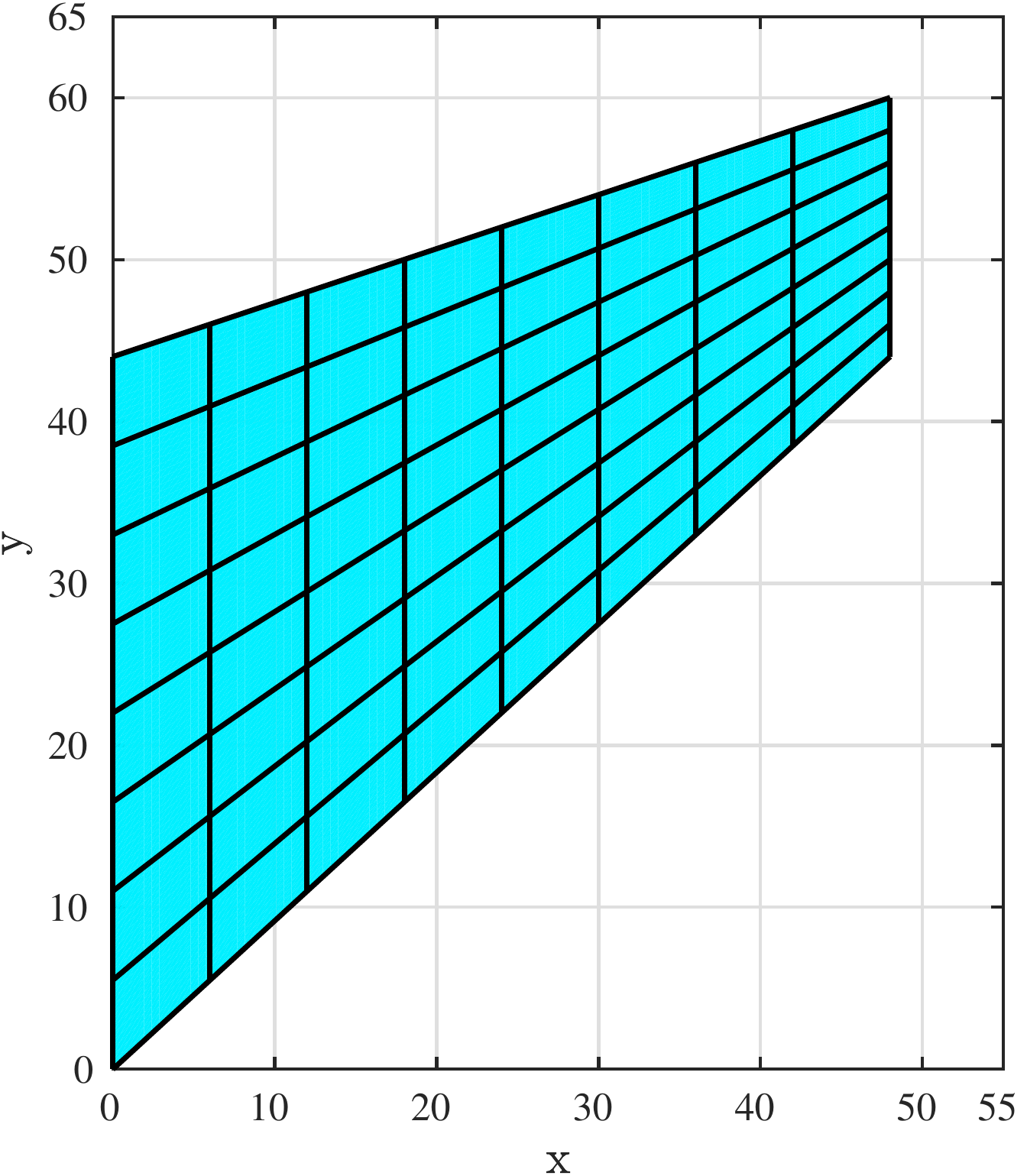}
\caption{Discretized domain}
\label{element_discretization_cooks_membrane}
\end{subfigure} 
\newline
\begin{subfigure}{0.31\columnwidth}
\centering
\includegraphics[width=0.9\columnwidth]{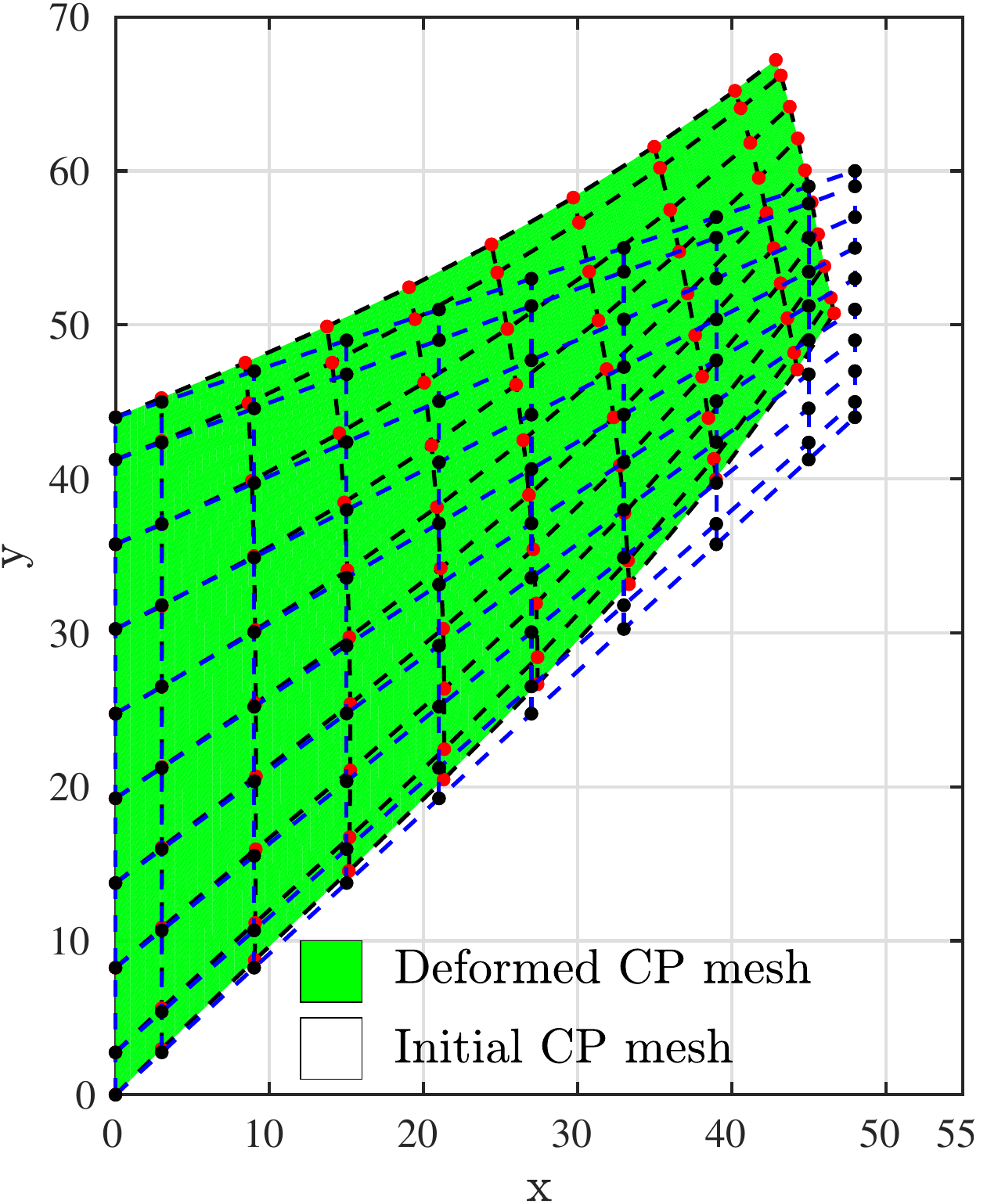}
\caption{Deformed CP mesh}
\label{Deformed_CP_mesh_cooks_membrane}
\end{subfigure}
\begin{subfigure}{0.31\columnwidth}
\centering
\includegraphics[width=0.9\columnwidth]{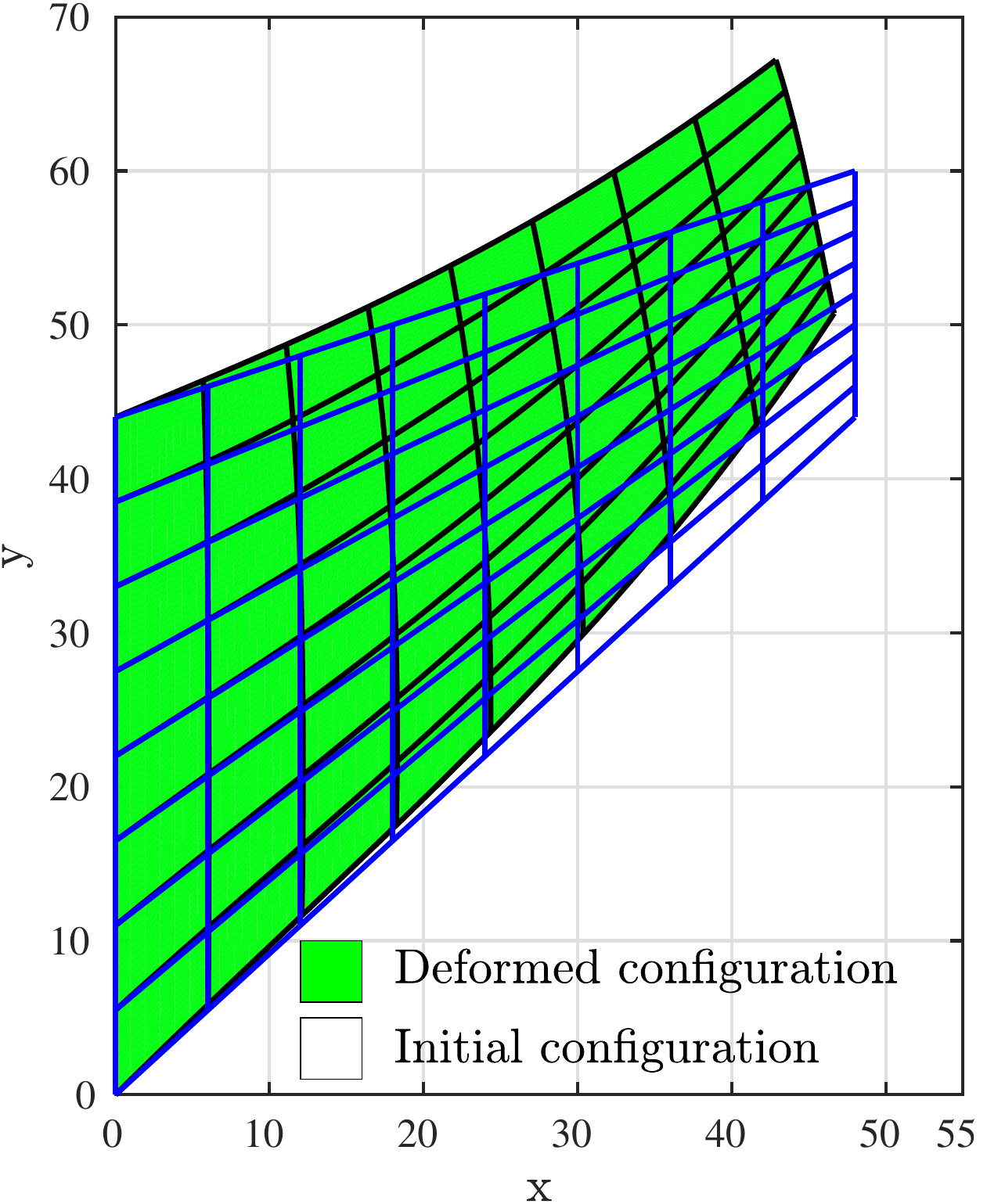}
\caption{Deformed discretized domain}
\label{deformed_element_discretization_Cooks_membrane}
\end{subfigure} 
\begin{subfigure}{0.33\columnwidth}
\centering
\includegraphics[trim={12cm 0cm 8cm 1cm},clip,width=1\columnwidth]{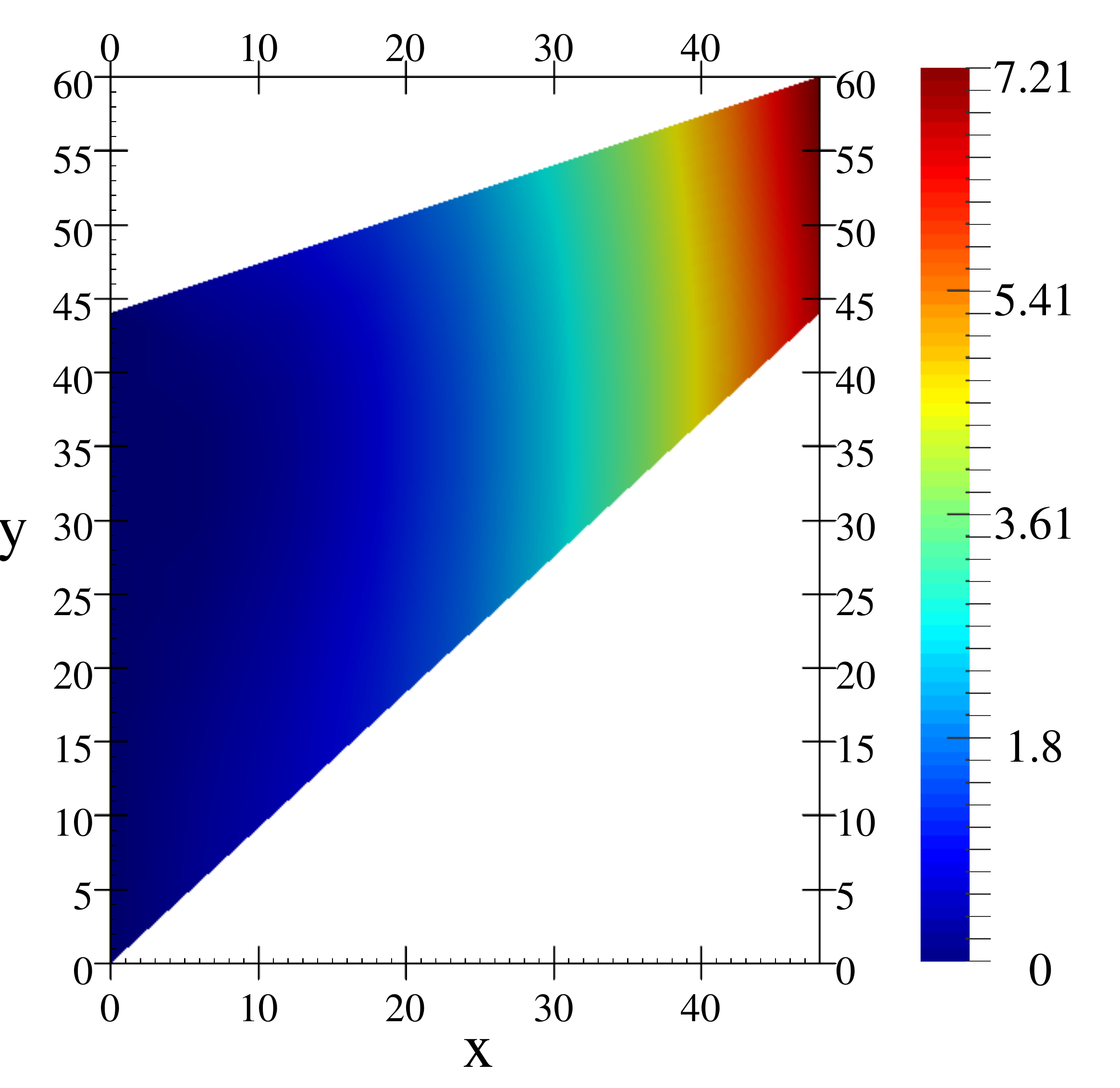}
\caption{$u_y^{H-IGA}$}
\label{Contour_plot_Cooks_membrane}
\end{subfigure}
\caption{(a) The problem setup and boundary conditions for a Cook's membrane problem, (b-c) Geometric description consisting 8$\times$8 quadratic NURBS elements, (d-e) respective deformed configuration, and (f) the contour plot for vertical displacement using hybrid IGA formulation }
\end{figure}

Next, the Cook's membrane problem is simulated \cite{Elguedj2008}. The problem setup and the boundary conditions are illustrated in Figure~\ref{cooks_membrane_problem_definition} where $f_y$ is the load per unit length. Setting $\nu=0.4999$, the problem becomes a typical case of volumetric locking while investigating the nearly incompressible behavior of the domain under combined bending and shear deformation.

The geometric data to construct the coarsest possible mesh, representing the exact geometry, is provided in the appendix (Table~\ref{CP_and_weights_for_cooks_membrane_problem}). Once the initial mesh is generated, the sequence of refined meshes are modeled using the refinement techniques. One such mesh of $8 \times 8$ NURBS elements, having the quadratic basis functions along the $\xi$ and $\eta$ direction, is illustrated in Figure~\ref{control_point_mesh_cooks_membrane}-\ref{element_discretization_cooks_membrane}. The elaborated results are presented in the Figure~\ref{Deformed_CP_mesh_cooks_membrane}-\ref{Contour_plot_Cooks_membrane} which focuses on the deformed configuration of control point mesh and the discretized domain along with the contour plot for vertical displacement for the stated mesh.

\begin{figure}
\begin{subfigure}{0.48\columnwidth}
\centering
\includegraphics[width=1\columnwidth]{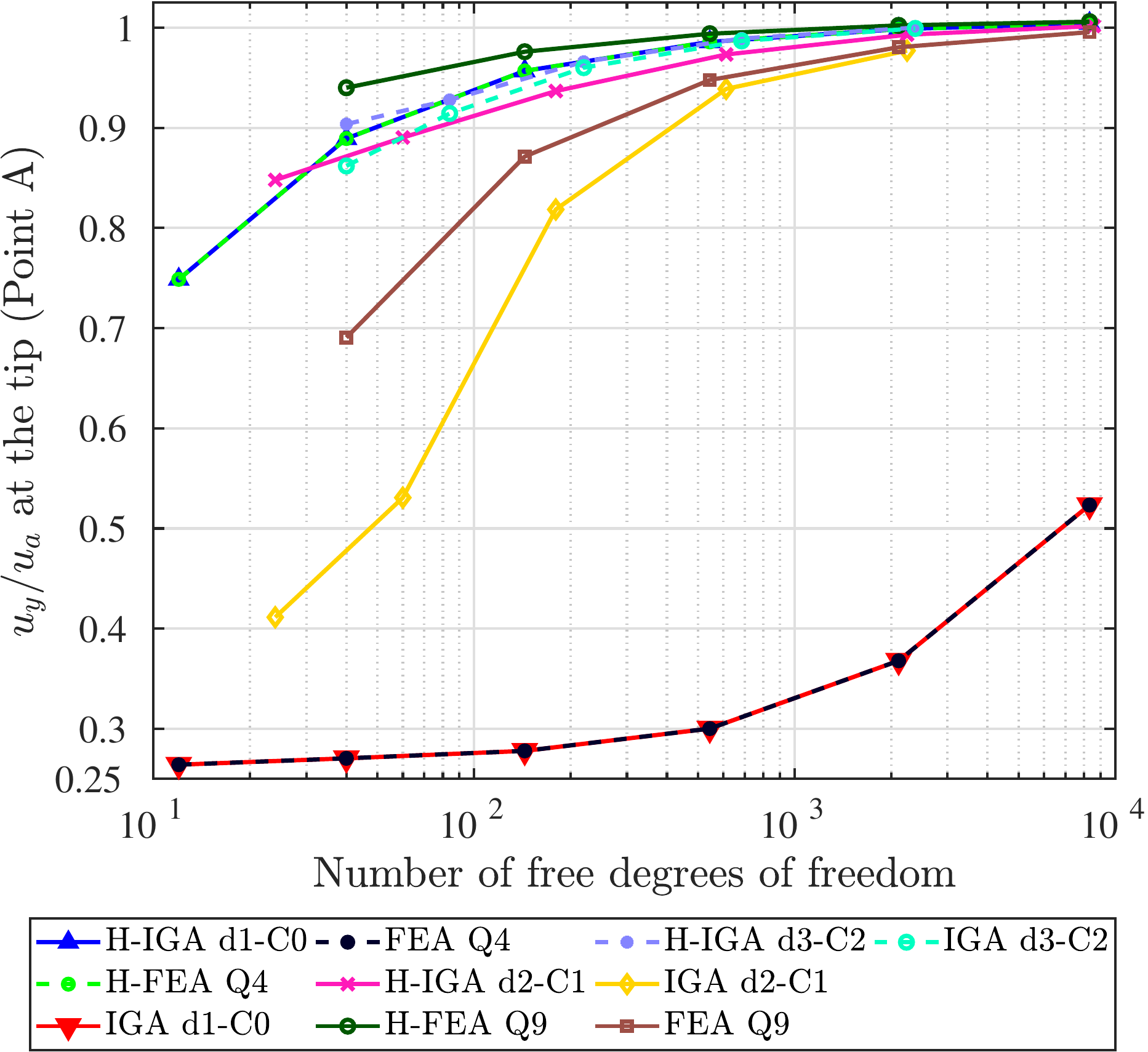}
\caption{$u_y/u_a$ }
\label{Normalized_displacement_Cooks_membrane}
\end{subfigure}
\begin{subfigure}{0.48\columnwidth}
\centering
\includegraphics[width=1\columnwidth]{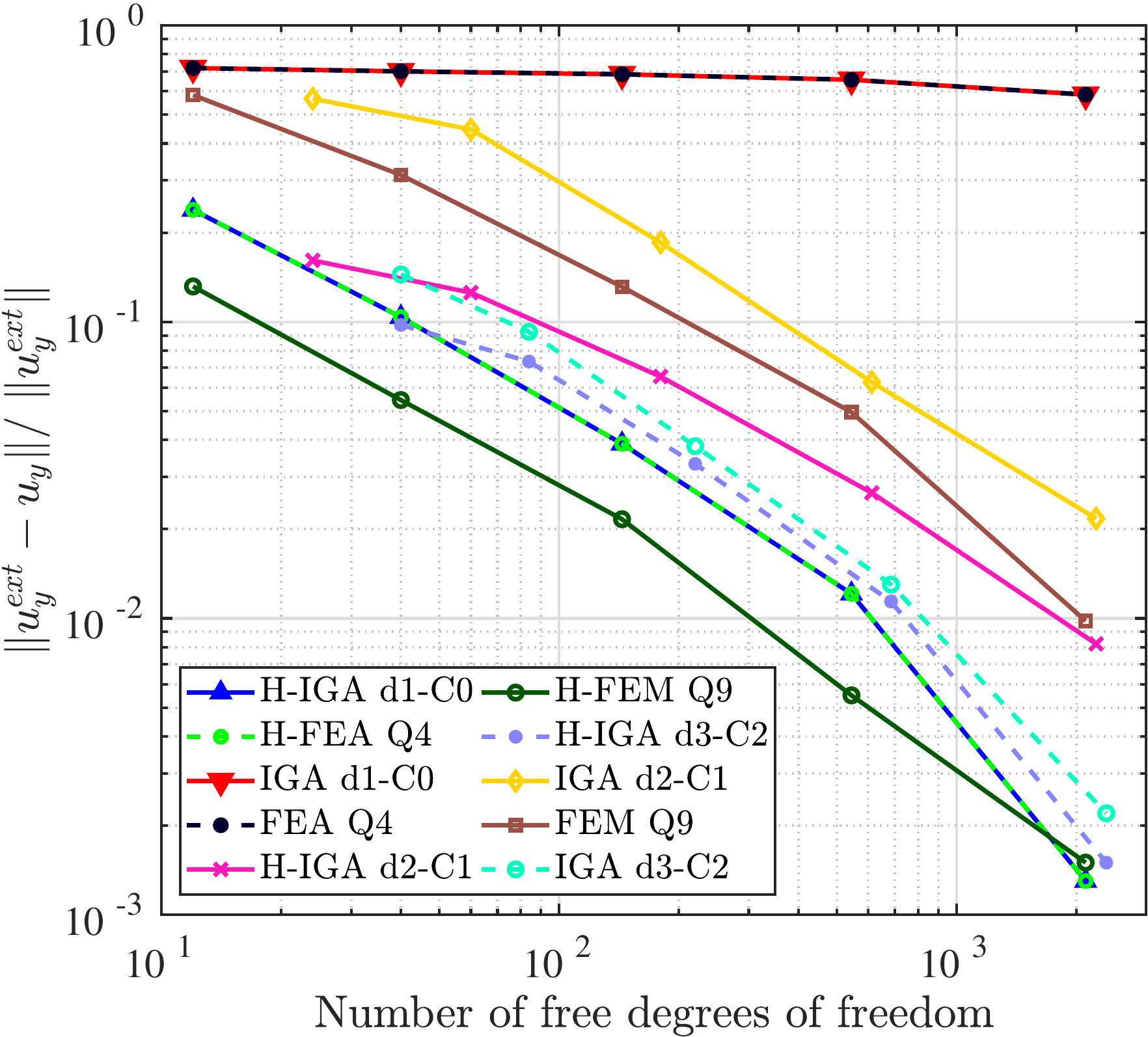}
\caption{$L_2$ error norm for vertical displacement}
\label{L_2_error_cooks_membrane_problem}
\end{subfigure}
\caption{(a) Normalized vertical displacement at point `A' and $L_2$ error norm for vertical displacement versus the active degrees of freedom for the Cook's membrane problem}\label{Cooks_membrane_solutions}
\end{figure}

The normalized vertical displacement at the point `A' against the reference solution of 7.7 is evaluated as shown in Figure~\ref{Normalized_displacement_Cooks_membrane}. Furthermore, the convergence of the relative $L_2$ error norm of displacement versus the number of active degrees of freedom is shown in Figure~\ref{L_2_error_cooks_membrane_problem}. As analytical expression for the displacements are not well established for the stated problem, the reference to evaluate the $L_2$ norm is the well converged solution of high-refined mesh of cubic degree basis functions.

It can be observed that the conventional FEA for four node quadrilateral elements and its equivalent IGA formulation locks severely. Even with the significant high refinement, the results only marginally improve. On the other hand, hybrid formulation can successfully alleviate locking to produce the superior results. Moving to higher degree basis functions, where the IGA basis are different from the Lagrangian basis functions, it can be seen that the quadratic NURBS elements locks significantly for lower mesh refinements. However, the proposed hybrid IGA works very well even with very low number of active degrees of freedom and providing the better coarse mesh accuracy. With further elevation in the degree of basis function, the gap between the conventional and hybrid IGA results gets insignificant yet hybrid IGA is marginally better that the conventional IGA formulation. 

\subsection{Infinite plate with a hole problem}

\begin{figure}[pos = ht]
\begin{subfigure}{0.38\columnwidth}
\centering
\includegraphics[width=1\columnwidth]{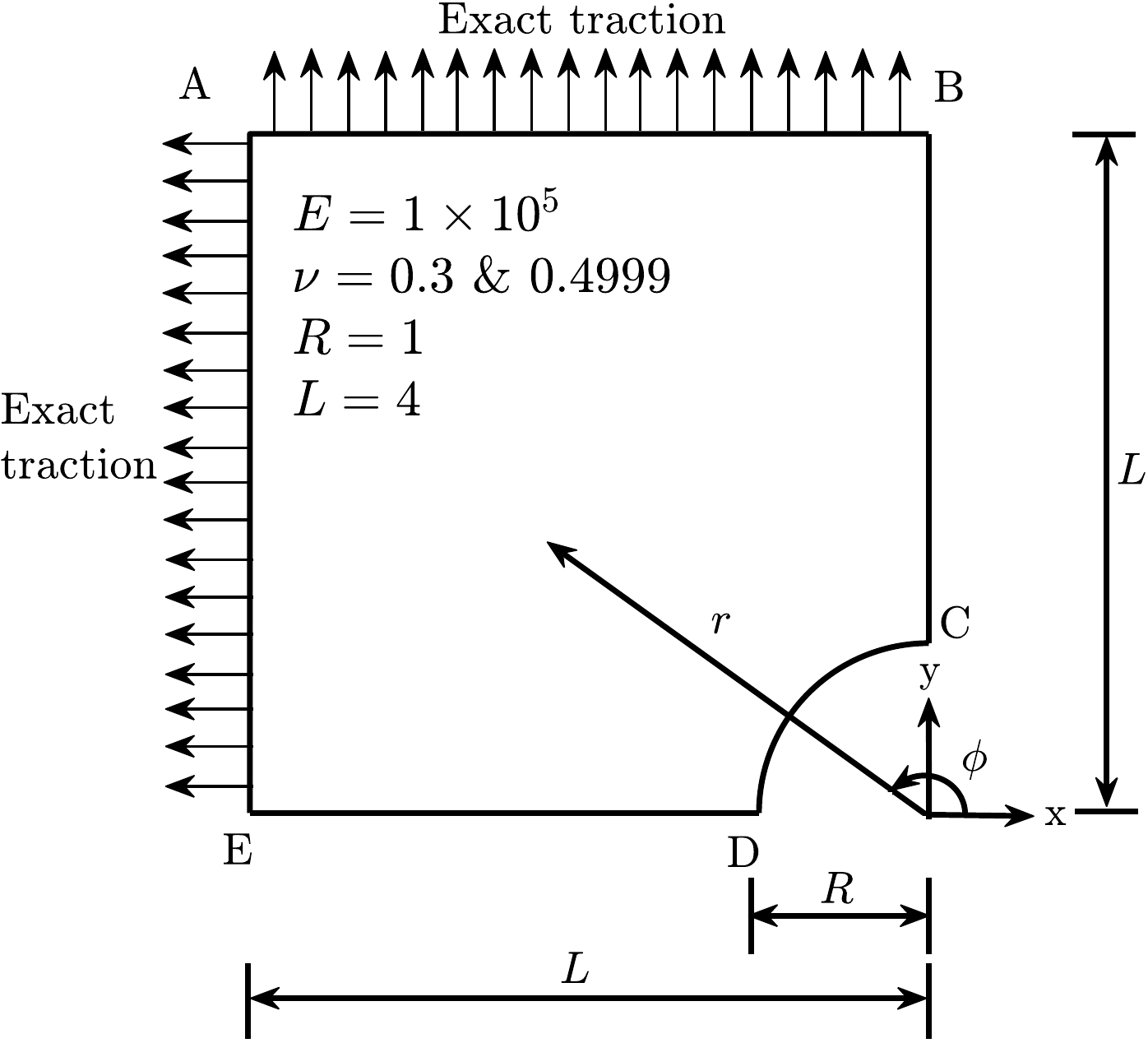}
\caption{The problem setup and boundary conditions}
\label{Plate_with_a_hole_problem_definition}
\end{subfigure}
\begin{subfigure}{0.3\columnwidth}
\centering
\includegraphics[width=1\columnwidth]{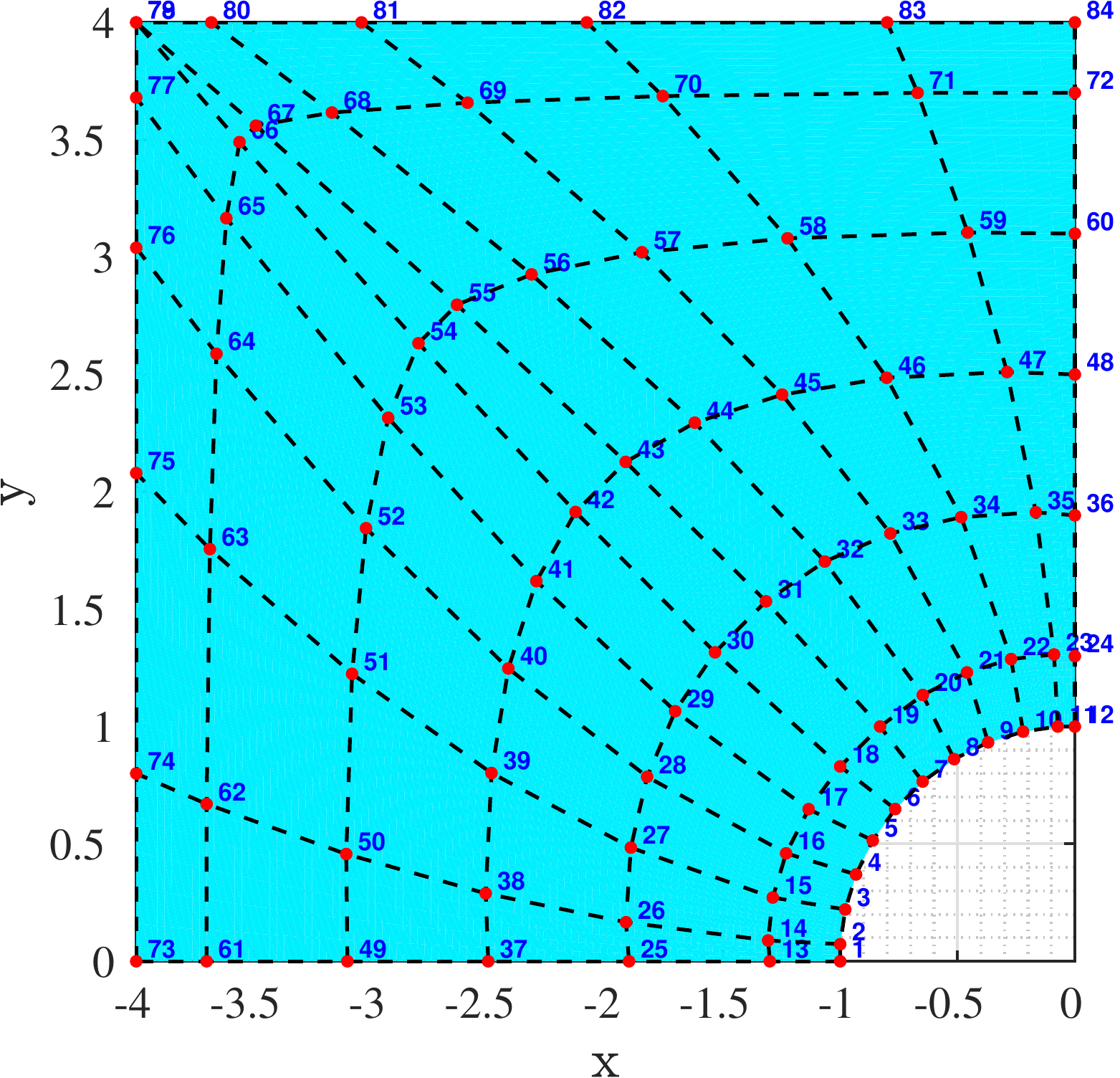}
\captionsetup{width=8.48cm}
\caption{Control point mesh}
\label{Control_point_mesh_Plate_with_a_hole_problem}
\end{subfigure}
\begin{subfigure}{0.3\columnwidth}
\centering
\includegraphics[width=1\columnwidth]{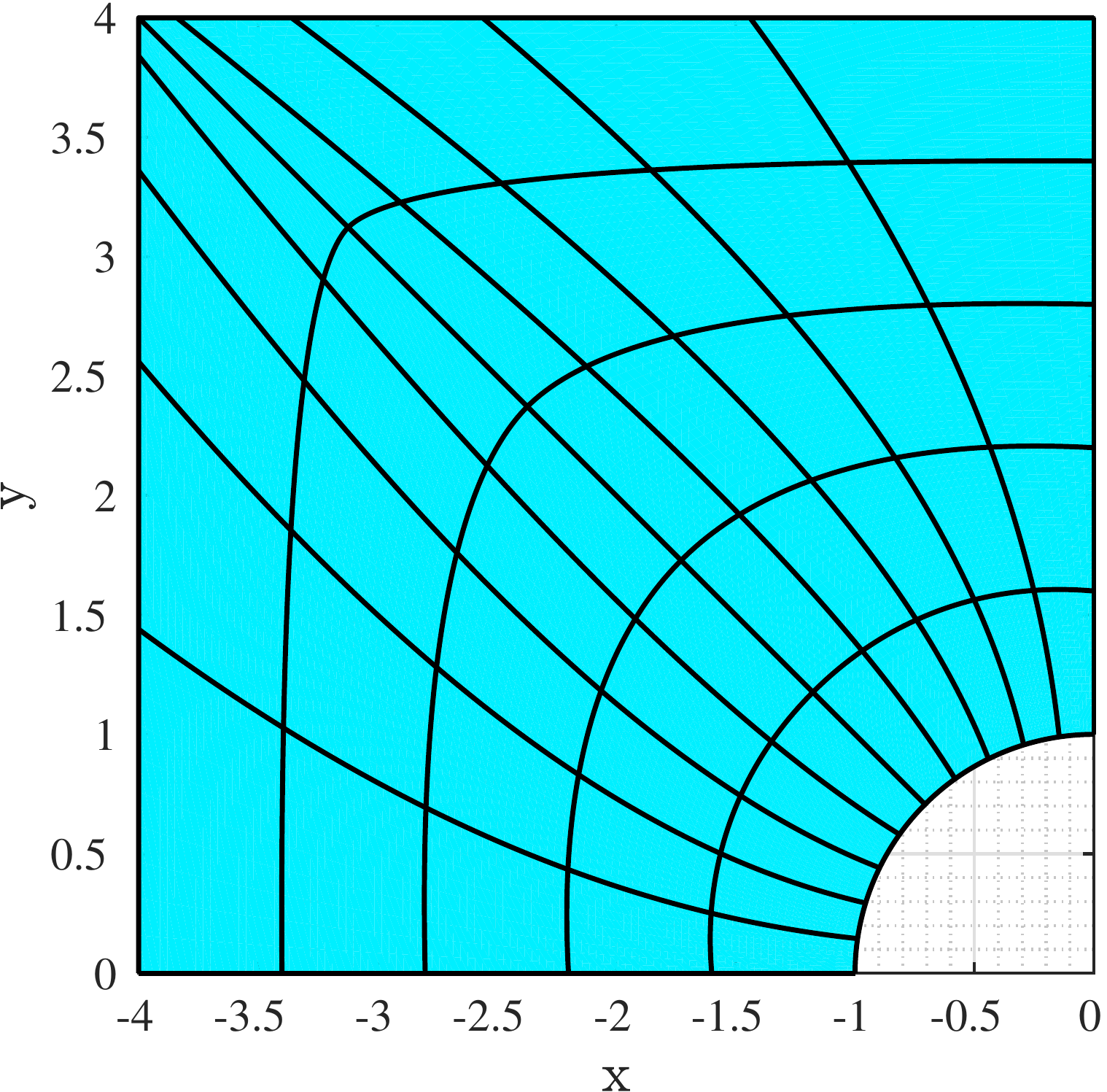}
\captionsetup{width=8.48cm}
\caption{Element discretization}
\label{Element_discretization_Plate_with_a_hole_problem}
\end{subfigure} 
\caption{(a) The problem setup and boundary conditions for an infinite plate with a circular hole problem under constant in-plane tension ($T = 1$) at infinity, and (b-c) Geometric description (a quarter portion) consisting 10$\times$5 NURBS elements with quadratic basis along $\xi$ and $\eta$ direction}
\end{figure}

The problem, see Figure~\ref{Plate_with_a_hole_problem_definition}, consists of a two-dimensional infinite plate with a circular hole under constant uni-axial in-plane tension ($T = 1$) at infinity \cite{Elguedj2008}. Owing to symmetry of the problem, only a quarter portion of the plate is considered, as shown in Figure~\ref{Plate_with_a_hole_problem_definition}.

The problem setup and the boundary conditions are illustrated in Figure~\ref{Plate_with_a_hole_problem_definition}, which includes the symmetric boundary condition on edge ED and BC, and the Neumann boundary condition on edge AB and AE. The exact traction, applied on the boundary AB and AE, is evaluated using the analytical expression of stresses (Eqn.~\ref{sigma_xx},~\ref{sigma_yy},and~\ref{sigma_xy}) for the stated problem setup.

A plane strain condition is assumed, and the analytical solution \cite{Elguedj2008} for the stress field is given as follows,
% ======================= stresses =============================
\begin{align}
\sigma_{xx} &= 1 - \frac{R^2}{r^2}\bigg( \frac{3}{2} \cos 2\phi + \cos 4\phi \bigg) + \frac{3R^4}{2r^4} \cos 4\phi \label{sigma_xx}\\
\sigma_{yy} &= - \frac{R^2}{r^2}\bigg( \frac{1}{2} \cos 2\phi - \cos 4\phi \bigg) - \frac{3R^4}{2r^4} \cos 4\phi \label{sigma_yy}\\
\tau_{xy} &=  -\frac{R^2}{r^2}\bigg( \frac{1}{2} \sin 2\phi + \sin 4\phi \bigg) + \frac{3R^4}{2r^4} \sin 4\phi \label{sigma_xy}
\end{align}
% -------------------------------------------------------------
where, $r = \sqrt{x^2 + y^2}$ and $\phi = \tan^{-1}(y/x)$. The analytical expression for the displacement field is given as, 
% ================== Displacement =============================
\begin{align*}
u_x(r,\phi) &= \frac{R}{8\mu}\bigg[ \frac{r}{R}(k+1)\cos \phi + \frac{2R}{r}((1+k)\cos \phi + \cos 3\phi) - \frac{2R^3}{r^3} \cos 3\phi \bigg] \\
u_y(r,\phi) &= \frac{R}{8\mu}\bigg[ \frac{r}{R}(k-3)\sin \phi + \frac{2R}{r}((1-k)\sin \phi + \sin 3\phi) - \frac{2R^3}{r^3} \sin 3\phi \bigg] 
\end{align*}
% ------------------------------------------------------------
where, $\mu = \dfrac{E}{2(1+\nu)}$ and $k = 3-4\nu$ (for plane strain condition).

\begin{figure}[pos=H]
\begin{subfigure}{0.49\columnwidth}
\centering
\includegraphics[trim={0cm 0cm 0cm 0cm},clip,width=0.8\columnwidth]{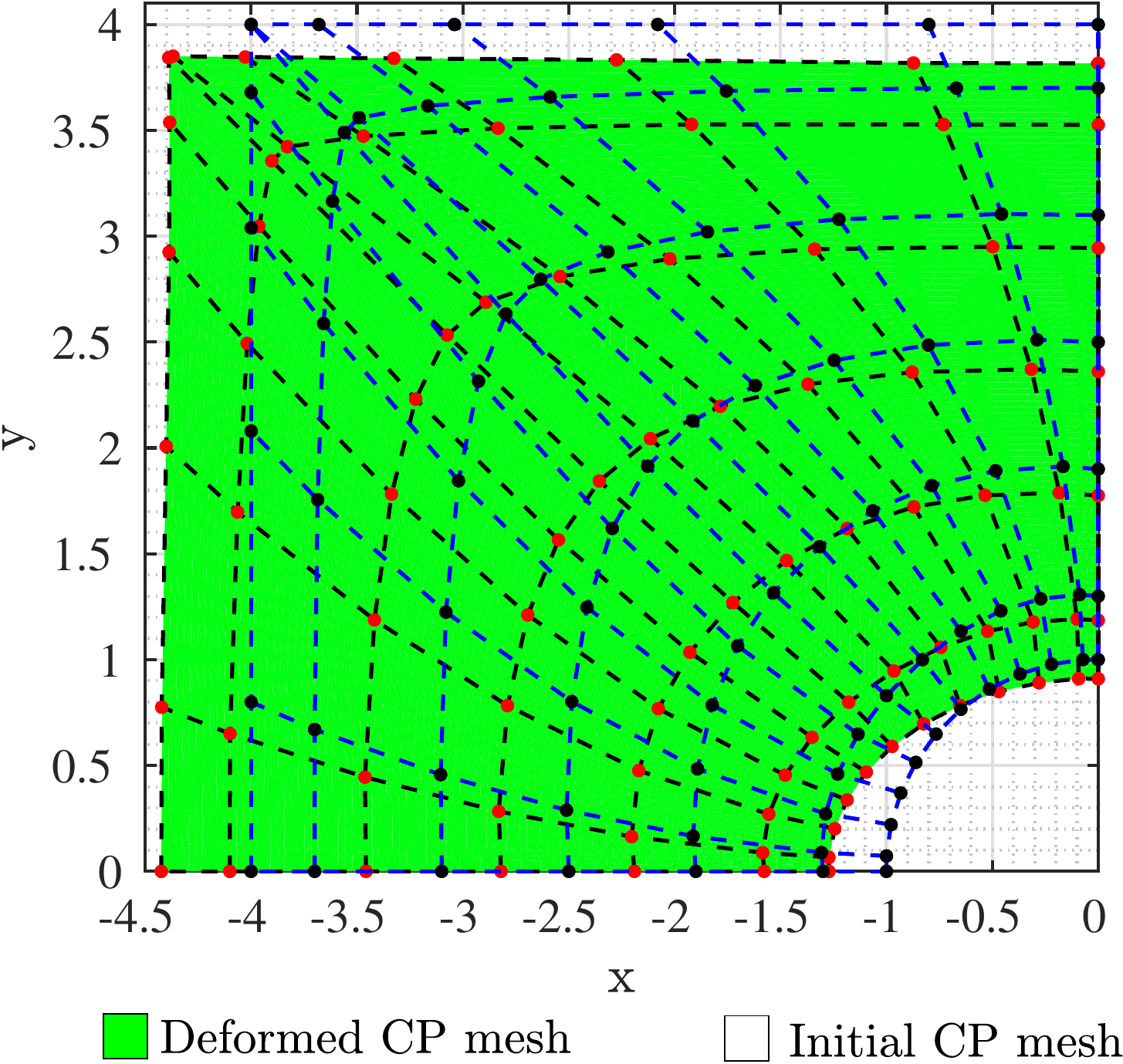}
\caption{Deformed CP mesh}
\label{Deformed_CP_mesh_Plate_with_a_hole_problem_03}
\end{subfigure}
\begin{subfigure}{0.49\columnwidth}
\centering
\includegraphics[trim={0cm 0cm 0cm 0cm},clip,width=0.8\columnwidth]{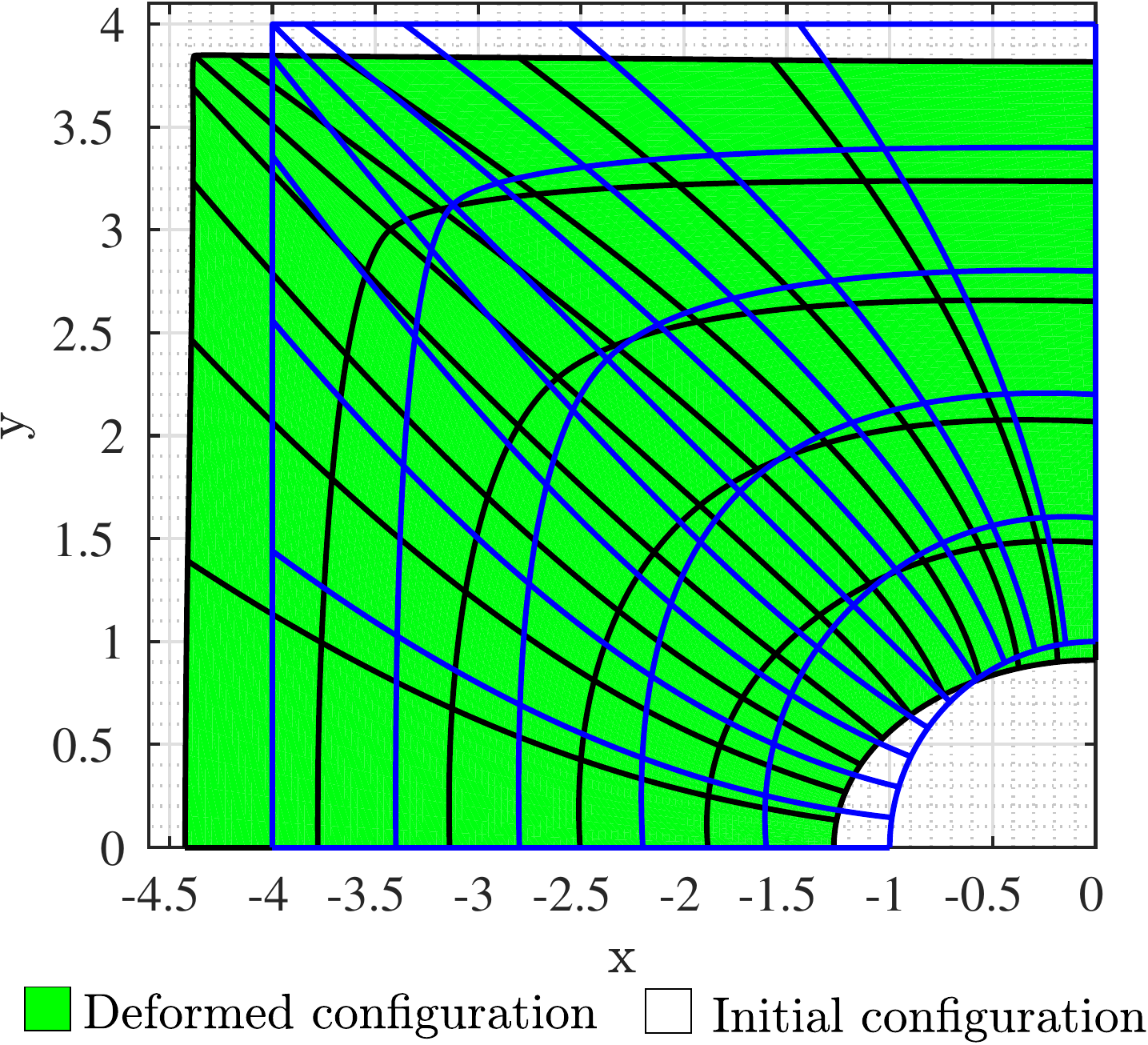}
\caption{Deformed discretized domain}
\label{Deformed_Element_discretization_Plate_with_a_hole_problem_03}
\end{subfigure} 
\newline
\begin{subfigure}{0.49\columnwidth}
\centering
\includegraphics[trim={0cm 0cm 0cm 0cm},clip,width=0.8\columnwidth]{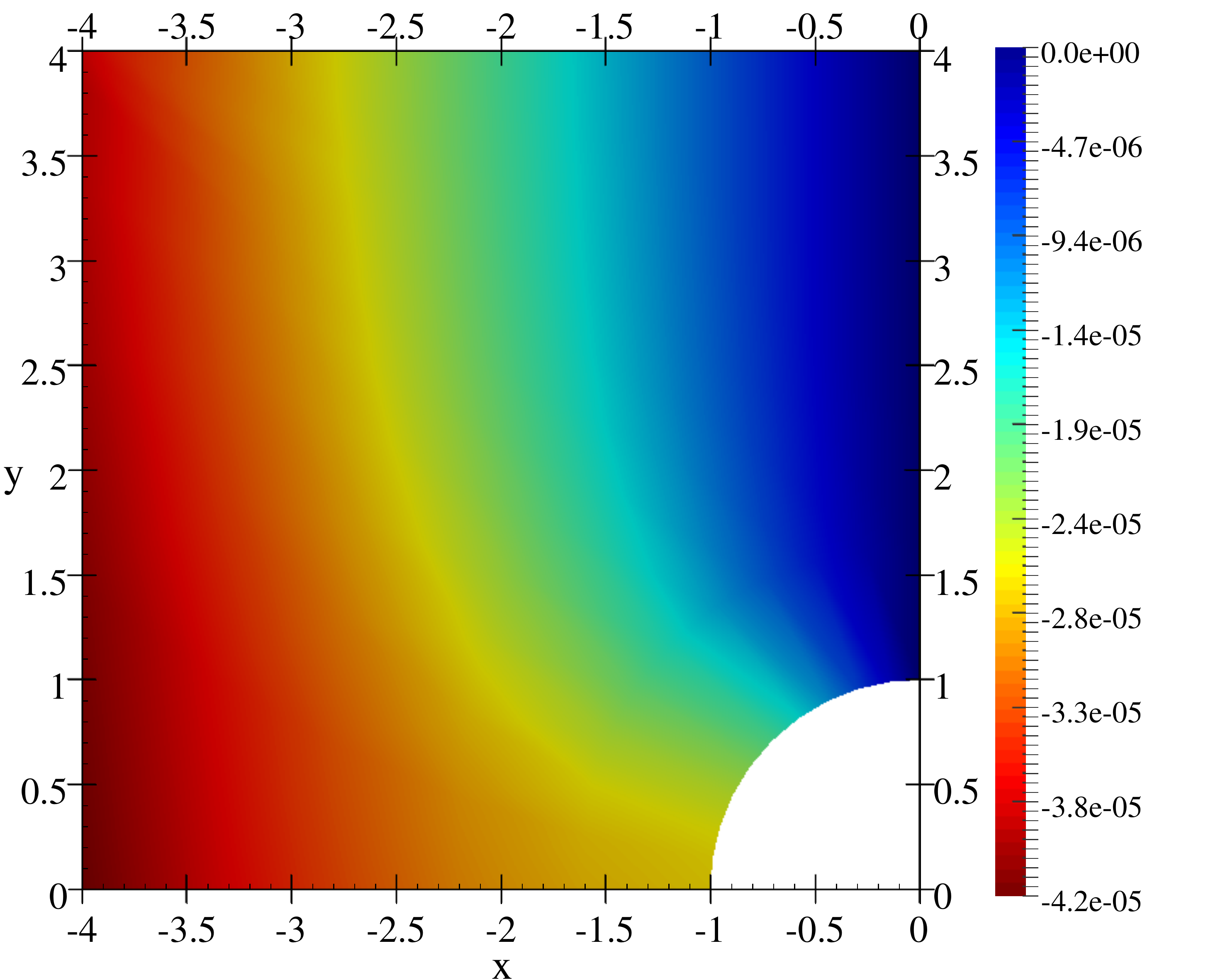}
\caption{$u_x^{H-IGA}$}
\label{u_x_numerical_plate_with_hole_nu_03}
\end{subfigure}
\begin{subfigure}{0.49\columnwidth}
\centering
\includegraphics[trim={0cm 0cm 0cm 0cm},clip,width=0.8\columnwidth]{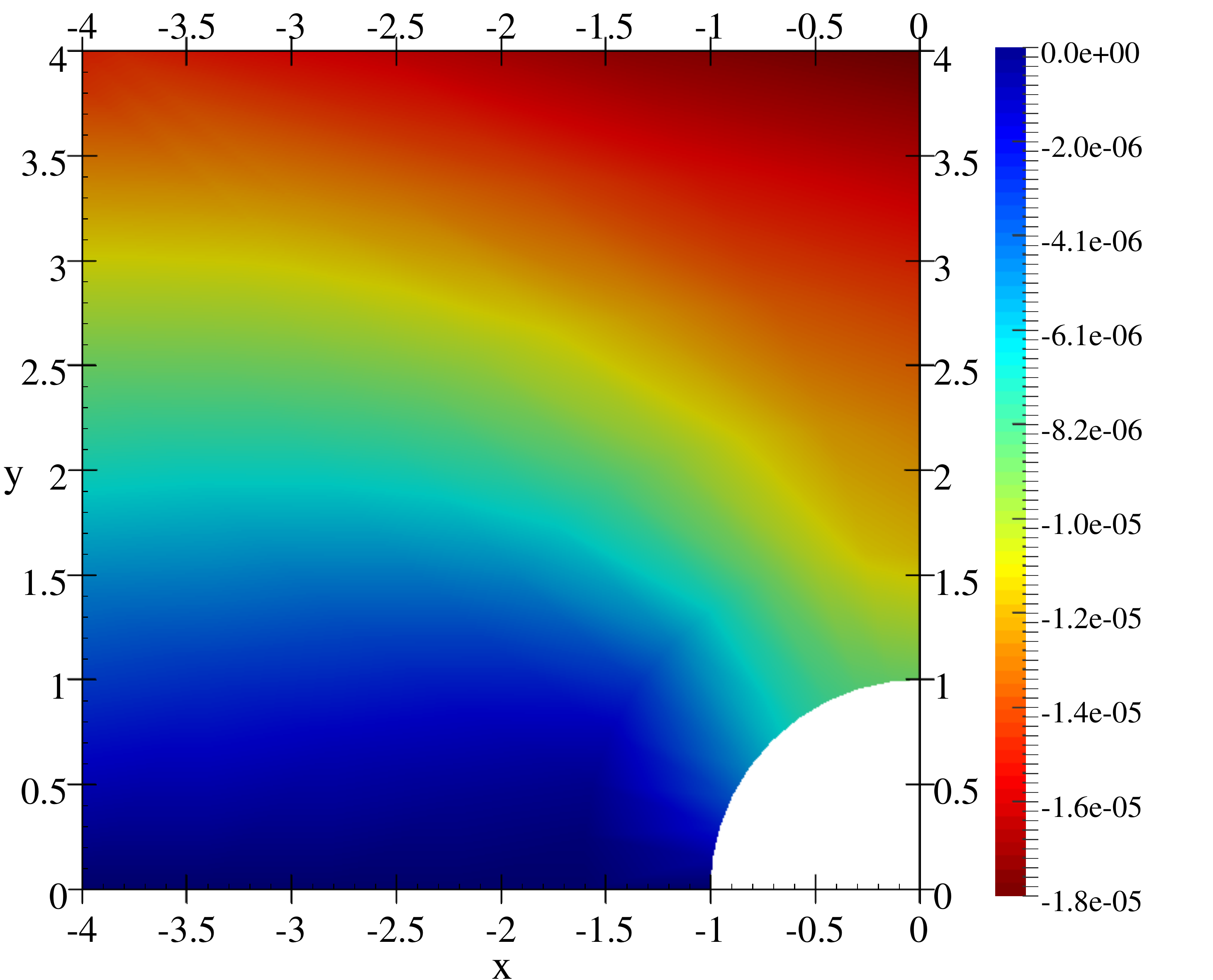}
\caption{$u_y^{H-IGA}$}
\label{u_y_numerical_plate_with_hole_nu_03}
\end{subfigure} 
\newline
\begin{subfigure}{0.49\columnwidth}
\centering
\includegraphics[trim={0cm 0cm 0cm 0cm},clip,width=0.8\columnwidth]{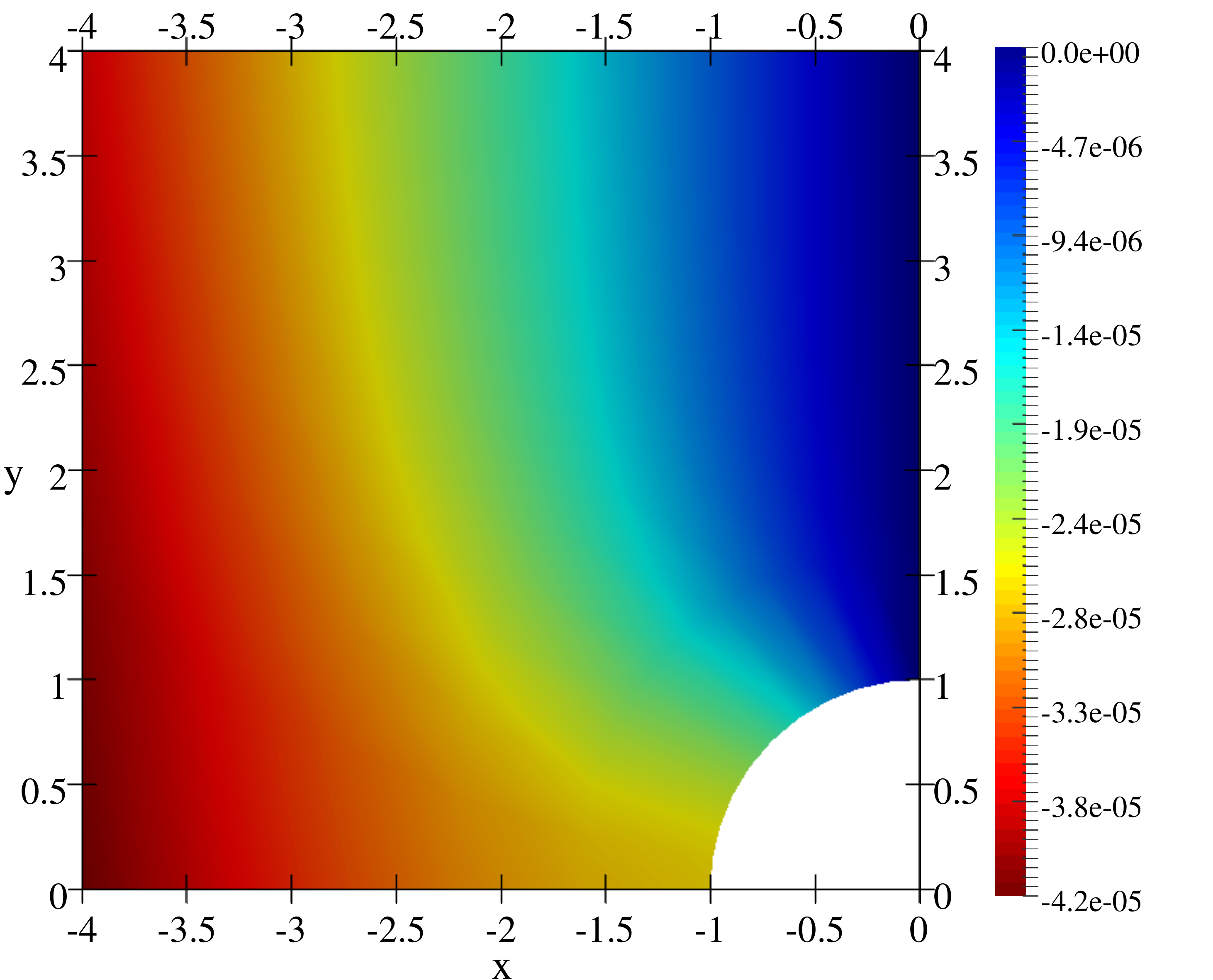}
\caption{$u_x^{exact}$}
\label{u_x_numerical_plate_with_hole_nu_03_exact}
\end{subfigure}
\begin{subfigure}{0.49\columnwidth}
\centering
\includegraphics[trim={0cm 0cm 0cm 0cm},clip,width=0.8\columnwidth]{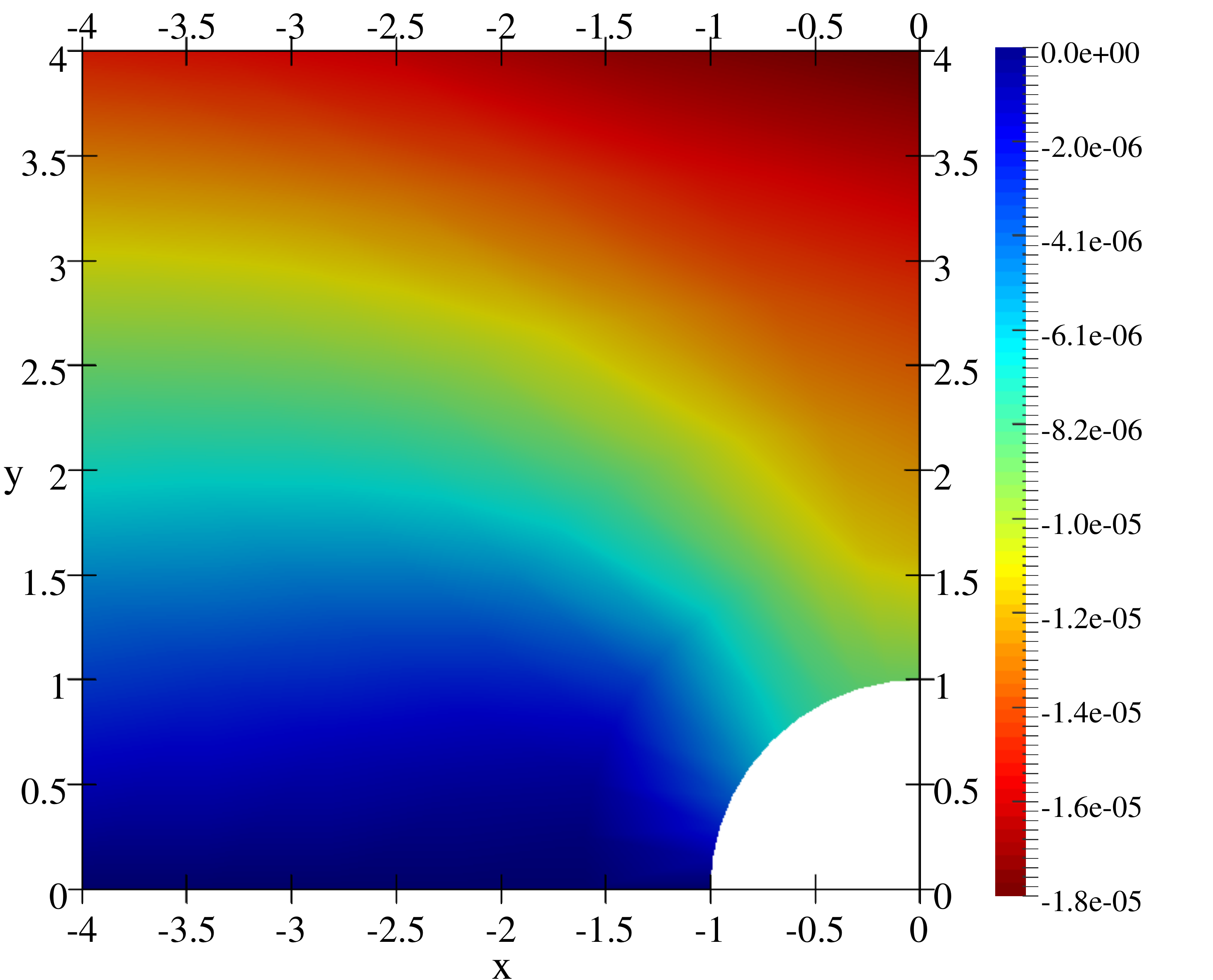}
\caption{$u_y^{exact}$}
\label{u_y_numerical_plate_with_hole_nu_03_exact}
\end{subfigure} 
\caption{(a-b) Deformed geometric description of an infinite plate with a hole problem (a quarter portion) for 10$\times$5 NURBS elements with quadratic basis along $\xi$ and $\eta$ direction (magnified by the factor of $1\times10^4$ for better visualization), and contour plots for horizontal and vertical displacement for the mesh illustrated in Figure~\ref{Control_point_mesh_Plate_with_a_hole_problem}-\ref{Element_discretization_Plate_with_a_hole_problem} using hybrid IGA formulation (c-d) alongside the analytical solution (e-f) for $\nu=0.3$}\label{comprehensive_results_plate_with_hole_03}
\end{figure}

\begin{figure}[pos=H]
\begin{subfigure}{0.49\columnwidth}
\centering
\includegraphics[trim={0cm 0cm 0cm 0cm},clip,width=0.8\columnwidth]{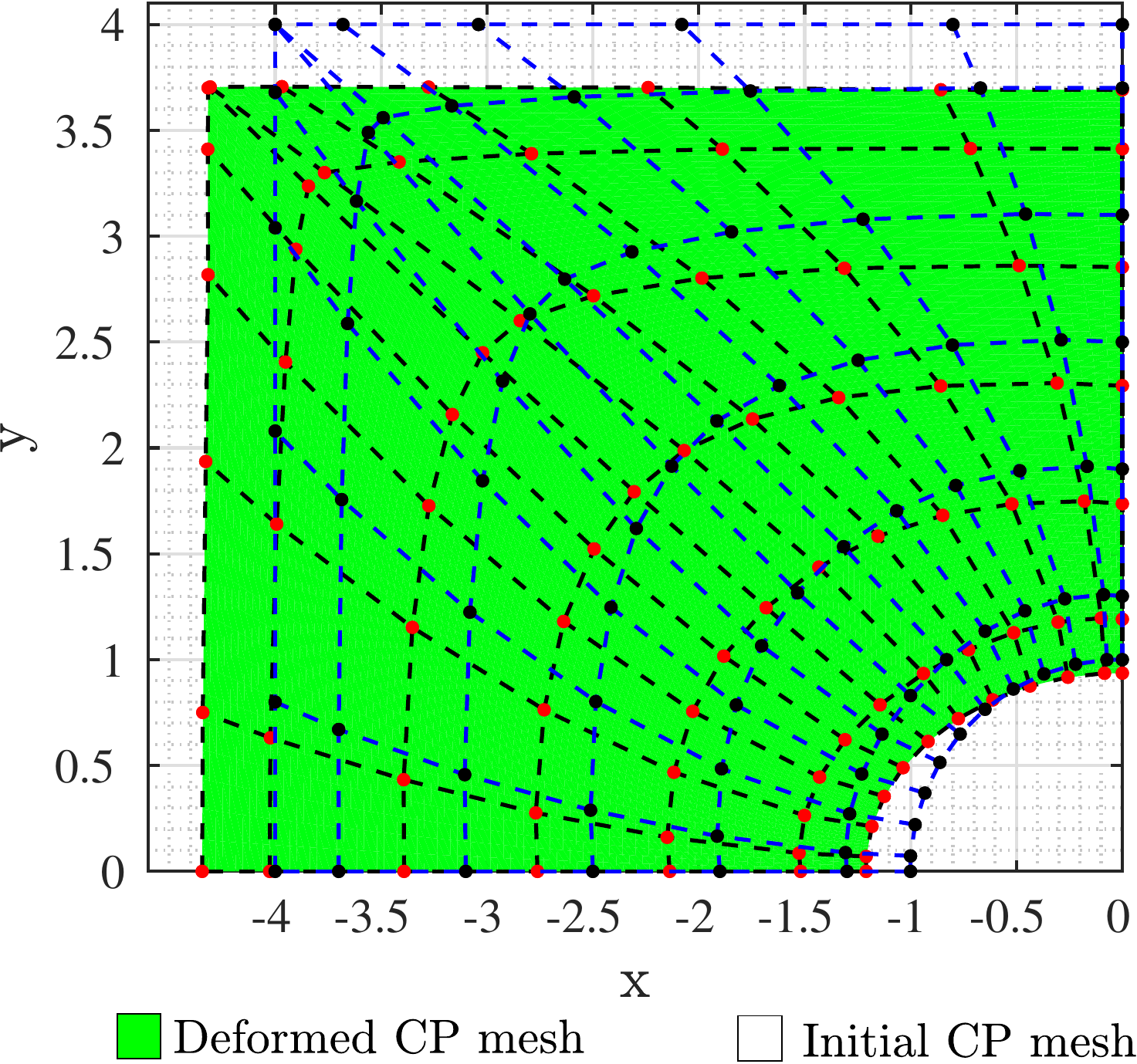}
\captionsetup{width=8.48cm}
\caption{Deformed control point mesh}
\label{Deformed_CP_mesh_Plate_with_a_hole_problem}
\end{subfigure}
\begin{subfigure}{0.49\columnwidth}
\centering
\includegraphics[trim={0cm 0cm 0cm 0cm},clip,width=0.8\columnwidth]{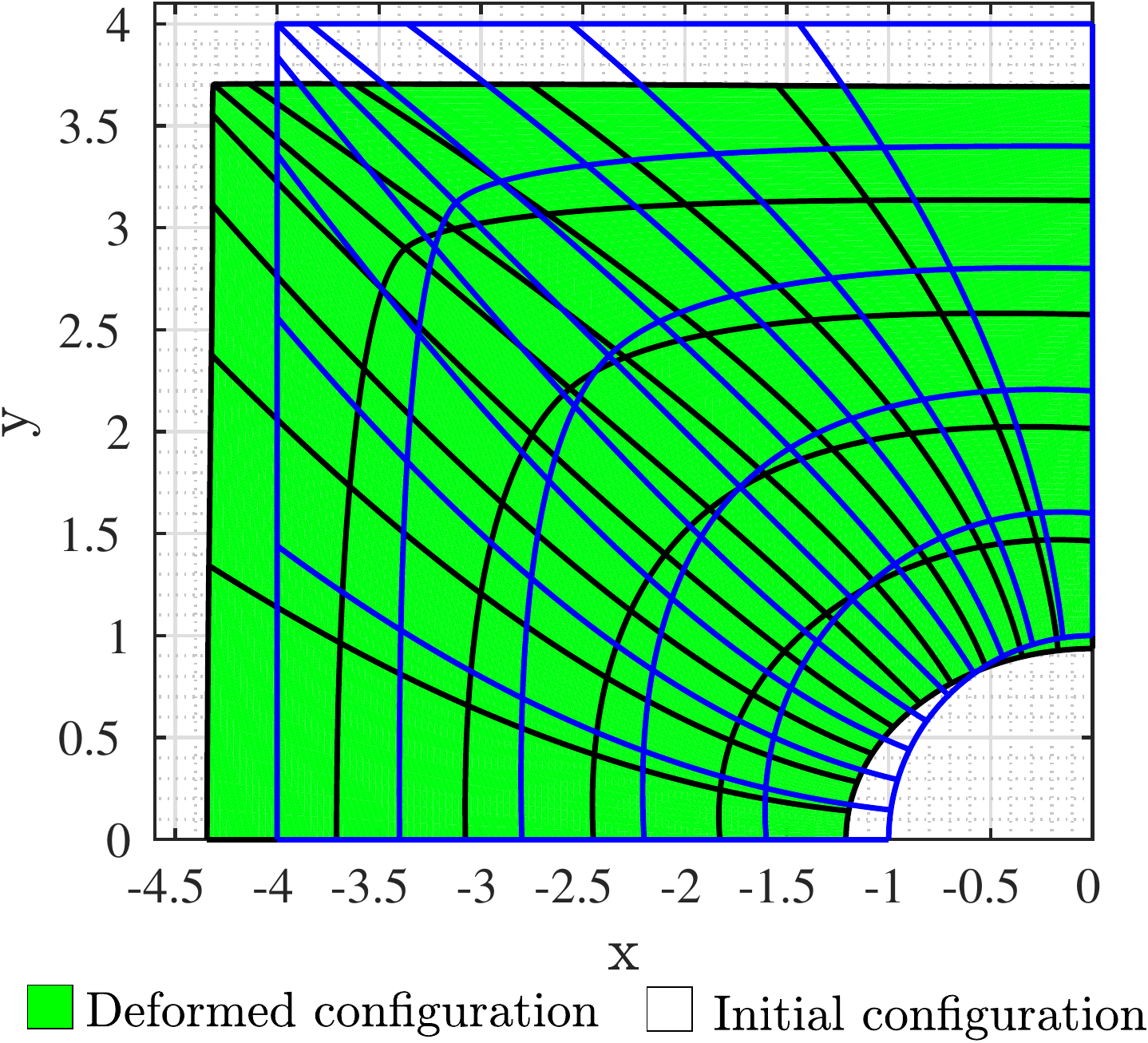}
\captionsetup{width=8.48cm}
\caption{Deformed discretized problem domain}
\label{Deformed_Element_discretization_Plate_with_a_hole_problem}
\end{subfigure} 
\newline
\begin{subfigure}{0.49\columnwidth}
\centering
\includegraphics[trim={0cm 0cm 0cm 0cm},clip,width=0.8\columnwidth]{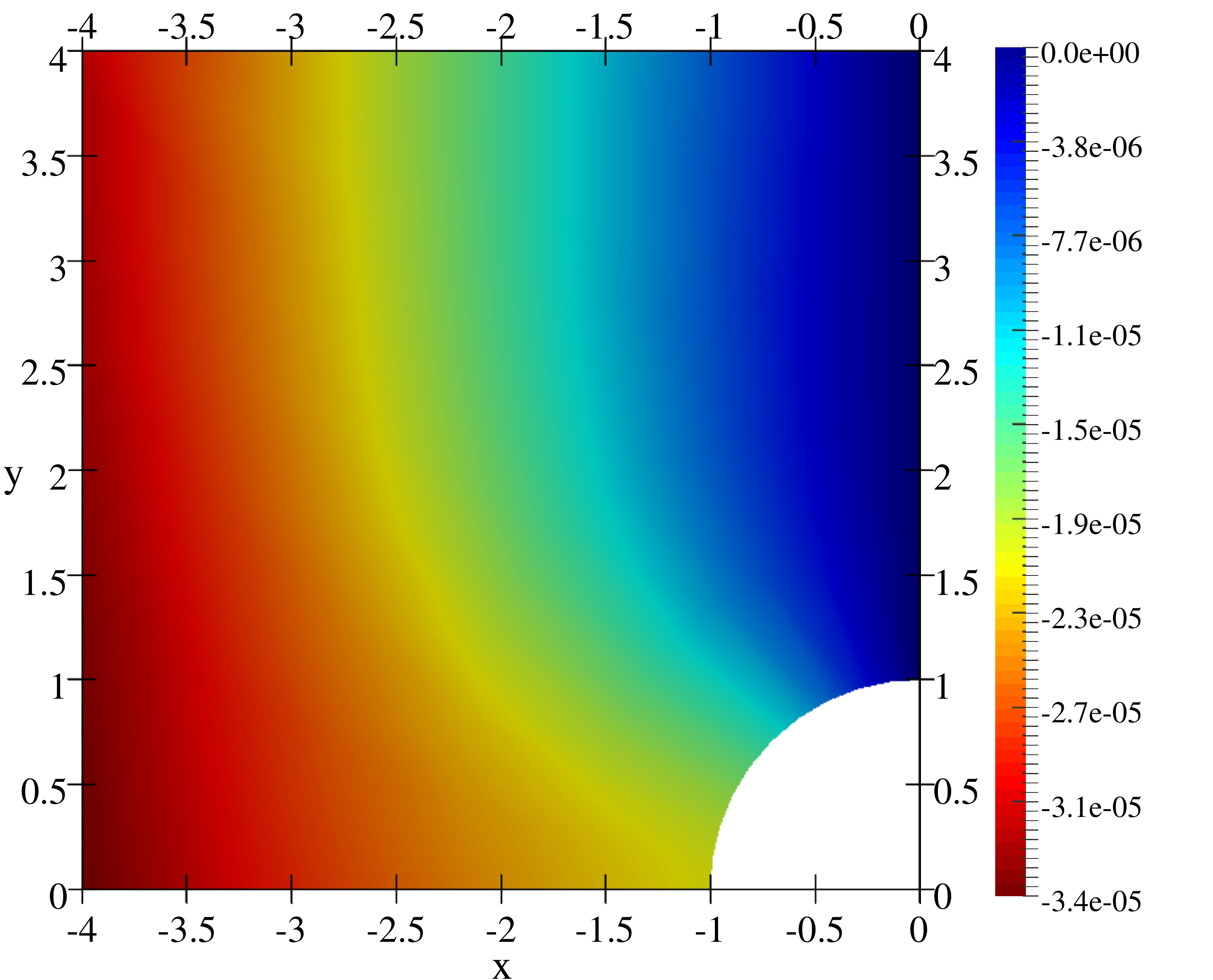}
\captionsetup{width=8.48cm}
\caption{$u_x^{H-IGA}$}
\label{u_x_numerical_plate_with_hole_nu_04999}
\end{subfigure}
\begin{subfigure}{0.49\columnwidth}
\centering
\includegraphics[trim={0cm 0cm 0cm 0cm},clip,width=0.8\columnwidth]{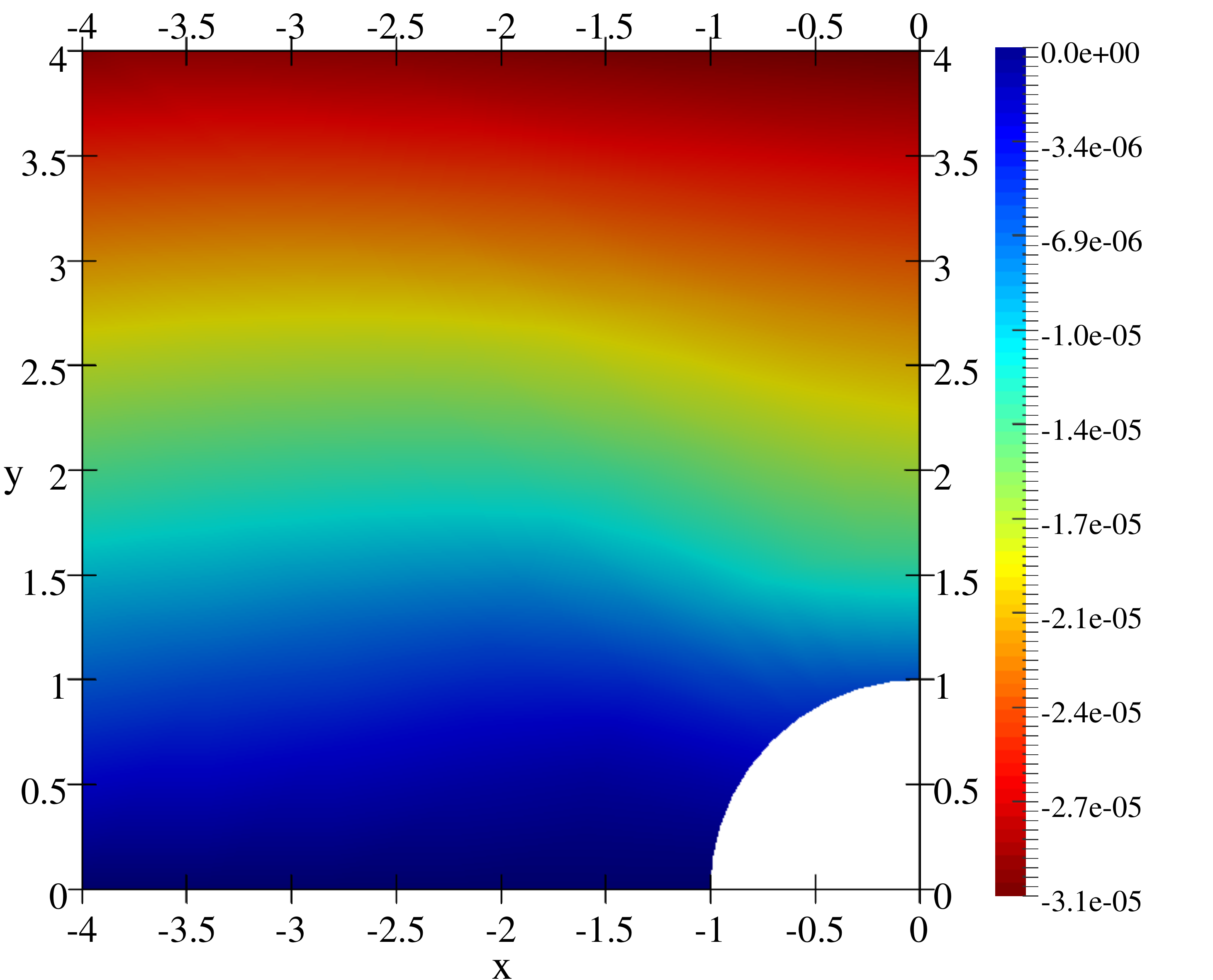}
\captionsetup{width=8.48cm}
\caption{$u_y^{H-IGA}$}
\label{u_y_numerical_plate_with_hole_nu_04999}
\end{subfigure} 
\newline
\begin{subfigure}{0.49\columnwidth}
\centering
\includegraphics[trim={0cm 0cm 0cm 0cm},clip,width=0.8\columnwidth]{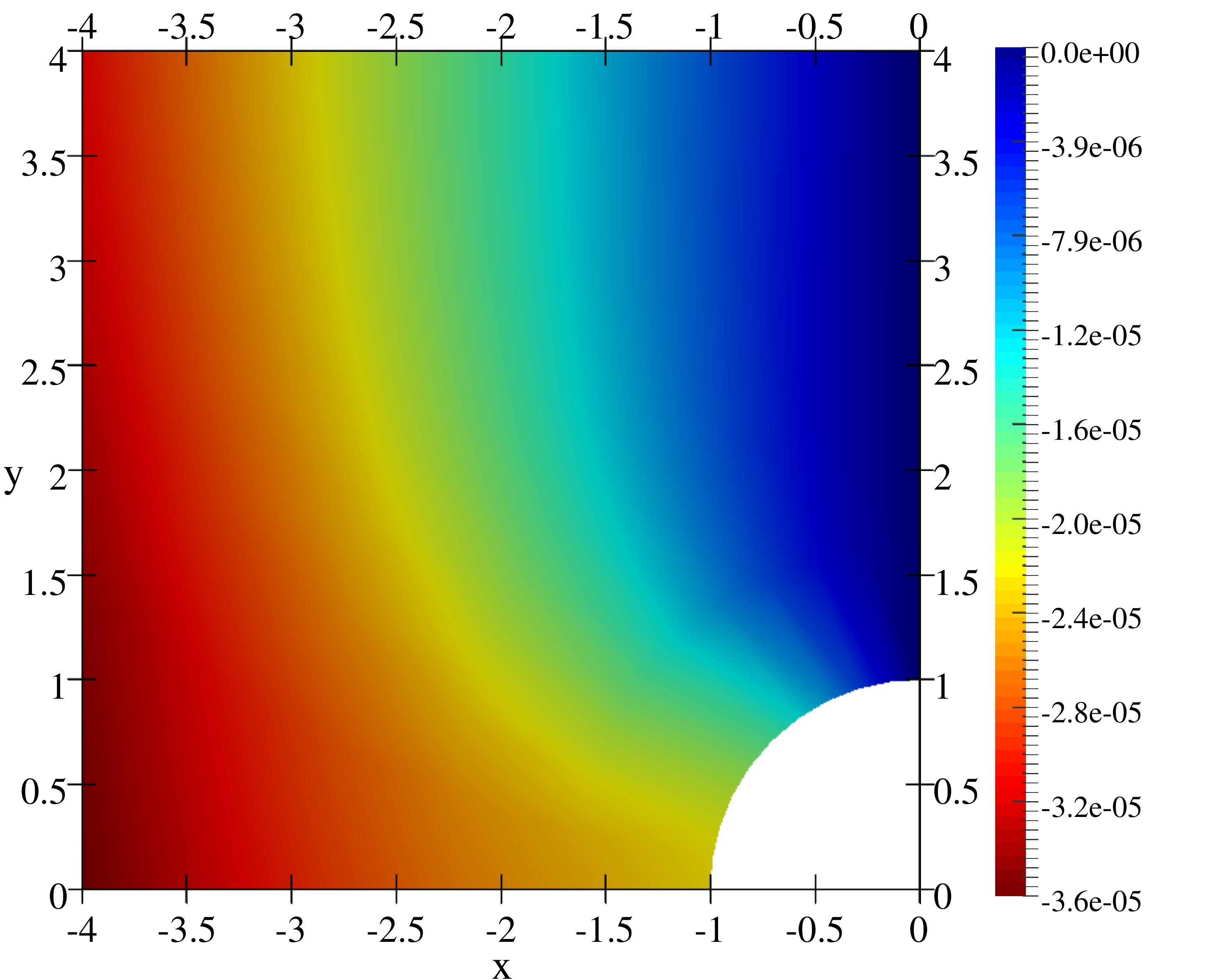}
\captionsetup{width=8.48cm}
\caption{$u_x^{exact}$}
\label{u_x_numerical_plate_with_hole_nu_04999_exact}
\end{subfigure}
\begin{subfigure}{0.49\columnwidth}
\centering
\includegraphics[trim={0cm 0cm 0cm 0cm},clip,width=0.8\columnwidth]{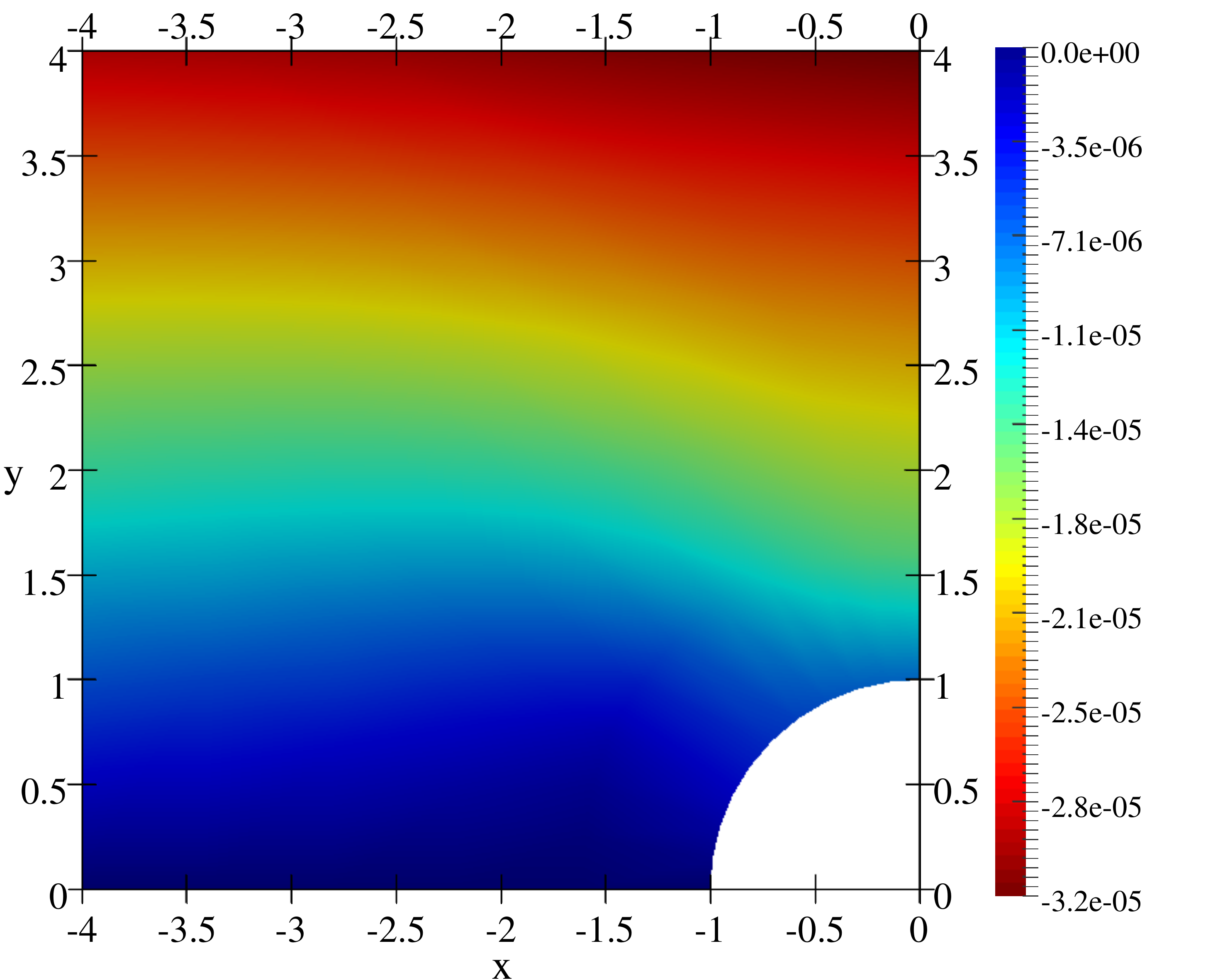}
\captionsetup{width=8.48cm}
\caption{$u_y^{exact}$}
\label{u_y_numerical_plate_with_hole_nu_04999_exact}
\end{subfigure} 
\caption{(a-b)Deformed geometric description of an infinite plate with a hole problem (a quarter portion) for 10$\times$5 NURBS elements with quadratic basis along $\xi$ and $\eta$ direction (magnified by the factor of $1\times10^4$ for better visualization), and contour plots for horizontal and vertical displacement for the mesh illustrated in Figure~\ref{Control_point_mesh_Plate_with_a_hole_problem}-\ref{Element_discretization_Plate_with_a_hole_problem} using hybrid IGA formulation (c-d) alongside the analytical solution (e-f) for $\nu=0.4999$}\label{plate_with_hole_nu_04999_combined_solutions}
\end{figure}

The minimal requirement to maintain the exact geometry is to use the quadratic NURBS basis functions along the $\phi$ direction. Furthermore, the coarsest possible mesh to represent the exact geometry consists of two quadratic elements, one in the radial direction and two in the $\phi$ direction. The CP and the respective weights are given in the appendix (Table~\ref{control_points_plate_with_hole_quad}-\ref{control_points_plate_with_hole_cub}). Once the initial mesh is generated, the sequence of meshes is constructed using the $h$-refinement. One such mesh of 10$\times$5 elements is illustrated in Figure~\ref{Control_point_mesh_Plate_with_a_hole_problem}-\ref{Element_discretization_Plate_with_a_hole_problem}, which focuses on the number of control points involved, the control point mesh and the respective element discretization of a domain. For a better understanding, the comprehensive results for the mesh described in Figure~\ref{Control_point_mesh_Plate_with_a_hole_problem}-\ref{Element_discretization_Plate_with_a_hole_problem} are presented in Figure~\ref{comprehensive_results_plate_with_hole_03}-\ref{plate_with_hole_nu_04999_combined_solutions} which focuses on the deformed configuration of the problem domain along with the numerical displacement field  in comparison with the analytical solution for different Poisson's ratios.

The problem is solved using the conventional and hybrid IGA alongside their FE counterparts. To test the efficiency and robustness of the method, the problem is studied for the two cases. For the first case, $\nu$ is considered as 0.3 so that the solution is unaffected by volumetric locking. Secondly, $\nu$ value is considered as 0.4999 to analyze the nearly incompressible behavior which is highly influenced by locking. 

\begin{figure}[pos=hb]
\begin{subfigure}{0.49\columnwidth}
\centering
\includegraphics[trim={0cm 0cm 0cm 0cm},clip,width=1\columnwidth]{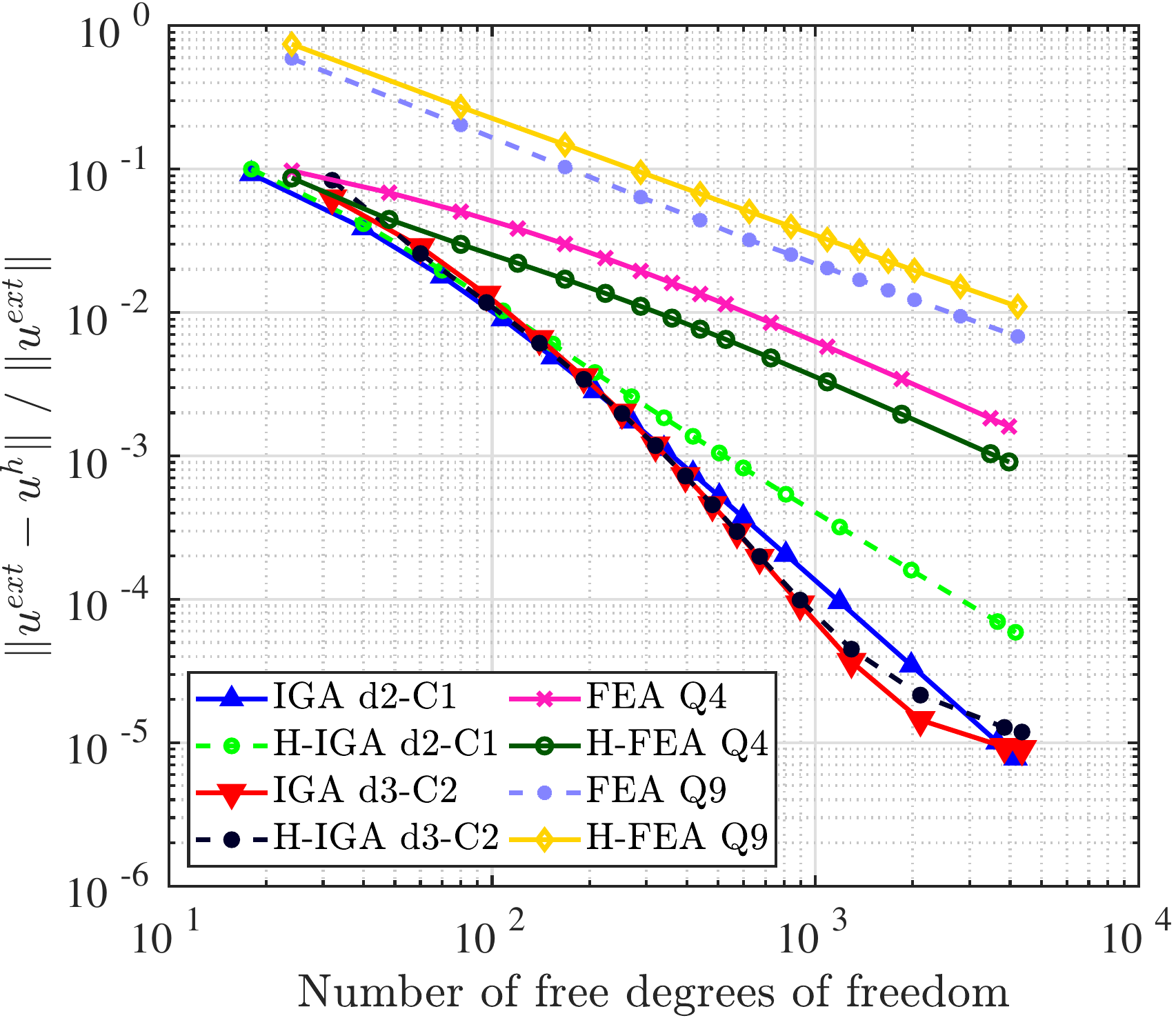}
\captionsetup{width=8.48cm}
\caption{$\nu = 0.3$}
\label{L_2_error_plate_with_hole_03}
\end{subfigure}
\begin{subfigure}{0.49\columnwidth}
\centering
\includegraphics[trim={0cm 0cm 0cm 0cm},clip,width=1\columnwidth]{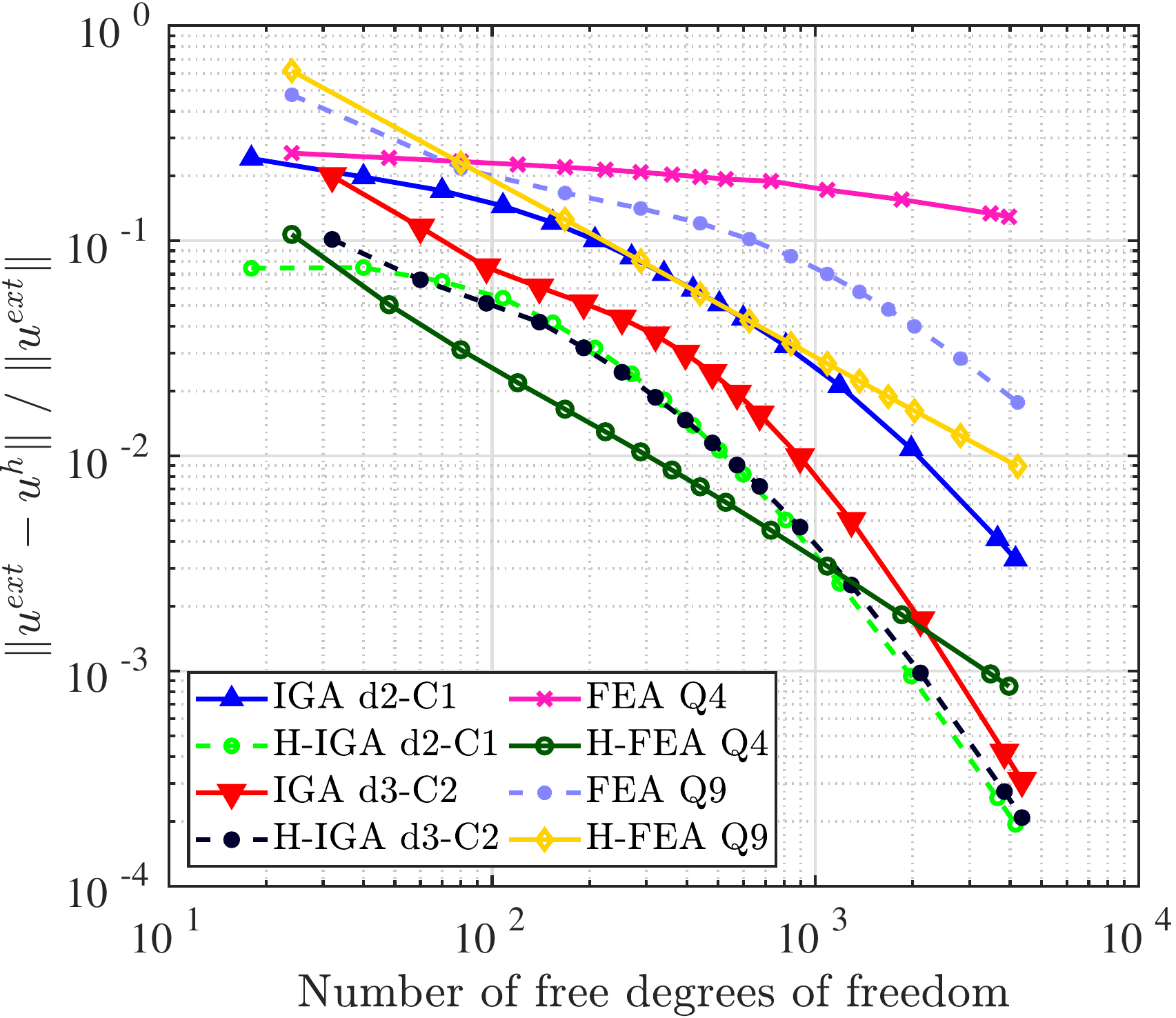}
\captionsetup{width=8.48cm}
\caption{$\nu = 0.4999$}
\label{L_2_error_plate_with_hole_04999}
\end{subfigure} 
\caption{Convergence study of a relative $L_2$ error norm of displacement versus the active degrees of freedom for a plate with hole problem } \label{convergence_l_2_plate_with_hole_combine}
\end{figure}

For the first case, the convergence solution for the relative $L_2$ error norm of displacement versus the number of active degrees of freedom is evaluated as shown in Figure~\ref{L_2_error_plate_with_hole_03}. It can be seen that the convergence rates for the conventional and hybrid IGA are closely identical, whether it is for quadratic or cubic basis functions. This authenticates the fact that the hybrid IGA formulation is not restricted to locking dominated problem domains but can also be effectively used for problems that are not influenced by the locking. Moreover, with a comparative perspective, the IGA or proposed hybrid IGA performs superior to its FEA or hybrid FEA counterpart.

In the second case, the nearly incompressible behavior is investigated by setting the Poisson's ratio close to the value $\frac{1}{2}$. The convergence curves specific to $\nu = 0.4999$ are presented in the Figure~\ref{L_2_error_plate_with_hole_04999}. It can be seen that the conventional IGA locks while using the quadratic basis functions. Though refinement considerably reduces the error, the proposed hybrid IGA outperforms the conventional IGA formulation in terms of coarse mesh accuracy. Similar results are obtained with cubic basis functions, conventional IGA is less sensitive to locking, but the hybrid IGA provides better accuracy at a relatively low number of elements. Furthermore, the proposed IGA is significantly better than the FE counterparts. Either its Q4 or Q9 element, the conventional FEA locks severely; however, the hybrid FE considerably alleviates the locking yet the results are inferior to IGA.

\section{Conclusion}
In the present work, a novel class of hybrid elements is proposed to alleviate the locking in NURBS-based IGA using a two-field Hellinger-Reissner variational principle. The proposed elements exercise the independent interpolation schemes for displacement and stress field. The displacement field is approximated using the standard NURBS-based interpolations; however, the special treatment is followed in approximating the stress field to ensure locking-free results. The principle concept follows the assessment of normal stress components in relation to derivatives of the respective displacement interpolations and the choice of higher-order terms in shear components are defined such that they ensure correct stiffness rank and suppress the spurious zero-energy mode, which is necessary to avoid locking. To asses the performance, the proposed formulation along with the conventional single-field IGA and Lagrangian-based FEA and hybrid FEA formulation, have been implemented on several two-dimensional linear elastic examples. The results for typical benchmark numerical examples authenticate the potency of the proposed formulation. The two-field hybrid IGA tends to perform well, for analyzing nearly incompressible problem domains that are severely affected by volumetric locking as well as for thin plate and shell problems where the shear locking is dominant, by alleviating the different types of locking. The same formulation can also be implemented for standard problems (i.e.\, those without the locking effect) without affecting the solution accuracy. The effectiveness of the method can be clearly seen while using the lower order NURBS basis functions. However, higher-order basis functions are less affected by locking, but the proposed method outperforms the results in an aspect of coarse mesh accuracy, which eventually leads to lower computational efforts. In all the presented examples, the results obtained by the hybrid IGA formulation are superior to that of the conventional IGA or FEA formulations. Though the results for higher-order NURBS basis are only marginally better for locking free problem domains, the robustness of the method to perform well in all situations can not be ignored. The extension of the proposed method in a three-dimensional regime and further for a non-linear analysis of the locking-dominated problems would be of interest for future study.

\section*{Declaration of competing interest}

The authors declare that they have no known competing financial interests or personal relationships that could have appeared to influence the work reported in this paper.

\appendix
\section{Geometric data for modeling the base coarse mesh for the presented problems}
The section is intended to provide the required control points co-ordinates and the respective weights to construct the initial mesh for the stated problem domains. The sequence of the refined meshes are modeled by employing the knot insertion or the degree elevation algorithms in one or both direction. For the reader's interest, the MATLAB codes based on these algorithms can be found in an open source code library called NURBS toolbox \cite{NURBS_toolbox}. The subroutines named \textbf{\textit{nrbkntins}} and \textbf{\textit{nrbdegelev}} are of the particular interest for the refinement strategies.
% ===================== table: Cp for rectangular beam ==============================
\begin{table}[width=1\linewidth,cols=7,pos=H]
\caption{Control points for modeling the initial coarse mesh (single element) for a rectangular beam problem with linear basis along $\xi$ and $\eta$ direction (weights associated with each control point is 1, knot vector along $\xi$ and $\eta$ direction is $\left[0,0,1,1 \right]$).} \label{Table_CP_rectangular_beam}
\begin{tabular*}{\tblwidth}{@{} LLL|LL|LL@{} }
\toprule 
\multirow{2}{*}{$j$} & \multicolumn{2}{L|}{Slenderness = 10} & \multicolumn{2}{L|}{Slenderness = 100} & \multicolumn{2}{L}{Slenderness = 1000} \\
\cline{2-7}
 & $P_{1,j}$ & $P_{2,j}$ & $P_{1,j}$ & $P_{2,j}$ & $P_{1,j}$ & $P_{2,j}$ \\
\midrule
1 & (0,0,0) & (100,0,0) & (0,0,0) & (100,0,0) & (0,0,0) & (100,0,0) \\ 
2 & (0,10,0) & (100,10,0) & (0,1,0) & (100,1,0) & (0,0.1,0) & (100,0.1,0)\\
\bottomrule
\end{tabular*}
\end{table}
% -----------------------------------------------------------------------------
%\subsection{Curved beam}
% ===================== table: Cp for rectangular beam ==============================
\begin{table}[width=1\linewidth,cols=10,pos=h]
\caption{Control points for modeling a single element representing the problem domain of a curved beam with quadratic basis along $\xi$ and linear basis along $\eta$ direction, the respective knot vectors are $\boldsymbol\Xi = \left[ 0,0,0,1,1,1 \right]$ and $\boldsymbol H = \left[0,0,1,1 \right]$.} \label{control_points_curved_beam}
\begin{tabular*}{\tblwidth}{@{} LLLL|LLL|LLL@{} }
\toprule 
\multirow{2}{*}{$j$} & \multicolumn{3}{L|}{Slenderness = 10} & \multicolumn{3}{L|}{Slenderness = 100} & \multicolumn{3}{L}{Slenderness = 1000} \\
\cline{2-10}
 & $P_{1,j}$ & $P_{2,j}$ & $P_{3,j}$ & $P_{1,j}$ & $P_{2,j}$ & $P_{3,j}$ & $P_{1,j}$ & $P_{2,j}$ & $P_{3,j}$ \\
\midrule
1 & (0,9.5,0) & (9.5,9.5,0) & (9.5,0,0) & (0,9.95,0) & (9.95,9.95,0) & (9.95,0,0) & (0,9.995,0) & (9.995,9.995,0) & (9.995,0,0)\\
2 & (0,10.5,0) & (10.5,10.5,0) & (10.5,0,0) & (0,10.05,0) & (10.05,10.05,0) & (10.05,0.1,0) & (0,10.005,0) & (10.005,10.005,0) & (10.005,0,0)\\
\bottomrule
\end{tabular*}
\end{table}

\begin{table}[pos=h,width=0.45\textwidth]
\begin{minipage}{0.45\textwidth}
\captionsetup{width=0.45\textwidth}
\caption{Weights associated with the control points for all slenderness ratios of curved beam problem} \label{weights_curved_beam}
\begin{tabular*}{1\textwidth}{@{} LLLL@{} }
\toprule 
$j$ & $w_{1,j}$ & $w_{2,j}$ & $w_{3,j}$ \\
\midrule
1 & 1 & 0.7071 & 1 \\
2 & 1 & 0.7071 & 1 \\
\bottomrule
\end{tabular*}
\end{minipage}
\hfill
\begin{minipage}{0.45\textwidth}
\centering
\caption{Control points for modeling a single element representing the domain of a Cook's membrane problem with linear basis along $\xi$ and $\eta$ direction (weights associated with all the control points = 1 and knot vector along $\xi$ and $\eta$ direction is $\left[0,0,1,1 \right]$).} \label{CP_and_weights_for_cooks_membrane_problem}
\begin{tabular*}{1\textwidth}{@{} LLL@{} }
\toprule 
$j$ & $P_{1,j}$ & $P_{2,j}$ \\ 
1 & (0,0,0) & (48,44,0) \\
2 & (0,44,0) & (48,60,0) \\
\bottomrule
\end{tabular*}
\end{minipage}
\end{table}

\begin{table}[width=1\linewidth,cols=9,pos=H]
\caption{Control points and respective weights for modeling a quarter portion of a plate with a hole problem using two quadratic elements (two along $\xi$ direction and one along $\eta$ direction, knot vectors along $\xi$ and $\eta$ are given as; $\boldsymbol\Xi = \left[ 0,  0,  0.5,  1,  1 \right]$, $\boldsymbol H = \left[ 0, 0, 1, 1 \right]$, Radius $= r = 1$, Length $= L = 4$, weight $=w^{cp}=0.5(1 + 1/\sqrt 2)$} 
\label{control_points_plate_with_hole_quad}
\begin{tabular*}{\tblwidth}{@{} LLLLLLLLL@{} }
\toprule 
$j$ & $P_{1,j}$ & $P_{2,j}$ & $P_{3,j}$ & $P_{4,j}$ & $w_{1,j}$ & $w_{2,j}$ & $w_{3,j}$ & $w_{4,j}$\\
\midrule
1 & ($-r$,0,0) & ($-r$,0.414,0) & ($-0.414$,$r$,0) & (0,1,0)& 1 &$w^{cp}$ & $w^{cp}$& 1\\
2 & (-2.5,0,0) & (-2.5,0.75,0) & (-0.75,2.5,0) & (0,2.5,0) & 1 & 1& 1& 1\\
3 & ($-L$,0,0) & ($-L$,$L$,0) & ($-L$,$L$,0) & (0,$L$,0) & 1 & 1& 1& 1\\
\bottomrule
\end{tabular*}
\end{table}

\begin{table}[width=1\linewidth,cols=11,pos=h]
\caption{Control points and respective weights for modeling a quarter portion of a plate with a hole problem using two cubic elements (two along $\xi$ direction and one along $\eta$ direction, knot vectors along $\xi$ and $\eta$ are given as; $\boldsymbol\Xi = \left[ 0,  0, 0, 0.5,  1, 1, 1 \right]$, $\boldsymbol H = \left[ 0, 0,0, 1, 1,1 \right]$} 
\label{control_points_plate_with_hole_cub}
\begin{tabular*}{\tblwidth}{@{} LLLLLLLLLLL@{} }
\toprule 
$j$ & $P_{1,j}$ & $P_{2,j}$ & $P_{3,j}$ & $P_{4,j}$ & $P_{5,j}$ & $w_{1,j}$ & $w_{2,j}$ & $w_{3,j}$ & $w_{4,j}$ & $w_{5,j}$\\
\midrule
1 & (-1,0,0) & (-1,0.2612,0) & (-0.7929,0.7929,0) & (-0.2612,1,0) & (0,1,0) & 1 & 0.9024 & 0.8047 & 0.9024 & 1 \\
2 & (-2,0,0) & (-2.0696,1.5942,0) & (-2.0219,2.0219,0) & (-1.5942,2.0696,0) & (0,2,0) & 1 & 0.9349 & 0.8698 & 0.9349 & 1 \\
3 & (-3,0,0) & (-3.0673,2.8376,0) & (-3.0798,3.0798,0) & (-2.8376,3.0673,0) & (0,3,0) & 1 & 0.9675 & 0.9349 & 0.9675 & 1 \\
4 & (-4,0,0) & (-4,4,0) & (-4,4,0) & (-4,4,0) & (0,4,0) & 1 & 1 & 1 & 1 & 1 \\
\bottomrule
\end{tabular*}
\end{table}
% -----------------------------------------------------------------------------

%\bibliographystyle{model1a-num-names}
\bibliographystyle{cas-model2-names}%\bibliographystyle{model1-num-names}
\bibliography{library_edited}% Loading bibliography database

\end{document}